\documentclass[11pt]{amsart}

\usepackage[T1]{fontenc}
\usepackage[protrusion=false]{microtype} % Micro-typography for better spacing
\usepackage{mathrsfs}   % For \mathscr

\usepackage{amsmath, amssymb, amsfonts, amsthm}
\usepackage{mathtools} % Extension of amsmath, fixes bugs and adds features
\mathtoolsset{showonlyrefs} %remove unreferenced equation numbers
\usepackage{xparse}    % For advanced command definition
\usepackage{bm}

\usepackage[a4paper, top=1.8cm, bottom=1.8cm, left=1.8cm, right=1.8cm]{geometry}
\usepackage{paralist}  % Inline lists
\usepackage[dvipsnames, table]{xcolor}
\usepackage{caption}
\usepackage{subcaption}
\usepackage{tikz}
\usetikzlibrary{shapes, arrows.meta, positioning, calc}
\tikzset{
    Block/.style = {rectangle, rounded corners, minimum width=3cm, minimum height=1cm, text centered, font=\normalsize, color=black, draw=black, line width=1pt, fill=white},
    DA/.style    = {thick, draw=black, line width=1pt, ->, >=stealth}
}

\usepackage[backend=biber, style=alphabetic, giveninits, backref=true]{biblatex}
\usepackage{thmtools}
\usepackage{thm-restate}

\declaretheorem[numberwithin=section, style=plain]{theorem}
\declaretheorem[numberwithin=section, numberlike=theorem, style=plain]{proposition, lemma, corollary, conjecture}
\declaretheorem[numberwithin=section, numberlike=theorem, style=definition]{definition, remark}

\usepackage{url}
\usepackage{hyperref}
\hypersetup{
    colorlinks=true,
    linktoc=all,
    linkcolor=blue,
    citecolor=green,
    urlcolor=blue
}

\numberwithin{equation}{section}
\newcommand{\mbb}{\mathbb}
\newcommand{\mc}{\mathcal}
\newcommand{\mf}{\mathbf}
\newcommand{\mk}{\mathfrak}
\newcommand{\mrm}{\mathrm} % Added definition (used in your macros but missing in original)

\newcommand{\psmb}{\left( \begin{smallmatrix}}
\newcommand{\psme}{ \end{smallmatrix} \right)}
\newcommand{\smat}[4]{\left( \begin{smallmatrix} #1 & #2 \\ #3 & #4 \\ \end{smallmatrix} \right)}

\newcommand{\tp}[1]{#1^t}

\newcommand{\lan}{\langle}
\newcommand{\ran}{\rangle}

\DeclareMathOperator{\tr}{tr}

\renewcommand*{\mod}{\operatorname{mod}}

\newcommand*{\GL}[2]{\operatorname{GL}_{#1}(#2)}

\newcommand*{\Sp}[2]{\operatorname{Sp}_{#1}(#2)}
\newcommand*{\GSp}[2]{\operatorname{GSp}^{+}_{#1}(#2)}

\newcommand{\z}{\mathbb{Z}}
\newcommand{\Q}{\mathbb{Q}}

\newcommand{\complex}{\mathbb{C}}

\newcommand{\h}{\mbb{H}}

\newcommand{\hn}{\mbb H_n}

\newcommand{\sptwo}{\mrm{Sp}_2( \mf Z)}

\newcommand{\glnz}{\mrm{GL}_n( \mbb Z)}
\newcommand{\glnc}{\mrm{GL}_n(\mbb C)}

\newcommand{\slnz}{\mrm{SL}_n(\mbb Z)}
\newcommand{\spnz}{\mrm{Sp}_n(\mbb Z)}

\newcommand{\q}{\quad}
\newcommand{\qq}{\qquad}

\newcommand{\sumn}{\sum \nolimits}

\newcommand{\dia}{\mrm{diag}}

\newcommand{\Mrho}[2]{M_\rho(\Gamma_\pm^{(#1)}(#2, L)}

\author{Pramath Anamby }
\address{School of Arts and Sciences\\ 
Ahmedabad University\\ 
Ahmedabad -- 380052, India.}
\email{pramath.anamby@gmail.com, pramath.anamby@ahduni.edu.in}

\author{Soumya Das}
\address{Department of Mathematics\\ 
Indian Institute of Science\\ 
Bengaluru -- 560012, India.}
\email{soumya@iisc.ac.in}

\date{}
\subjclass[2020]{Primary 11F30, 11F46,  Secondary 11F50, 11F37} 
\keywords{Fourier coefficients, Siegel modular forms, fundamental discriminant, vector
valued, non-vanishing, Jacobi forms}
\begin{document}
\title[Fundamental Fourier coefficients]{Fundamental Fourier coefficients of Siegel modular forms of higher degrees and levels}

\begin{abstract}
    We prove the following statement about any Siegel modular form $F$ of degree $n$ and arbitrary odd level $N$ on the group $\Gamma_{0}^{(n)}(N)$. Let $A(F,T)$ denote the Fourier coefficients of $F$ and write $T=(T(i,j))$. 
    
    Suppose that $F$ has  a non-zero Fourier coefficient $A(F,T_0)$ such that $(T_0(n,n),N)=1$. Then there exist infinitely many odd and square-free (and thus fundamental) integers $m$ such that $m=\mrm{discriminant}(T)$ and $A(F,T)\neq 0$. In the case of odd degrees, we prove a stronger result by replacing odd and square-free with odd and prime. We also prove quantitative results in this direction.

    As a consequence, we can show in particular that the statement of the main result in \cite{eischen2024algebraicity} about the algebraicity of certain critical values (and the expected functional equation) of the spinor $L$-functions of holomorphic newforms (in the ambit of Deligne's conjectures) on congruence subgroups of $\mrm{GSp}(3)$ is unconditional.
\end{abstract}
\maketitle

\section{Introduction}

This paper is a compendium of results about distinguishing a Siegel modular form by its arithmetically interesting Fourier coefficients, especially in the level aspect. Let $A(F,T)$ denote the Fourier coefficients of a Siegel modular form on some congruence subgroup (see Section \ref{prelim}). We say that a Fourier coefficient $A(F,T)$ of  $F$ has `property $P$' if the indexing matrix $T$ has the property $P$. For instance, we say that $A(F,T)$ is  `fundamental' (resp. `odd and square-free, resp. prime') if $T$ is such that $\mrm{discriminant}(T)$ is fundamental (resp. odd and square-free, resp. a prime). We make this convention throughout this paper (see \eqref{disc}).

When the property $P$ is arithmetic, such results are known to be of significance in the analytic theory of $L$-functions, and also in their arithmetic properties. This connection is the main motivation behind this article. For instance, in one of the first non-trivial results in this direction, Saha \cite{saha2013siegel} showed the existence of non-zero `fundamental' Fourier coefficients of any non-zero Siegel cusp form of level one -- this implied the existence of good Bessel models for $\GSp{2}{\z}$, which has lots of applications to special values of $L$-functions originating in the works of Furusawa and also played an important role in the proof of the transfer of automorphic representations from $\mrm{GSp}_2$ to $\mrm{GL}_4$. 
Similarly, in degree $3$ and level one, `prime' Fourier coefficients imply, upon using the integral representation for the spinor $L$-function $Z_F(s)$ of an eigenform $F$ obtained by Pollack \cite{pollack2017spin}, an unconditional proof of its functional equation. This was proved in \cite{boche-das}, and it was remarked that the argument therein might work for higher levels, ``{\it but we expect more complicated answers (compare \cite{anamby2024non} for $n = 1$). One may have to take into account the Fourier expansions at all cusps simultaneously, and one may expect new difficulties concerning primes dividing the level}.'' Our work (and \cite{washburn2025certain}) shows the cusp at $\infty$ is enough. In some sense, this paper may be viewed as a sequel to \cite{boche-das}.

Apart from its own intrinsic interest (cf. the above remark), one of the main reasons for our interest in this topic is the recent result by Eischen et al \cite{eischen2024algebraicity}. In that paper, the algebraicity of the special values of the spinor $L$-function $Z_F(s)$ of a newform $F$  on the group $\Gamma^{(3),0}(N)$, of scalar weight $2k \ge 12$ (upto periods, in the sense of Deligne's conjectures), were obtained at the subset  $2 \le r \le k-2$ of the set of critical points $\{3-k,\ldots, k-2\}$. This result was conditional (see  \cite[Condition~1.1.1]{eischen2024algebraicity}) on the existence of a non-zero Fourier coefficient of $F$ corresponding to a maximal order in a quaternion algebra over $\Q$. Prior to this paper, Washburn \cite{washburn2025certain} recently made notable progress on this topic; however, his results are insufficient to settle the above issue.  We will see below that the results of this paper make \cite{eischen2024algebraicity} unconditional when specialized to scalar weights and newforms.

\subsection{Main results}\label{sec:mainresults}
Let us now describe the main results of the paper. Even though our goal is to prove results on the Hecke congruence subgroups $\Gamma_0^{(n)}(N)$, we first define (and work with) a congruence subgroup which sits between $ \Gamma_1^{(n)}(N)$ and $\Gamma_0^{(n)}(N)$ (see section~\ref{prelim-smf}). This is for serious technical reasons that will be discussed later.

We put $\displaystyle \Gamma_1^{(n)}(N, L):= \{\smat{A}{B}{C}{D}\in\mrm{Sp}_n(\mbb Z): C\equiv 0\bmod N, A\equiv 1_n\bmod L  \}$,  where $N$ is any integer and $L$ is its square-free part.
Let $\mc E \subset\sptwo$ be the group of order $2^n$ consisting of all matrices of the form $\displaystyle \smat{A}{0}{0}{A^{-1}}$ such that $A $ is a diagonal matrix with diagonal entries $\pm 1$.
We define the group $\Gamma_\pm^{(n)}(N,L)$ by the internal semi-direct product 
\begin{align} \label{gamma-pm}
    \Gamma_\pm^{(n)}(N,L):= \mc E \ltimes \Gamma_1^{(n)}(N,L),
\end{align}
where $\mc E$ acts on $\Gamma_1^{(n)}(N,L)$ by conjugation. It is clear that every element of $\Gamma_\pm^{(n)}(N,L)$ can be uniquely written as $\varepsilon \cdot M$ where $\varepsilon \in \mc E$ and $M \in \Gamma_1^{(n)}(N,L)$. Note that $\mc E \cap \Gamma_1^{(n)}(N,L)=\{I\}$ and
\begin{align}
     \Gamma_1^{(n)}(N)\subset\Gamma_{\pm}^{(n)}(N,L)\subset \Gamma_0^{(n)}(N). 
\end{align}

Let $\rho$ be a polynomial, not necessarily irreducible,  representation of $ \mrm{GL}_n(\mbb C)$. Denote by $\Mrho{n}{N}$ the vector space of holomorphic vector-valued Siegel modular forms on $\Gamma_1^{(n)}(N)$ with automorphy factor $\rho$ (see Section~\ref{prelim} for more details) and `determinantal' weight $k(\rho) $ with $k(\rho) - \frac{n-1}{2} \ge 2 $ (see Section~\ref{prelim} for the definition). We need one more piece of notation. Let $\Lambda_n$ denote the set of all $n\times n$ symmetric, positive semi-definite, half-integral matrices.
For $M \in \Lambda_n$, denote by $d(M)$ its `absolute discriminant' (i.e., ignoring the usual sign), where the discriminant $\mrm{disc}(2M)$ of $2M$ is defined by 
\begin{equation} \label{disc}
 \mrm{disc}(2M)= \begin{cases}    (-1)^{n/2}\det(2M) &\text{ if } n \text{ is even,} \\
(-1)^{(n-1)/2}\frac{1}{2} \det(2M) &\text{ if } n \text{ is odd}; \end{cases} \text{
and we put 
 \, } d(M):=|\mrm{disc}(2M)|.
\end{equation}
Its well known that $\mrm{disc}(2M) \equiv 0,1 \mod 4$. Thus if $d(M)$ is odd and square-free, then $\mrm{disc}(2M)$ being odd must be a fundamental discriminant. By a slight abuse of notation, we call $M$ to be `fundamental' if $d(M)$ is a fundamental discriminant.

Further, put
\begin{align} \label{sfx}
{\mk S}_F (X) := \{ d \leq X, \,  d \, \mbox{\rm{odd, square-free}}, \, (d,N)=1 \, \mid  d(T)=d \, \mbox{\rm{ for some }} T \mbox{\rm{ and }}  A(F,T) \neq 0 \}.  
\end{align}

\begin{theorem}\label{mainthm}
Let $N>1$ be an odd integer. Let $F \in \Mrho{n}{N}$ be non-zero and $k(\rho) - \frac{n-1}{2} \ge 2 $. Suppose that there exists a $T_0\in \Lambda_n$ with $(T_0(n,n), N)=1$ such that $A(F,T_0)\neq 0$. Then, the following  quantitative result holds: for any given $\epsilon>0$,
\begin{align*}
\# {\mk S}_F (X)  \gg\begin{cases} X(\log X)^{-1/2} &\text{ if } n \text{ is odd };  \\ 
X^{5/8 - \epsilon} &\text{ if } n \text{ is even }.
\end{cases}
\end{align*}
where the implied constant depends only on $F$ and $\epsilon$.
\end{theorem}

In fact, we establish a stronger result. Theorem~\ref{mainthm} follows from the more general Theorem~\ref{diagonal-sqrfree}, which additionally demonstrates that $T$ can be assumed to be diagonal modulo $L$.

Qualitatively, the theorem says that for any $F \in \Mrho{n}{N}$ with the same condition on the weights as above, there are infinitely many $T\in \Lambda_n$ with $(d(T), N)=1$ such that $A(F,T)\neq 0$ as $d(T) \to \infty$ along the set of fundamental discriminants.
The assumption $N>1$ is not a restriction, as the result is known for $N=1$ from \cite{boche-das}.

{\bf Throughout this paper we take $N>1$}.

We remark here that the $T_0$ in the hypothesis of Theorem \ref{mainthm} need not be positive-definite (see Remark \ref{rmk:Tpsd}). But the $T$ in the conclusion of the Theorem will be positive-definite by virtue of the condition $(d(T), N)=1$ and $N>1$.

\begin{remark}\label{rmk:Tpsd}
By not assuming $T_0$ to be positive in the hypothesis of Theorem \ref{mainthm} ( and Theorem \ref{diagonal-sqrfree} later), one reduces the extra work of checking for positive-definiteness at each step of the induction process. This helps us in dealing with cusp and non-cusp forms simultaneously. Also see Lemma \ref{lem:psd-pd} and Remark \ref{rmk:psd-pd}.
\end{remark}
\subsection{Some difficulties in higher degrees and levels}
The proof of Theorem~\ref{mainthm} requires a very careful selection of induction hypotheses and the also the associated congruence subgroup! Since this point is the crux of the matter in our approach, with degrees $n \ge 3$ being the main points of interest, we wish to discuss it in some detail (see also the introduction in \cite{washburn2025certain}). 
\subsubsection{Vector-valued setting} \label{vec-smf}
First, working in the setup of vector-valued modular forms is indispensable in our approach -- this was of course the case in \cite{boche-das} as well, which we follow in many aspects. Even if one wishes \textit{only} a result for scalar-valued modular forms, our induction-based approach makes it necessary to pass through vector-valued objects -- hence, starting from vector-valued objects is a necessity. Of course, treatment of these objects is much desirable in their own right, and has applications (see, e.g., \cite{marzec2018non}).

\subsubsection{Usage of the Fourier-Jacobi(FJ) model}
A large part of the analytic theory of Siegel modular forms in higher degrees has seen development via the so-called Fourier-Jacobi models. In the present context, it was used in  \cite{yamana2009determination},  \cite{saha2012determination}, \cite{saha2013siegel}, \cite{anamby2019distinguishing} and also in \cite{boche-das}. The main idea behind the success of the approach in \cite{boche-das} with respect to the FJ-model -- which we follow in this paper as well -- is that we work with the so-called (non-zero) `primitive' theta components $h_\mu$ of suitable FJ-coefficients $\phi_{F,T}$ (defined by $\displaystyle F(Z)=\sumn_T \phi_{F,T}e(TZ)$ and $\phi_{F,T}=\sumn_\mu h_\mu \Theta_{T,\mu}$, see subsection~\eqref{jacobi}); as opposed to the `Eichler-Zagier map' (essentially $\sumn_\mu h_\mu$, say, as used in \cite{saha2013yoshida}), which is commonly employed to pass to the half-integral weight setting. However, in general, this map may be identically zero, even in degree $2$ and non-square-free levels -- at least it is difficult to show that it is not zero in higher degrees and matrix indices. More recently, in \cite{washburn2025certain}, the problem was considered for the first time for higher levels in degree at least $3$,  and the same issue about the `Eichler-Zagier map' arose. This issue also appeared in \cite{anamby2019distinguishing}, and it requires quite some amount of work to overcome this -- by considering character-twisted Eichler-Zagier maps and controlling the levels, or working away from the level.
\subsubsection{Correct induction setup}
This is perhaps the most serious point in our approach.
Even while working with the vector-valued setup, in the case of higher levels, especially when the degree  $n \ge 2$, there is the issue of using some notion of `new-ness' of the modular form either in the sense of \cite{ibukiyama2012atkin} (which we denote by the property of being ``{\bf content}-new'' or ``{\bf content}-old''), or assuming that $F$ is an eigenfunction of the Hecke operators $U(p)$, $p|N$ (cf. \cite{saha2013yoshida}), or something else along these lines. Otherwise, old forms provide obvious counterexamples to the desired result about fundamental discriminants. The point being: these properties do not usually pass along to the smaller degrees, which is required to invoke the induction hypothesis.
This issue is inevitable -- it was discussed in some detail in \cite{washburn2025certain}, and it was partially overcome by restricting the setup to the case of primitive nebentypus (in which case there is no {\bf content}-oldspace) using the relevant results in \cite{ibukiyama2012atkin}. This setup then behaves well with the induction argument.

For instance, if one wants to start with the `minimal' assumption: say that any $F$ on $\Gamma_0^{(n)}(N)$ with a non-zero `primitive' Fourier coefficient has a fundamental Fourier coefficient, the corresponding induction argument doesn't seem to work in any straightforward manner! (cf. Remark~\ref{der}) 

In the context of applications, even if one could prove that the newforms coming from local representation theory, say as considered in \cite{eischen2024algebraicity}, are newforms in the sense of \cite{ibukiyama2012atkin}, we see no immediate way to attack the problem just from this information. Usually, the only input one uses about the newforms is that they have a non-zero primitive Fourier coefficient, which, as explained in the previous paragraph, is unlikely to give workable inroads.

One faces the same issue if the initial hypothesis was the existence of at least one non-zero Fourier coefficient whose content or the discriminant is away from $N$. However, after widening the induction setting, in some special cases one can indeed prove statements similar to the above-mentioned statements; see Corollary~\ref{cor:Nchi}.

\subsubsection{Choosing the correct setting of levels} \label{intro-lev}
Quite surprisingly, it seems to us that there is not much hope via our method if one starts with, say, odd and square-free levels, as was the case in e.g., \cite{saha2013yoshida} or \cite{washburn2025certain}. In degree $2$ (\cite{saha2013yoshida}), this was successful; however, the final fundamental discriminants that one obtained may not have been co-prime to the level. In higher degrees (\cite{washburn2025certain}), this was also successful, but the primitive nebentypus of $F$ was a somewhat strict condition. For instance, this condition doesn't render the result in \cite{eischen2024algebraicity} as unconditional.
For the result in \cite{eischen2024algebraicity} to become unconditional, one requires this result for the trivial nebentypus (albeit for the groups $\Gamma^{(n),0}(N)$).

There are two reasons why working with a specific type of level or nebentypus is not good for us:

(i) When we take congruences, the level inflates, and the congruence subgroup might change.

(ii) Unless the discriminants involved are co-prime to the level, the central technical machinery of producing `primitive' non-zero theta components of the FJ coefficients breaks down (see  Proposition~\ref{prop:thetacom}).

\subsubsection{Working with the correct congruence subgroup} \label{sbgp}
Tied to the previous point is the choice of which congruence subgroup of level $N$ to start with. The reason is the following.
At various points in the proof, we employ congruences $\mod N$ to pass from a condition on a specific Fourier coefficient to the same condition on the full support of the Fourier expansion. This is crucial for us -- say while passing onto degrees $n-1$ via the Taylor-series based construction $F \mapsto F_{(n-1)}$ (Section~\ref{sec:Fcirc}).
Now, if one attempts to proceed by induction from the setup of $M_\rho(\Gamma_0^{(n)}(N))$ by our approach, then it doesn't work as these types of subgroups are not preserved by the congruence procedure. It turns out that the groups $\Gamma_1^{(n)}(N)$ work well with respect to congruences (see section~\ref{prelim-smf} for the definition). It should not be confused with another related, bigger congruence subgroup $\Gamma_{1, \det}^{(n)}(N)$, which is not suitable for us to start with. 

However, the subgroups $\Gamma_1^{(n)}(N)$ do not work for us, there are two reasons behind this:

(i) In the definition of the congruence subgroups, if $A \equiv 1 \mod N$, then $N$ becomes the `lattice'-level of  the associated FJ-coefficient $\phi_{F,T}$. Since in our approach $N$ can't be assumed to be square-free -- the technical heart of our argument involving Gauss sums  and the ``maximal-rank arguments'' break down (see (ii) below). It turns out that this can be overcome by considering the groups $\Gamma_1^{(n)}(N,L)$, with $L$ being the square-free part of $N$ -- but they are still not enough.

(ii) For the arguments with Gauss sums to work, we additionally need a ``symmetry'' property of the $\phi=\phi_{F,T}$: $c_\phi(n,r)= \pm c_\phi(n,-r)$.
The reason can be explained in simple terms as follows: we encounter matrices whose entries are Gauss sums $G(r)$ essentially in the variables $r$, and $G(r)=G(-r)$ will mean that the matrices will never have maximal rank. We have to choose only one of $\pm r$, and that the symmetry property precisely ensures this. Unfortunately the groups $\Gamma_1^{(n)}(N,L)$ do not endow $\phi$ with this property (clearly $\Gamma_0^{(n)}(N)$ does, but we can't use them). We overcome this by adjoining the `sign' matrices $\mc E$, and
this gives us the groups $\Gamma_\pm^{(n)}(N,L)$ defined in \eqref{gamma-pm}. In particular, these groups behave well with respect to the induction setup. This choice works well with congruences: $(n,N)=1$ if and only if $(n,L)=1$. Moreover, $L$, which is the ``lattice--level'' of $\phi_{F,T}$  appears as (part of the) modulus of the Gauss sums. That it is square-free is crucial to the ``maximal-rank" arguments to work (see subsection~\ref{sec:thetamain})

{\it Therefore, we need to simultaneously start with \textit{arbitrary} odd levels, specific type of congruence subgroups, and work with suitable conditions away from the level.  We have to pay a price for this, though.} These are discussed in the next two subsections.

\subsubsection{Theta components of vector-valued Jacobi forms}
In order to argue by induction, we need to pass on from vector-valued Fourier-Jacobi coefficients $\phi_{F,T}$ to scalar-valued ones. As will be discussed later, our approach would have been significantly simplified if we could have worked with those of the form $(n-1,1)$ decomposition -- which means working with Jacobi forms that have a scalar index; in fact, odd prime $p$ indices would have been enough. We also need good knowledge of the theta-decomposition for these Jacobi forms. This is known for scalar-valued Jacobi forms (of all types). But unfortunately we can not `extract' a scalar-valued Jacobi form by the strategy as in \cite{boche-das} since in this case one is led to consider automorphy with respect to the matrices $\psmb  CZ+D & 0 \\ 0 & 1\psme$ ($CZ+D \in \GL{n-1}{\complex}$), as $C,D$ vary, can not be assumed to be upper-triangular if $n-1 \ge 2$ -- which is required, via the Lie-Kolchin theorem to  ensure that $\rho (\psmb  CZ+D & 0 \\ 0 & 1\psme) $ is also upper-triangular -- in order to `extract' a scalar-valued Jacobi form out of $\phi_{F,p}$. The same issue rules out all the decompositions $(r,n-r)$ with $r \ge 2$.

We only need to know that $c_{\phi_{F,p}}(S,0) \neq 0$ for some $S$. If a theta-decomposition existed for $\phi_{F,p}$, this would have amounted to showing that the `$0$'-th theta-component $h_0=h_{0}(\phi_{F,p})\neq 0$.  This is known for scalar-valued Jacobi forms, but
unfortunately not for vector-valued Jacobi forms. We have a very recent pre-print  \cite{raum25}, which does relate the vector-valued Jacobi forms with a collection of scalar-valued Jacobi forms of scalar-index by means of covariant differential operators. But as far as we are aware, this differential operator is not explicit, and moreover, it is not immediately clear how to deduce the non-vanishing result we wish -- from this isomorphism. The arguments using the $(n-1,1)$ type Jacobi forms, conditional on the assumption $h_{0}(\phi_{F,p})\neq 0$, will be treated elsewhere. %{\color{red} explain more -- why scalar matrices do not suffice to deal with the n vector components}

Therefore we have to work with the more popular $(1,n-1)$ type Jacobi forms, which lead to matrix-index Jacobi forms, from which we do know how to `extract' scalar-valued ones --- as here one is led to consider automorphy with respect to the matrices $\psmb  cz+d & 0 \\ 0 & 1_{n-1}\psme$, which are indeed upper-triangular. But this leads one to a very subtle question related to the `Eichler-Zagier map' for higher levels. This point was the main focus in the work of \cite{washburn2025certain}, where this was handled by considering twisted maps and restricting to primitive nebentypus. We avoid this by using `primitive'-theta components $h_\mu$ of $\phi_{F,T}$. However, we have to work with Jacobi forms with non-trivial lattice-index, which is quite delicate and is discussed in section \ref{sec:JacobiIntro}. Finally, we overcome all of these by reducing the question to the existence of non-vanishing of $h_{\mu}$ with $\mu \equiv 0 \mod L$ after reduction to scalar-valued forms.
In fact there is still another subtle point regarding the technical congruence conditions in the statements of Theorem~\ref{diagonal-sqrfree} and Proposition~\ref{jacobi-theta-main}, which allow us to overcome the difficulties mentioned above. See Remark~\ref{deg-freedom}.

\subsubsection{Jacobi forms with non-trivial lattices and non-vanishing of their theta-components}\label{sec:JacobiIntro}
As mentioned above, the Jacobi forms that arise in this way are on a proper sub-lattice of the form $M \z \times \z$ of $\z \times \z$, say when $n=1$. Many existing results in the literature need to be adapted to this setting for our method to work. 

First, we observe that the theta decomposition for this class of Jacobi forms differs from that of trivial-lattice index forms (see \eqref{theta-decomp}). The issue is compounded by the property that now the set of indices in the theta decomposition feature {\bf  both} the lattice-level $L$ and the index $T$ of $\phi_{F,T}$, see Proposition~\ref{prop:thetacom}. The condition $(L, d(T))=1$ helps untangle them.

A significant part of the article is dedicated to the study of non-vanishing properties of theta components of these Jacobi forms (see Section \ref{jacobi}). As far as we are aware, these are not available in the literature, certainly not for levels. 
%{\color{red} why couldn't we prove $h_0$ nonzero and went for mod L}

%{\color{red} To some extent, we transfer the questions to Jacobi forms with the standard lattice $\z^n \times \z^n$. This induces the support of the resulting theta expansion on multiples of $N$, which we deal with by taking recourse to the theory index-old Jacobi forms, which may be of independent interest, see Sections~\ref{index1-sec},~\ref{index2-sec}. A significant part of the article is dedicated to the study of non-vanishing properties of theta components of these Jacobi forms (see Section \ref{sec:jacobiTheta}). As far as we are aware, these are not available in the literature, certainly not for levels.

%We prove two results. First, we establish the existence of non-zero 'primitive' theta components of odd, square-free matrix index Jacobi forms. Second, we show that for Jacobi forms of scalar index $p$ ($p$ is a prime such that $(p,N)=1$), lattice $N\z^n \times \z^n$ and higher degrees, the theta component $h_0$ does not vanish. The first is required for patching the Fourier expansions of degree $n-1$ with degree $n$ (see proof of Theorem~\ref{thm:Gamma1} given in Section~\ref{sec:thm-gamma1}). The second is also required for the same reason, but while working $\mod N$ (see proof of Proposition~\ref{diagonal-minor} given in section~\ref{sec:diag-minor}). These complications do not arise if one could work with the relatively bigger groups $\Gamma_{1, \det}^{(n)}(N)$.

\subsection{Intermediary results and discussion of the proof}
\label{preps}

In our approach, we need to invoke induction a few times, with varying induction hypotheses. Each induction step facilitates the next. For instance, we start with a mild condition: $F$ has a non-zero Fourier coefficient $A(F,T_0)$ whose right lower entry is co-prime to $N$. This, via induction (with the correct congruence subgroups, see Section~\ref{sbgp}), implies the same property for all the diagonals of another index $T_1$. We next assume this information as an induction hypothesis and proceed. This is illustrated in Section~\ref{traverse} and the flow-chart therein. It seems very difficult for us to argue in any different manner.

%Towards the proof of Theorem \ref{mainthm}, we first prove two intermediary results, as discussed above.

%\begin{restatable}{proposition}{tnndiagonal}\label{tnn-diagonal}
%Let $F\in {\color{red} M_\rho(\Gamma_1^{(n)}(N))}$ be such that $A(F, T_0)\neq 0$ for some $T_0\in\Lambda_n$ with $(T_0(n,n), N)=1$. Then there exists a $T_1\in \Lambda_n$ with $(T_1(j,j), N)=1$ for all $1\le j\le n$ such that $A(F, T_1)\neq 0$.   
%\end{restatable}

The principal technical result of the paper is the following, from which Theorem \ref{mainthm} follows as a direct consequence.
\begin{restatable}{theorem}{diagonalsqrfree} \label{diagonal-sqrfree}
    Let $F\in \Mrho{n}{N}$ be such that $A(F, T_0)\neq 0$ for some $T_0\in\Lambda_n$ with $(T_0(j,j), L)=1$ for all $1\le j\le n$. Then there exist infinitely many $T\in \Lambda_n$ with the following properties.
    \begin{enumerate}
        \item $T\equiv \mrm{diag}(t_1, t_2,..,t_n)\bmod L$;
        \item $d(T)$ is odd and square-free with $(d(T), L)=1$;
        \item $A(F, T)\neq 0$.
    \end{enumerate} 
\end{restatable}
The structure of the proof of Theorem \ref{mainthm} is illustrated schematically below. In this diagram, each set signifies the existence of at least one $T$ with $A(F,T) \neq 0$.

\begin{figure}[!htbp]
\centering
\begin{tikzpicture}[
    % GLOBAL STYLES
    node distance = 0.7cm and 1.5cm, % Vertical and Horizontal spacing
    Block/.style = {
        rectangle, 
        draw, 
        rounded corners, 
        align=center, 
        minimum height=3.5em, 
        text width=3.2cm, % Fixed width ensures rows align perfectly
        inner sep=5pt,
        font=\small
    },
    BlockA/.style = {
        rectangle, 
        draw, 
        rounded corners, 
        align=center, 
        minimum height=1.5em, 
        %text width=3.2cm, % Fixed width ensures rows align perfectly
        inner sep=5pt,
        font=\small
    },
    Line/.style = {
        ->, 
        >={Stealth[length=3mm]}, 
        thick,
        rounded corners=3pt
    },
    Label/.style = {
        font=\tiny,
        text=black,
        align=center,
        auto
    }
]

    % --- ROW 1 (Left to Right) ---
    \node (LC) [BlockA] {All right lower corners of $T \in \Lambda^+_{n}$ away from $L$};
    \node (AD) [BlockA, right=of LC] {All diagonals away from $L$};
    \node (DiagInd) [Block, below=of AD, text width=4.5cm, xshift=1.5cm] {Induction hypothesis: $T$ is ``nice'' \& diagonal $\mod L$ (assuming for size $\le n-1$)};
    
    %--- Row 2 (Right to Left)
    \node (AM) [Block, left=of DiagInd] {$T$ diagonal $\mod L$ except for $T(1,n), T(n,1)$};
    \node (N1N) [Block, left=of AM]  {$T(1,n), T(n,1) \equiv 0 \mod L$};
    \node (DiagL) [Block, below=of N1N, yshift=-0.5cm] {All $T$ with $T$ diagonal $\mod L$};

     % --- ROW 3 (Left to Right) ---
    \node (AwL) [Block,right=of DiagL, yshift=1cm] {All $T$ with $d(T)$ away from $L$};
    \node (AwD) [Block, below=of AwL, yshift=0.5cm] {All $T$ with $d(T)$ away from $d(S)$ for some ``nice'' $S \in \Lambda^+_{n-1}$};
    \node (Thm13) [Block] at (DiagL -| DiagInd) {Theorem \ref{diagonal-sqrfree}};

    % --- ROW 4 (Right to Left) ---
    \node (FUND) [Block, below=of AwD] {Infinitely many fundamental $T$ away from $L$ (Theorem \ref{mainthm})};

    % --- ARROWS ---

    % Row 1
    \draw [Line] (LC) -- (AD) node[midway, above, Label] {Prop \ref{tnn-diagonal}};
    \draw [Line] (AD.east) -- ++(2,0) |- (DiagInd.east);
    % Row 1 to Row 2
    
    \draw [Line] (DiagInd) -- (AM)  node[midway, above, align=center, text width=1.5cm, font=\tiny] { Thm~\ref{diagonal-sqrfree}: Step 2};
    % Row 2
    \draw [Line] (AM) -- (N1N) 
        node[midway, above, align=center, text width=1.5cm, font=\tiny] {Thm~\ref{diagonal-sqrfree}: Step 3} node[midway,below, font=\tiny] { Prop \ref{jacobi-theta-main}};
    \draw [Line] (N1N) -- (DiagL) node[midway, left, align=center, text width=1.5cm, font=\tiny] {Thm \ref{diagonal-sqrfree}: Step 4};
    % Row 2 to 3
    \draw [Line] (DiagL.east) -- (AwL.west) node[midway, below, yshift=-0.25cm]{\tiny{Prop \ref{alldiag-fund-prop}}} node[midway, sloped, above]{\tiny{Stage I}};

    % Row 3 
    \draw [Line] (DiagL.east) -- (AwD.west)  node[midway, sloped, below]{\tiny{Stage II}};
    \draw [Line] (AwL.east) -- (Thm13.west);
    \draw [Line] (AwD.east) -- (Thm13.west);
    % Row 3 to 4
    \draw [Line] (Thm13) |- (FUND.east);
\end{tikzpicture}
\caption{ {\small Flow of the proof: the property of $T=(T(i,j)$ in each box is different and carries with it the property of the box preceding it. ``Nice''-ness of $M$ refers to $d(M)$ being odd, square-free and away from $L$.} }
\label{fig:flow}
\end{figure}

%The implication of the first arrow in Figure \ref{fig:flow} is achieved by the following proposition. 

\subsection{Traversing through Figure~1 to understand the proof of Theorem~\ref{diagonal-sqrfree}} \label{traverse}

We give a brief summary of how the flowchart works. The first basic principle, adopted at {\it each} Step, is to ensure that a property (P) of a fixed Fourier coefficient $A(F,T_0)$ of a non-zero modular form $F$ holds for {\it all} Fourier coefficients of $F$. This is achieved by congruences $\mod T_0$. We  start with $F$ having  the property the $A(F,T_0)\neq 0$ such  that $(T_0(n,n),L)=1$. This is upgraded to a property of all diagonals of another matrix $T_1$ in the next  Step. We thereafter postulate that such an $F$ should have infinitely many $T$ which are ``odd, square-free'' and ``diagonal  $\mod L$'' and away from $L$ such that $A(F,T)\neq 0$ and proceed to prove it by induction on $n$. We abbreviate ``odd, square-free, away from $L$'' with ``nice''.

Induction hypothesis applied to the modular forms created by the left upper and right lower $n-1$ size blocks, ensures that (again by taking congruences) all Fourier coefficients of $F$ are supported on $T$ whose diagonals are away from $L$ and all but the right upper and left lower corners are $\equiv 0 \mod L$. The goal is to get hold of a $T$ with $d(T)$ away from $L$. To achieve this, i.e., to take care of the above two entries of $T$, we prove a result about Jacobi forms (with non-trivial lattice-index, cf. Proposition~\ref{jacobi-theta-main}) ensuring this. This is also illustrated in Figure~2. This Jacobi form has usual index, say $S$, of size $n-1$.

The next procedure is to ensure that $T$ is such that $d(T)$ is away from $d(S)$ as well. This will ensure that $F$ is supported on $T$ with $(d(T), Ld(S))=1$ -- this is what is required when we reduce the problem to elliptic modular forms. In order to get hold of $T$  such that $d(T)$ is away from $d(S)$, we prove the existence of ``primitive'' theta components of $\phi_{F,S_1}$ where possibly $S_1 \equiv S \mod L$ needs to be considered. Special ``nice''-property of $S$ is also crucial here. The proofs of the above two Steps leading to 
$(d(T), Ld(S))=1$ go via Gauss sum calculations and form the technical heart of the paper.

Once we get the Fourier expansion to be supported on $T$ with $(d(T), Ld(S_1))=1$, we simply use the usual strategy $F \to \phi_{F,S_1} \to h_\mu$ to finish the proof as in \cite{boche-das}.

%We prove two results. 

%(i) First, we establish the existence of non-zero 'primitive' theta components of odd, square-free matrix index Jacobi forms. This is required for patching the Fourier expansions of degree $n-1$ with degree $n$ (see proof of Theorem~\ref{thm:Gamma1} given in Section~\ref{sec:thm-gamma1}). 

%(ii) Second, we show that for  scalar-valued Jacobi forms of matrix index, lattice with level structure $L\z^n \times \z^n$ ($L=$ the square-free part of $N$) and higher degrees, the theta components $h_\mu$ ($\mu \equiv 0 \mod L$) do not vanish.  This is a crucial input in the final part of the proof of Theorem~\ref{mainthm} to obtain, from the hypotheses of loc. cit., a $T$ such that $(d_T,N)=1$ and $a(F,T)\neq 0$. 

%{\color{red} basically briefly write about the flow chart}

By virtue of the main results of this paper we can now more or less completely settle, in simple terms, the question of the existence of `fundamental' Fourier coefficients for higher levels. The following theorem, and other similar results can be derived from our main results, see section~\ref{applns}. A slightly more general version is proved in Theorem~\ref{eq-cond}.

\begin{theorem} \label{genres-intro}
    Let $F \in M_{\rho}(\Gamma_{0}^{(n)}(N), \chi)$ be non--zero,  $N$ being odd and $\chi$ even. Then $F$ has a non-zero `fundamental' Fourier coefficient if and only if it has a non-zero `primitive' Fourier coefficient. 
\end{theorem}

Thus, the existence of \textit{fundamental} Fourier coefficients is reduced to a rather simple criterion of ``primitive'' Fourier coefficients. This improves upon or generalizes all the known results about fundamental Fourier coefficients of Siegel modular forms in the level aspect (e.g., \cite{saha2013yoshida}, \cite{boche-das}, \cite{anamby2024non}, \cite{washburn2025certain}, etc).
As examples, one can consider the case of Eisenstein series in $ M_k(\Gamma_{0}^{(n)}(N), \chi)$. 

\subsection{Non-cuspforms in Theorem~\ref{mainthm}}

We make a few remarks about the non-cuspforms treated in Theorem~\ref{mainthm}. In \cite[Remark~2.6]{saha2013siegel},  it was asked whether its main result is true for non-cuspforms, and it was indicated that the similarity of the growth properties for Fourier coefficients of Siegel-and Klingen–Eisenstein series makes the situation rather delicate. Even more so in higher levels.
Incidentally, such a result in level one was first obtained in \cite[Proposition~7.7]{das-boech-FJ2018}, which was generalized to higher degrees in level one in \cite{boche-das}. Whereas in \cite{das-boech-FJ2018}, this was obtained as an application of an asymptotic formula for the Fourier coefficients of non-cusp forms, with our approach, we can settle this issue without working with an explicit description for the Fourier coefficients of Eisenstein series, cf. the above section. 

%The only difference with the of the level one situation is for the modular forms to have a non-zero ``primitive'' Fourier coefficient; whereas this always holds in level one, in higher levels this doesn't -- and we have it as a necessary condition.

\subsection{Some applications}
We give two applications of our main results. First, when the degree $n$ is odd, we can strengthen Theorem~\ref{mainthm} (as in \cite{boche-das}) to odd and prime discriminants (see Theorem~\ref{oddthm}). This is achieved by the ideas from \cite{boche-das}.
When this is applied for $n=3$ in Corollary~\ref{spl-val}, we can prove that the special value result in \cite{eischen2024algebraicity} is unconditional. To do this, we use the close relationship between newforms on $\Gamma^{(n),0}(N)$ and $\Gamma_0^{(n)}(N)$, see Lemma~\ref{U^0-U_0}. This shows that $F$ is an eigenfunction of $U(p)$ for all $p|N$, which implies the existence of a non-zero primitive Fourier coefficient of $F$, which allows one to invoke Theorem~\ref{genres-intro}. Therefore, we can state that in particular (with notations as in \cite{eischen2024algebraicity}):
\begin{theorem}[\cite{eischen2024algebraicity}, Corollary~\ref{spl-val}] \label{eis-thm}
   Let $M$ be an integer.
Let $\pi$ be a cuspidal automorphic representation associated to a holomorphic cuspidal Siegel
eigenform $\phi$ of scalar weight $2k\ge  12$ and level $\Gamma^{(3),0}(N)$ on $\mrm{GSp}_6$. Let $s_0\in\mbb Z$ be such that $4\le s_0\le k-2$. Then
\begin{align}
    \frac{L^{(M)}(s_0, \pi, \mrm{Spin})}{\pi^{4s_0+6k-6}\lan \phi^\sharp, \phi\ran}\in \mbb Q(\phi).
\end{align}
\end{theorem}
In fact, our result also implies, as in the case of the work of pollack \cite{pollack2017spin}, the expected functional equation and other analytic properties of the spinor $L$-function in question by the integral representation of $L^{(M)}(s_0, \pi, \mrm{Spin})$ obtained in \cite{eischen2024algebraicity}. 

The second application is to the representation numbers of one quadratic form by another. This is possible since we have included non-cusp forms in our setup. We show that any `primitive' lattice of odd level represents infinitely many `fundamental' lattices of smaller orders. See Corollary~\ref{spl-val}. This result may have uses in the theory of quadratic forms. Certainly this seems not easy to obtain without the use of modular forms.

It would be definitely very interesting to investigate the analogous question for other kinds of congruence subgroups of $\spnz$ like the paramodular groups, or for the unitary (Hermitian modular) groups or the orthogonal groups $O(2,n)$ etc. We expect that the ideas in this paper to carry over, but one may expect some specific obstacles to be overcome. We would expect an analogue of Theorem~\ref{genres-intro} to hold in these cases.
One might also expect that primitive non-zero Fourier coefficients might imply the existence of a non-zero diagonal Fourier coefficient, perhaps with pairwise distinct odd primes on the diagonal.

\subsection*{Acknowledgments}
{\small
It is a pleasure for the authors to thank Siegfried B\"ocherer, Jan H. Bruinier, Ellen Eischen, T. Ibukiyama and Sydney Washburn for conversations on the topic of the paper and for their interest and encouragement.
S.D. thanks IISc. Bangalore, UGC Centre for Advanced Studies, DST, India, 
%only 
for financial support.
P.A. thanks Ahmedabad University for providing excellent infrastructure and research support.
}

\section{Notations and Preliminaries}
\label{prelim}
\begin{enumerate}
    \item We use $\z, \Q, \mbb R$, and $\complex$ to denote the integers,
rationals, reals, and complex numbers, respectively.

\item For any $z\in \mbb C$, we write $e(z):= e^{2\pi i z} $ and $e_c(z):= e^{2\pi i z/c} $ for any real $c \neq 0$.

\item For a commutative ring $R$ with unit, $M_{n,m}(R)$ and $M_n( R)$ denote the set of $n \times m$  and $n\times n$ matrices over $R$, respectively. When $n=1$, we write $M_{1, m}(R)= R^m$.  $\mrm{GL}_n (R)$ denotes the group of invertible elements in $M_n(R)$. We will denote the transpose of $A$ by $\tp{A}$. For matrices $A$ and $B$ of appropriate size, we write $A[B]:=B^tAB$.

\item  For any $1\le r< n$, we use the following embeddings.
\begin{align}\label{GL-Emb-Up}
    \mrm{GL}_{r}(\mbb R) \hookrightarrow \mrm{GL}_{n}(\mbb R): g\mapsto \smat{g}{0}{0}{1_{n-r}}=g^*,
\end{align}
 \begin{align}\label{GL-Emb-Down}
    \mrm{GL}_{r}(\mbb R) \hookrightarrow \mrm{GL}_{n}(\mbb R): g\mapsto \smat{1_{n-r}}{0}{0}{g}=g_*,
\end{align}
\begin{align}\label{Sp-Emb-Up}
    \mrm{Sp}_{r}(\mbb R) \hookrightarrow \mrm{Sp}_{n}(\mbb R): g=\smat{A}{B}{C}{D}\mapsto \begin{psmallmatrix}
            A&0&B&0\\
            0&1_{n-r}&0&0\\
            C&0&D&0\\
            0&0&0&1_{n-r}
    \end{psmallmatrix}=g^{\uparrow},
\end{align}
\begin{align}\label{Sp-Emb-Down}
    \mrm{Sp}_{r}(\mbb R) \hookrightarrow \mrm{Sp}_{n}(\mbb R): g=\smat{A}{B}{C}{D}\mapsto \begin{psmallmatrix}
            1_{n-r}&0&0&0\\
            0&A&0&B\\
            0&0&1_{n-r}&0\\
            0&C&0&D
    \end{psmallmatrix}=g^{\downarrow}.
\end{align}
\item We say a matrix $M=(M(i,j))$ is half--integral if $2M(i,j),\;M(i,i)\in \z$. Let  $N\ge 1$ and $A, B$ be two half--integral matrices. We say $A\equiv B \bmod N$ if $A=B+N\cdot C$ for some half--integral $C$.
\end{enumerate}
   
\subsection{Siegel Modular Forms} \label{prelim-smf} 
Let $N$ be any integer and $L$ be its square-free part. Let $M$ be any divisor of $N$. Then we consider the following congruence subgroups of $\sptwo$.
\begin{align}
    \Gamma_0^{(n)}(N)&:=\{\smat{A}{B}{C}{D}\in\mrm{Sp}_n(\mbb Z): C\equiv 0\bmod N\};\\
\Gamma_1^{(n)}(N, M)&:= \{\smat{A}{B}{C}{D}\in\mrm{Sp}_n(\mbb Z): C\equiv 0\bmod N, A\equiv 1_n\bmod M  \};\\
\Gamma_1^{(n)}(N)&:= \{\smat{A}{B}{C}{D}\in\mrm{Sp}_n(\mbb Z): C\equiv 0\bmod N, A\equiv 1_n\bmod N \}.
\end{align}
Let $\mc E \subset\sptwo$ be the group of order $2^n$ as in subsection~\ref{sec:mainresults}. %consisting of all matrices of the form $\displaystyle \smat{A}{0}{0}{A^{-1}}$ such that $A $ is a diagonal matrix with diagonal entries $\pm 1$.
We define the group $\Gamma_\pm^{(n)}(N,M)$ by the internal semi-direct product 
\begin{align}
   \Gamma_\pm^{(n)}(N,M):= \mc E \ltimes \Gamma_1^{(n)}(N,M),
\end{align}
where $\mc E$ acts on $\Gamma_1^{(n)}(N,M)$ by conjugation. It is clear that every element of $\Gamma_\pm^{(n)}(N,M)$ can be uniquely written as $\varepsilon \cdot \gamma$ where $\varepsilon \in \mc E$ and $\gamma \in \Gamma_1^{(n)}(N,M)$. Note that $\mc E \cap \Gamma_1^{(n)}(N,M)=\{I\}$ and
\begin{align}
   \Gamma_1^{(n)}(N)\subset\Gamma_{\pm}^{(n)}(N,M)\subset \Gamma_0^{(n)}(N). 
\end{align}
When $M=L$, the square-free part of $N$, the group $\Gamma_\pm^{(n)}(N,L)$ is already defined in subsection~\ref{sec:mainresults}.

Let $\rho$ be a polynomial, and not necessarily irreducible representation $\rho \colon \mrm{GL}_n(\mbb C)\rightarrow \mrm{GL}(V)$ where $V$ is finite-dimensional with $m=\dim(V)$. The largest nonnegative integer $k$ such that $\det^{-k}\otimes \rho$ is
still polynomial, will be called the (determinantal) weight $k(\rho)$. We tacitly use the fact that this weight does not decrease if we tensor $\rho$
with another polynomial representation or restrict it to some 
$\mrm{GL}_{n'}(\mbb C)$ sitting inside $\mrm{GL}_{n}(\mbb C)$ as an algebraic  subgroup. This follows easily by looking at the entries of $\rho(g)$ in any matrix realization of $\rho$.

Let $\mbb H_n$ denote the Siegel upper half space of degree $n$. The symplectic group $\mrm{Sp}_n(\mbb R)$ acts on $\mbb H_n$ by
$Z\mapsto \gamma \langle Z \rangle=(AZ+B)(CZ+D)^{-1}$; for a polynomial representation $\rho$ 
with values in $\mrm{GL}(V)$ we define the stroke operator action on $V$-valued functions $F$ on $\mbb H_n$ by
\[ (F\mid_{\rho}\gamma)(Z):=\rho(CZ+D)^{-1}F(\gamma \langle Z \rangle). \]
A Siegel modular form of degree $n$ and automorphy factor $\rho$ for the group $\Gamma_\pm^{(n)}(N, M)$
is then a $V$-valued holomorphic function $F$ on $\mbb H_n$ satisfying
$F\mid_{\rho}\gamma= F$ for all $\gamma\in \Gamma_\pm^{(n)}(N, M)$
with the standard additional condition of \textit{boundedness at cusps} in degree $1$.
Any such $F$ has a Fourier expansion given by
\begin{align}
    F(\mc Z)= \sumn_{T\in \Lambda_n}A(F, T)e(\tr(T\mc Z)),
\end{align}
where $\Lambda_n$ (resp. $\Lambda_n^+$) denotes the set of all $n\times n$ symmetric, positive semi--definite (resp. positive--definite), half--integral matrices. If $F$ is cuspidal, then this summation is supported on $\Lambda_n^+$. 

We denote by $M_{\rho}(\Gamma_\pm^{(n)}(N, M))$ the vector space of all such functions and by
$S_{\rho}(\Gamma_\pm^{(n)}(N, M))$ the subspace of cusp forms; if $\rho=\det^k$ (that is, scalar valued), we write
as usual $M_k(\Gamma_\pm^{(n)}(N, M))$ and $S_k(\Gamma_\pm^{(n)}(N, M))$ . 

Let $U\in \glnz$ be such that $\smat{(U^t)^{-1}}{0}{0}{U}\in \Gamma_1^{(n)}(N)$. Then we have  $F(Z[U^{-1}])= \rho(U) F(Z)$. Comparing the Fourier coefficients, we get
\begin{align}\label{glnz-eq}
A(F, T[U])= \rho(U) A(F, T).
\end{align}

For any $T\in \Lambda_n$, define the content $\mf c(T)$ as $\mf c(T):=\max\{d\in \mbb N \;| \;d^{-1}T \in \Lambda_n\}$. We say that $T$ is primitive if $\mf c(T)=1$.

\subsection{Jacobi forms for the group \texorpdfstring{$\Gamma_\pm^{(n)}(N, M)$}{gamma}}\label{jacobi}

\subsubsection{Vector-valued Jacobi forms}\label{jacobi-vector}
Let us decompose
${\mathcal Z}\in \mbb  H_n$ into blocks as follows:
\begin{equation}\label{decompo}
\mathcal Z=\begin{pmatrix}\tau &z\\
z^t & Z \end{pmatrix}  \qq (z \in M_{1, n-1}(\mbb C), Z\in \mbb H_{n-1}).
\end{equation}
Clearly, every $F\in M_\rho( \Gamma_\pm^{(n)}(N, M))$ has a 
Fourier-Jacobi expansion with respect to the decomposition above: 
\[F({\mathcal Z})=\sumn_{S\in \Lambda_{n-1}} \phi_S(\tau,z) e(SZ).\]
The $\phi_S$ are then ``Jacobi forms'' of automorphy factor $\rho$ and
index $S$, i.e., the functions $\psi(\mc Z)=\psi_S({\mathcal Z}):= \phi_T(\tau,z)e (TZ)$
on $\mbb H_n$ are holomorphic, satisfy $\psi\mid_{\rho} \gamma=\psi$ for all 
$\gamma\in C_{n, n-1}(N, M;\mbb Z)$, where
\begin{align}
  C_{n, n-1}(N, M;\mbb Z) &:=  C_{n, n-1}(\mbb Z) \cap  \Gamma_\pm^{(n)}(N, M) \\ 
  &= \{\left(\begin{smallmatrix} A & B\\
C & D\end{smallmatrix}\right)\mid (C,D)= 
\left(\begin{smallmatrix} * & 0 &{}& * & D_2\\
0 & 0 & {} & 0 & D_4\end{smallmatrix}\right), C\equiv 0\bmod N, D\equiv 1\bmod M\}; \label{kl-congr}
\end{align}
and the Klingen parabolic subgroup $C_{n, n-1}(\mbb Z)$ is defined as in \eqref{kl-congr}, but without the congruence conditions.
In \eqref{kl-congr}, $*$ denotes some scalar entries. 

In addition, $\psi$ satisfies the boundedness conditions (Fourier expansion at cusps). For example, at the cusp $\infty$, we have an expansion of the form
\[ \psi(\mc Z)= \sum_{T = \left( \begin{smallmatrix} n & r/2 \\ r^t/2 & S \end{smallmatrix}\right) \in \Lambda_n} a_\psi(T) e(T \mc Z).\]
Note that this definition of vector-valued Jacobi forms does not agree
with the one in \cite{ziegler1989jacobi}; for degree 2, our definition is the same as in 
\cite{ibukiyama2011generalization}.

\subsubsection{The scalar-valued Jacobi forms}\label{jacobi-scalar}
The case $\rho=\det^k$, i.e., the scalar-valued case, we encounter Jacobi forms of weight $k$, index $S\in \Lambda_g$ and lattice $\mc L_M:=M\cdot M_{g,n}(\mbb Z)\times M_{g,n}(\mbb Z)$ for the group $\Gamma^{n,g,J}_\pm(N, M):=\Gamma_\pm^{(n)}(N,M)\ltimes \mc L_M$. Since these objects will play a central role in this paper, we recall their definition and some related notation. 

A holomorphic function $\phi$ on $\hn\times M_{g,n}(\mbb C)$ is called a Jacobi form of weight $k$, index $S$ and lattice $\mc L_M$ for the group $\Gamma^{n,g,J}_\pm(N, M)$ if
\begin{enumerate}
    \item $\phi|_{k,S}\gamma = \phi \q \forall \gamma\in \Gamma^{n,g,J}_\pm(N, M)$.
    \item $\phi$ is bounded at the cusps of $\Gamma^{n,g,J}_\pm(N, M)\backslash \hn$ (automatic for $n\ge 2$, due to the K\"{o}cher principle).
\end{enumerate}

Any such $\phi$ as a Fourier expansion given by
\begin{equation}\label{jacobiFE}
    \phi(\tau,z)=\underset{4T-S^{-1}[R]\ge 0}{\sumn_{T\in \Lambda_n}\sumn_{R\in M_{g,n}(\z)}}C_\phi(T,R) e(\tr (T\tau+R^tz)).
\end{equation}
If the Fourier coefficients survive only for $4T-S^{-1}[R]>0$, then $\phi$ is called a Jacobi cusp form. $J_{k, S}(\Gamma^{n,g,J}_\pm(N, M), \mc L_M)$ (resp. $J^{cusp}_{k, S}(\Gamma^{n,g,J}_\pm(N, M),\mc L_M)$) denotes the space of Jacobi forms (resp. Jacobi cusp forms) of weight $k$ and index $S$ on $\hn \times M_{g,n}(\mbb C)$. We refer to $M$ as the ``lattice-level'' of $\phi$.

Let $\phi_S\in J_{k, S}(\Gamma_\pm^{n,g,J}(N,M), \mc L_M)$, then for any $\lambda\in M_{g,n}(\mbb Z)$ we have
\begin{align}
   \phi_S(\tau, z+M\lambda \tau)=  e(-\tr(M^2S[\lambda]\tau+2Mz^tS\lambda)) \phi_S(\tau,z).
\end{align}
Comparing the Fourier coefficients, we see that 
\begin{align} \label{fc-equiv}
 C_{\phi_S}(T_1, r_1)=C_{\phi_S}(T_2, r_2)   \text{ if } r_1\equiv r_2\bmod 2MS \text{ and } T_1-S^{-1}[r_1/2]=T_2-S^{-1}[r_2/2] .
\end{align} 

As a consequence, by rearranging the Fourier expansion in \eqref{jacobiFE}, we see that the scalar-valued Jacobi forms  $\phi_S\in J_{k, S}(\Gamma_\pm^{n,g,J}(N,M),\mc L_M)$ admit a ``theta expansion''
\begin{align}\label{theta-decomp}
    \phi_S(\tau, z)=\sumn_{\mu\in \mc N_{2MS}^{g,n}} h_{\mu}(\tau)\cdot \Theta_{MS, \mu}(M\tau,z),
\end{align}
where $\mc N_{2MS}^{g,n}:= M_{g,n}(\mbb Z)/(2MS)\cdot  M_{g,n}(\mbb Z)$
and
\[\Theta_{S, \mu}(\tau,z)=\sumn_{R\in M_{g,n}(\z)}
e(\tr( S[R+{\overline\mu}]\tau+ 2 z^t S(R+{\overline\mu})).\]
Here we use ${\overline\mu}:= (2S)^{-1}\cdot \mu$. The Fourier expansion of $h_{\mu}$ is of shape
\begin{equation} \label{hmu-fe}
h_{\mu}(\tau)=\sumn_{T -(2S)^{-1}[\mu]\ge 0} c_{\phi_S}(T ,\mu) e(\tr(T- (2S)^{-1}[\mu])\cdot \tau)
\end{equation}
and its Fourier coefficients are given  by 
\begin{equation}
c_{\phi_S}(T ,\mu) = A\left(F, \left(\begin{smallmatrix} T & \mu/2\\
\mu^t/2 & S\end{smallmatrix} \right) \right)
\label{relation}
\end{equation}
provided that $\phi_S$ is the Fourier-Jacobi coefficient of
some scalar-valued Siegel modular form $F\in M_k(\Gamma_\pm^{(n)}(N,M))$. In this case, the functions $\phi_S \in J_{k, S}(\Gamma_\pm^{n,g,J}(N,M), \mc L_M)$, which follows in a standard manner from the automorphy of the associated function $\psi_S$ and \eqref{kl-congr}.

% The automorphy properties of $h_\mu$ will be discussed later.

We have the following transformation properties of $\Theta_{S, \mu}(\tau,z)$ (see \cite{ziegler1989jacobi})
\begin{align}\label{theta-trans}
\Theta_{S, \mu}(\tau+1,z)&=  e(\tr(S^{-1}[\mu/2])) \Theta_{S, \mu}(\tau,z)\\
\Theta_{S, \mu}(-\tau^{-1},z \tau^{-1})&=\det(S)^{-n/2}\det(\tau/i)^{g/2}e(S[z/2]\tau^{-1}))\sum_{\eta\in \mc N_{2S}^{g,n}}e_2(-\tr(\eta^t S^{-1} \mu))\Theta_{S, \eta}(\tau,z). \label{S-theta-trans}
\end{align}
Let $\gamma=\smat{A}{B}{C}{D}\in \mrm{Sp}_n(\mbb Z)$. Since $\mc S=\smat{0_n}{-1_n}{1_n}{0_n}$ and $\mc T=\smat{1_n}{1_n}{0_n}{1_n}$ generate $\mrm{Sp}_n(\mbb Z)$, there exist complex numbers $\varepsilon_S(\mu, \nu; \gamma)$ such that
\begin{align}\label{thetatransSpn}
\Theta_{S, \mu}(\gamma\tau,z(C\tau+D)^{-1})= \det(C\tau+D)^{g/2}e(\tr(S[z](C\tau+D)^{-1}C))\sum_{\nu\in \mc N_{2S}^{g,n}}\varepsilon_S(\mu, \nu; \gamma)\Theta_{S, \nu}(\tau,z).
\end{align}

Let $\phi \in J_{k,T}(\Gamma_\pm^{1,g,J}(N,M), \mc L_M)$ and consider the theta decomposition of $\phi$ as in \eqref{theta-decomp}. We define the `primitive' theta components of $\phi$ as below. First, let
\begin{align}\label{def:d_T}
d_T:=   \begin{cases}
        d(T) & \text{ if   }  n  \text{ is even }; \\
        4 d(T) & \text{ if   } n \text{ is odd }.
    \end{cases}
\end{align}

\begin{definition}[Primitive theta components] \label{prelim:prim}
   For $ \mu\in \mc N_{2T}^{g,1}$, we say that the theta component $h_\mu$ is $T$-primitive if $(2T)^{-1}[\mu]$ has exact denominator  $d_T$.
\end{definition}

\subsection{Construction of a lower degree form from degree \texorpdfstring{$n$}{n} form.} \label{sec:Fcirc}

(1) Let $F\in \Mrho{n}{N}$. We now show how to construct a degree $r<n$ modular form from $F$ as in \cite{boche-das}. For any $\mc Z\in \mbb H_n$, let us write
\begin{equation}
    \mc Z=\smat{Z}{ z}{z^t}{Z'}; \text{ with } Z\in \mbb H_{n-r}, z\in M_{n-r,r}(\mbb C), Z'\in \mbb H_{r}.
\end{equation}
Then $F$ has a Fourier-Jacobi expansion given by
\begin{align}
    F(\mc Z)=\sumn_{S\in \Lambda_{r}} \varphi_S(Z, z) e(\tr(SZ')).
\end{align}
Let $z=(z_{ij})$ with $1\le i\le n-r$, $n-r+1\le j\le n$ and $z_{ij}\in \mbb C$. Consider the Taylor expansion of $F$ w.r.t $z$ and write
\begin{equation}
    F(\mc Z)= \sumn_{\lambda\in \mbb N^{n-r,r}} F_{\lambda}(Z, Z') z^\lambda,
\end{equation}
where for $\lambda=(\lambda_{ij})\in \mbb N^{n-r, r}$, $z^\lambda=\prod z_{ij}^{\lambda_{ij}}$.

Let $\nu(\lambda)=\sum_{i, j} \lambda_{ij}$ and $\nu_0=\min \{\nu(\lambda): F_\lambda\neq 0\}$. We now look at the Taylor coefficients of homogeneous degree $\nu_0$ and study a
polynomial of homogeneous degree $\nu_0$ in variables $X=(X_{ij})$, where $1\le i \le r$, $n-r+1\le j\le n$. 
\begin{align}
    F_{(r)}(Z, Z')&:= \sum_{\lambda: \nu(\lambda)=\nu_0} F_{\lambda}(Z, Z') X^\lambda= \sum_{S\in \Lambda_{r}} \sum_{\lambda: \nu(\lambda)=\nu_0} \left.\frac{\partial^\lambda}{\partial z^\lambda}\varphi_S(Z, z)\right |_{z=0} e(SZ') X^\lambda.
\end{align}
For $X=(X_{ij})$, let $\mbb C[X]_{\nu_0}$ denote the $\mbb C$ vector space of homogeneous polynomials of degree $\nu_0$. We can then view  $F_{(r)}$ as a function on $\mbb H_{n-r} \times \mbb H_r$ with values in $V\otimes \mbb C[X]_{\nu_0}$.

(2) We also construct a new modular form $F^{(r)}$ as below. Let us write
\begin{equation}
    \mc Z=\smat{Z}{ z}{z^t}{Z'}; \text{ with } Z\in \mbb H_r, z\in M_{r,n-r}(\mbb C), Z'\in \mbb H_{n-r}.
\end{equation}
We write $T= \smat{S}{s/2}{s^t/2}{S'}$ with $S\in \Lambda_r$, $S'\in \Lambda_{n-r}$ and define functions $\psi_S(z, Z')$ as below.
\begin{align}
    F(\mc Z)= \sumn_{S\in \Lambda_{r}}\psi_S(z, Z')e(\tr(SZ)).
\end{align}
We consider the Taylor expansion of $F$ w.r.t $z$ and write
\begin{equation}
    F(\mc Z)= \sumn_{\lambda\in \mbb N^{r,n-r}} F_{\lambda}(Z, Z') z^\lambda.
\end{equation}

Let $\nu(\lambda)=\sum_i \lambda_i$ and $\nu_0=\min \{\nu(\lambda): F_\lambda\neq 0\}$. As above, we now define
\begin{align}
   F^{(r)}(Z, Z')&= \sum_{\lambda: \nu(\lambda)=\nu_0} F_{\lambda}(Z, Z') X^\lambda= \sum_{S\in \Lambda_{r}} \sum_{\lambda: \nu(\lambda)=\nu_0} \left.\frac{\partial^\lambda}{\partial z^\lambda}\psi_S(z, Z')\right |_{z=0} e(SZ) X^\lambda.
\end{align}

Let $\rho^{(r)}$ and $\rho_{(r)}$ denote the restrictions of $\rho$ to $\mrm{GL}_{r}$ w.r.t the embedding in \eqref{GL-Emb-Up} and \eqref{GL-Emb-Down}, respectively. Further, let $\mrm{sym}^{t}_r:\mrm{GL}_r(\mbb C)\longrightarrow \mrm{GL}(\mbb C[X_1,...X_r]_t)$ denote the $t$-th symmetric power representation of the identity representation of $\mrm{GL}_r(\mbb R)$ on the space of homogeneous polynomials in $r$ variables of degree $t$; and define
\begin{align}
    \sigma_{\nu_0}:=\bigoplus_{\sum_m t_m=\nu_0} \left(   \mrm{sym}^{t_1}_r \otimes \mrm{sym}^{t_2}_r \cdots \otimes \mrm{sym}^{t_{n-r}}_r \right).
\end{align}
Then we have the following (see \cite[Proposition 3.1]{boche-das} for $r=n-1$ case).
\begin{proposition}\label{prop:Fcirc}
    Let $F\in \Mrho{n}{N}$. Then $F_{(r)}(Z, Z')$ (resp. $F^{(r)}(Z,Z')$), as a function of $Z'$ (resp. $Z$), is in $M_{\rho_{(r)}\otimes \sigma_{\nu_0}}(\Gamma_\pm^{(r)}(N,L))$ (resp. $M_{\rho^{(r)}\otimes \sigma_{\nu_0}}(\Gamma_\pm^{(r)}(N,L))$) . If $F\neq 0$, then $F_{(r)}(Z_0, Z')\neq 0$ and $F^{(r)}(Z, Z_0')\neq 0$ for some $Z_0, Z'_0\in \mbb H_{n-r}$ respectively.
\end{proposition}
\begin{proof}
The proof follows as in \cite[Proposition 3.1]{boche-das} once we understand the representation and use the embeddings in \eqref{Sp-Emb-Up} and \eqref{Sp-Emb-Down}. Towards this, it is best to write 
    \begin{align}
        z= (a_1, a_2,...,a_{n-r})^t; \q g\cdot z= (g \cdot a_1,...., g \cdot a_{n-r} )^t
    \end{align}
for the natural action of $g \in \mrm{GL}_r(\mbb C)$ where $a_j \in M_{r,1}(\mbb C)$. The assertion then follows by decomposing each polynomial in $z$ as a sum of homogeneous polynomials in $a_j$ of degrees $0 \le t_j \le \nu_0$ such that $\sum_{j=1}^{n-r} t_j=\nu_0$. We use the fact that direct sums distribute over the tensor products.
\end{proof}
\begin{remark}
    Note that we did not discuss and also do not need any automorphy properties of the functions $\psi_S$ mentioned above. Most likely, they are Jacobi forms.
\end{remark}
%%%%
\section{Preparatory results}
In this section, we collect some preparatory results that will be used in the proofs of the main results. At various points in the proofs, we need to construct new modular forms from old by restricting the Fourier expansion to be congruent to a fixed (matrix) index $\mod M$ for some integer $M$. We need to know the structure of the level and the congruence subgroup after this procedure, crucially, for induction to work.

\begin{lemma}\label{lem:FTp}
    Let $F \in \Mrho{n}{N}$ and $T\in \Lambda_n$. For any prime $p$, define 
    \begin{align}
        F_{T,p}(Z):= \sum_{S\equiv T \bmod p} A(F,S) e(\tr(SZ)).
    \end{align}
Then $F_{T,p}\in M_\rho(\Gamma_\pm^{(n)}(Np^2/(N,p), Lp/(L,p)))$.
\end{lemma}
\begin{proof}
    Let $F_{T,p}$ be defined as above. Then we can write
    \begin{align}
        F_{T,p}= p^{-n(n+1)/2} \sum_{B=B^t, B \bmod p} e(-\tr(TB/p))F|_\rho\smat{I_n}{B/p}{0}{I_n}.
    \end{align}
It is enough to show that $F_{T,p}$ is invariant w.r.t $\Gamma_1^{(n)}(Np^2/(N,p), Lp/(L,p))$ and $\mc E$ separately.

For each $B$, $F|_\rho\smat{I_n}{B/p}{0}{I_n}$ is a modular form for the group 
\begin{align}
\Gamma_B'&:=\smat{I_n}{B/p}{0}{I_n}^{-1}\Gamma_1^{(n)}(N, L)\smat{I_n}{B/p}{0}{I_n}\cap \Gamma_1^{(n)}(N,L)\\
&=\{\smat{A}{*}{C}{D}\in \Gamma_1^{(n)}(N,L): C\equiv 0\bmod p^2/(N,p), AB\equiv BD\bmod p\},
\end{align}
for the same representation $\rho$.
Further, for each $B\mod p$, $\Gamma_1^{(n)}(Np^2/(N,p), Lp/(L,p))\subseteq \Gamma'_B$ and thus we get that $F_{T,p}$ is invariant w.r.t $\Gamma_1^{(n)}(Np^2/(N,p), Lp/(L,p))$.

Let $\varepsilon=\smat{A}{0}{0}{A^{-1}}\in\mc E$, then we see that
\begin{align}
    \smat{I_n}{B/p}{0}{I_n}\cdot \varepsilon= \varepsilon\cdot \smat{I_n}{A^{-1}BA^{-1}/p}{0}{I_n}.
\end{align}
Note that $B':=A^{-1}BA^{-1}$ is symmetric, and as $B$ varies $\mod p$, $B'$ also varies $\mod p$. We thus get the required invariance w.r.t $\mc E$, and this completes the proof of the lemma.
\end{proof}
\begin{lemma}\label{lem:FTN}
     Let $F \in \Mrho{n}{N}$ and $T\in \Lambda_n$. Then $F_T:=F_{T, N}\in \Mrho{n}{N^2}$.
\end{lemma}
\begin{proof}
    Let $N=p_1.....p_t$ ($p_i$ not necessarily distinct). The proof now follows from Lemma \ref{lem:FTp} by induction on $t$.
\end{proof}

An important step in our proofs is to descend from vector-valued Jacobi forms in the sense of Section~\ref{jacobi-vector} to scalar-valued Jacobi forms in the sense of Section~\ref{jacobi-scalar}. The following proposition facilitates this step. In fact, we prove this more generally for the smaller Jacobi group $\Gamma_1^{1,g,J}(N):=\Gamma_1^{(1)}(N)\ltimes \mc L_N$.
\begin{proposition} \label{sca-j}
If $\varphi_T \neq 0$ is a vector-valued Jacobi form of index $T$ with respect to $\rho$ for the group $\Gamma_1^{1,g,J}(N)$, then there exist a component of $\varphi^{[r]}_T \neq 0$ of $\varphi_T$ which is a scalar-valued Jacobi form for the group $\Gamma^{1,g,J}_1(N)$ (of integral weight $k \geq k(\rho)$).
\end{proposition}

\begin{proof}
  The proof is the same as in   \cite[Proposition 3.4]{boche-das}, with the discrete subgroup replaced with that of ours. However, the only point in the argument where any specific feature of the discrete subgroup was used was the choice of 
  $M=\begin{psmallmatrix} 
  U^t & 0\\
0 & U^{-1}\end{psmallmatrix} \in \spnz $, with $U=\begin{psmallmatrix} 1 & N\ell\\
0 & 1_{n-1}
\end{psmallmatrix}$,  ($\ell \in \z^{n-1}$);
which leads to the transformation law (w.r.t. the lattice $N \z^{n-1} \times \z^{n-1}$)
\begin{equation}
\varphi_T(\tau,N\ell\tau+z)=\rho(U^{-1})\varphi_T(\tau,z)
e\left(-N^2T[\ell^t]\tau -2N\tau z T\ell^t\right) , \label{Lie}
\end{equation}
from which the proof followed. However since $M \in \Gamma^{(n)}_1 (N)$, we are fine.
For the rest of the proof, only the properties of the representation theory of the underlying Lie group $\glnc$ were used. Therefore, we omit the details.
\end{proof}
As an immediate corollary, we have the following.
\begin{corollary}
  If $\varphi_T \neq 0$ is a vector-valued Jacobi form of index $T$ with respect to $\rho$ for the group $\Gamma_\pm^{1,g, J}(N,L)$, then there exist a component of $\varphi:=\varphi^{[r]}_T \neq 0$ of $\varphi_T$ which is a scalar-valued Jacobi form for the group $\Gamma^{1,g,J}_\pm(,L)$ (of integral weight $k \geq k(\rho)$). Further, the Fourier coefficients $c_\varphi(n,r)$ of $\varphi$ satisfy $c_\varphi(n,r)= (-1)^{k}c_\varphi(n,-r)$.
\end{corollary}
\begin{proof}
    The existence of $\varphi$ follows from Proposition \ref{sca-j} since $\Gamma_1^{1,g,J}(N)\subset\Gamma_\pm^{1,g, J}(N,L)$.  For the second part, we note that since $-1 \in \Gamma^{(1)}_\pm(N, L)= \big \lan \{ \pm \}, \Gamma_1^{(1)} (N, L)\big \ran$, the symmetry property of $\varphi$, that is, $c_\varphi(n,r)= (-1)^{k}c_\varphi(n,-r)$ holds. This can be seen from the Fourier expansion of $\varphi_T$. 
\end{proof}

Next, we have the following non-vanishing result, which will serve as the base case of induction in the proofs that follow.
\begin{lemma}\label{ind-basecase}
    Let $f\in M_\rho(\Gamma_1^{(1)}(N))$ ($k(\rho) \ge 2$) be such that $a(f, n_0)\neq 0$ for some $(n_0,N)=1$. Then there exist infinitely many odd, square-free $n$ with $(n,N)=1$ such that $a(f, n)\neq 0$.
\end{lemma}
\begin{proof}
    Let $f=(f_1,...,f_m)^t \in M_\rho(\Gamma_1^{(1)}(N))$, where $m=\dim(\rho)$ so that $f_j$ are necessarily of some scalar weight $k_j \ge 2$. See \cite[Section~4.3]{boche-das}. By the hypothesis, there exists an $f_j$ such that $a(f_j, n_0)\neq 0$ for some $j$ and $(n_0, N)=1$. Let $g:=\sum_{(n,N)=1}a(f_j, n) q^n$. Then $g\neq 0$ and $a(g,n)=0$ for $(n, N)>1$. The proof now follows from \cite[Theorem 4.6]{boche-das}.
\end{proof}
We also make the following observation for two matrices $T_1, T_2\in \Lambda_n$ that are congruent modulo $N$.
\begin{lemma}\label{lem:detcong}
    Let $N\ge 1$ and $T_1, T_2\in \Lambda_n$ be such that $T_1\equiv T_2 \bmod N$. Then
    \begin{enumerate}
        \item $(\mf c(T_1), N)= (\mf c(T_2), N)$.
        \item $d(T_1)\equiv d(T_2) \bmod N$.
    \end{enumerate}
\end{lemma}
\begin{proof}
The proof is immediate by writing  $T_1= T_2+N\cdot T_3$.
\end{proof}
%%===================================

%%===================================
\section{Gauss sums and coefficients of the Theta representation}
For $M\in \mrm{Sp}_n(\mbb Z)$, $S\in \Lambda^+_g$ and $\mu,\nu\in \mc N_{2S}^{g,n}$, let $\varepsilon_S(\mu, \nu; M)\in \mbb C$ be as in \eqref{thetatransSpn}. We have the following relation for $\varepsilon_S(\mu, \nu; M)$.
\begin{lemma}
If $M_1, M_2\in \mrm{Sp}_n(\mbb Z)$, then there exists a constant $\omega(M_1, M_2)\in \mbb C$ with $|\omega(M_1, M_2)|=1$ such that
\begin{align}\label{epsilonmunu}
\varepsilon_S(\mu, \nu; M_1M_2)= \omega(M_1, M_2)\sumn_{\eta\in \mc N_{2S}^{g,n}}\varepsilon_S(\mu, \eta; M_1)\varepsilon_S(\eta, \nu; M_2).  
\end{align}
\end{lemma}
\begin{proof}
For $M=\smat{A}{B}{C}{D}\in \mrm{Sp}_n(\mbb Z)$ and $\tau\in \mbb H_n$, let $j(M; \tau)= \det(C\tau+D)$. Then using \eqref{thetatransSpn} for $M_1M_2$, $M_1$ and $M_2$ in sequence at $z=0$, we get
\begin{align}
\varepsilon_S(\mu, \nu; M_1M_2)= \frac{j(M_1; M_2\tau)^{g/2}j(M_2;\tau)^{g/2}}{j(M_1M_2; \tau)^{g/2}}    \sum_{\eta\in \mc N_{2S}^{g,n}}\varepsilon_S(\mu, \eta; M_1)\varepsilon_S(\eta, \nu; M_2).
\end{align}
Using the co-cycle condition, we see that the ratio at the front is independent of $\tau$ and is of absolute value $1$ (see \cite[Section 3]{freitag2006singular}). This completes the proof of the lemma.
\end{proof}
\subsection{Higher dimensional Gauss sums:}\label{sec:HDGauss}
Let $M\in \Lambda^+_g$, $b\in\mbb Z^g$ and $a\in \mbb Z$. We define the generalized quadratic Gauss-sum of degree $g$ associated to $a,b, M$ by
\begin{align} \label{gauss}
    G_g(a, b, 2M):= \frac{1}{\det(2M)}\sum_{\bm\eta\in \mc N_{2M}^{g,1}} e_2\left(a (2M)^{-1}[\eta]+2 \eta^t(2M)^{-1}b\right).
\end{align}
When $b=0$, we have the following non-vanishing property. First, we recall that $d(M)=|\mrm{disc}(2M)|$.
\begin{proposition}\label{gauss-nonvan}
  Let $M\in \Lambda_g^+$ be such that $d(M)$ is odd, square-free and $a$ be any integer such that $(a, \det(2M))=1$. Then $G_g(a, 0, 2M)\neq 0$. In particular, there exists an $m\in \mbb Z$ with $(m, \det(2M))=1$ such that
  \begin{align}
  G_g(a, 0, 2M)=\begin{cases}
      \frac{1}{2}G_1(am, 0, 4d(M)) & \text{ when } g \text{ odd};\\
      G_1(am, 0, d(M))&\text{ when } g \text{ even }.
  \end{cases}    
  \end{align}

\end{proposition}
   
\begin{proof}

\textbf{$g$ is odd:} As in \cite{boche-das}, let $U\in\slnz$ be such that $\widetilde{M}=(2M)[U]$ satisfies
\begin{align*}
    \widetilde{M}\equiv\mrm{diag}(*, *,..., \zeta d(M)) \bmod d(M)^f \q \text{ for } f\ge 2.
\end{align*}
Then the following set can be chosen as representatives for $\mbb Z^g/ \widetilde{M} \mbb Z^g$
\begin{align}\label{coset-g-odd}
    \{ \bm{ \tilde \eta}\}:= \{ (0,0,....,\eta): \eta\bmod 2d(M)\}.
\end{align}
Then the set $\{ U^t\bm{\tilde \eta}\}$ forms a set of coset representatives for $\mbb Z^g/(2M)\mbb Z^g$ and we get that
\begin{align}
\bm{ \eta}^t (2M)^{-1} \bm{\eta}= \bm{ \tilde \eta}^t \widetilde{M}^{-1}\bm{ \tilde \eta}.
\end{align}
Further, let $\widetilde{M}^*=(m_{i,j})$ denote the adjoint of $\widetilde{M}.$ That is $\widetilde{M}^*\widetilde{M}=2d(M)\cdot I_g$. From \cite{boche-das}, we have that $(m_{g,g}, 2d(M))=1$ and 
\begin{align}
\bm{ \tilde \eta}^t \widetilde{M}^{-1}\bm{ \tilde \eta}= \frac{1}{2d(M)}m_{g,g}\eta^2.    
\end{align}
Thus we can now write
\begin{align}
    G_g(a, 0, 2M)= \sumn_{\eta\bmod 2d(M)}e_2(am_{g,g}\eta^2/2d(M)).
\end{align}
Note that the set $\{\eta, \eta+2d: \eta\bmod 2d\}$ forms a set of coset representatives for $\bmod 4d$ and $e_2(am_{n,n}(\eta+2d(M))^2/2d(M))=e_2(am_{g,g}\eta^2/2d(M)) $. Thus, we get
\begin{align}
 G_g(a, 0, 2M)= \frac{1}{2}\sum_{\eta\bmod 4d(M)}e_2(am_{g,g}\eta^2/2d(M))=\frac{1}{2}G_1(am_{g,g}, 0, 4d(M)).    
\end{align}
The one-dimensional Gauss sum $G_1(am_{g,g}, 0, 4d(M))$ is non-zero whenever $(a, 2d(M))=1$. Moreover, in this case, we have
\begin{align}
 G_1(am_{g,g}, 0, 4d(M))= (1+i) \epsilon_{am_{g,g}}\sqrt{4d(M)} \left(\frac{am_{g,g}}{4d(M)}\right).  
\end{align}

\textbf{$g$ is even:} As in the previous case, let $U\in\slnz$ be such that $\widetilde{M}=(2M)[U]$ satisfies
\begin{align*}
    \widetilde{M}\equiv\mrm{diag}(*, *,..., d(M)) \bmod d(M)^f \q \text{ for } f\ge 2.
\end{align*}
Then the following set can be chosen as representatives for $\mbb Z^g/ \widetilde{M} \mbb Z^g$
\begin{align}\label{coset-g-even}
    \{ \bm{ \tilde \eta}\}:= \{ (0,0,....,\eta): \eta\bmod d(M)\}.
\end{align}
The set $\{U^t\bm{\tilde \eta}\}$ then forms a set of coset representatives for $\mbb Z^g/(2M)\mbb Z^g$. Thus, we get
\begin{align}
G_g(a, 0, 2M)= \sumn_{\eta\bmod d(M)}e_2(am_{g,g}\eta^2/d(M)).    
\end{align}
As above, let $\widetilde{M}^*$ denote the adjoint of $\widetilde{M}$. Let $\widetilde{M}_{g-1}$ denote the upper $g-1\times g-1$ block of $\widetilde{M}$. Then $m_{g,g}=\det(\widetilde{M}_{g-1})=\det((2M)_{g-1})$. But since $g$ is even, $\det((2M)_{g-1})$ is even. Thus $m_{g,g}\equiv 0\bmod 2$. Let $m_{g,g}=2m'_{g,g}$. Then the right side is just the quadratic Gauss sum $G_1(am'_{g,g}, 0, d(M)) $, which is now non-zero since $d(M)$ is odd. In this case, we have
\begin{equation}
G_1(am'_{g,g}, 0, d(M))= \epsilon_{d(M)}\sqrt{d(M)}   \left(\frac{am'_{g,g}}{d(M)}\right).\qedhere
\end{equation}
\end{proof}
%===============================

%===============================
\section{Theta components of Jacobi Forms}

In this section, we study the theta components of Jacobi forms, as recalled in \eqref{theta-decomp}, in detail. The main aim of this section is to prove the non-vanishing of certain theta components of Jacobi forms belonging to $J_{k,T}(\Gamma_\pm^{1,g,J}(N,L), \mc L_L)$. For the groups of the type $\Gamma_0^{(n)}(N)$ and the standard lattice $M_{g,n}(\mbb Z)\times M_{g,n}(\mbb Z)$, there are several results describing the non-vanishing of the associated theta components (see \cite{anamby2024non}, \cite{boche-das}, \cite{shankhadhar2023sign} etc). However, in our situation, we have to deal with the groups of type $\Gamma_\pm^{(n)}(N, L)$ and lattices of the type $\mc L_L=L\cdot M_{g,n}(\mbb Z)\times M_{g,n}(\mbb Z)$ and the results mentioned above do not fit into this setting. Hence, we give a detailed account of the calculations -- they are crucial in completing our induction procedure to prove Theorem~\ref{diagonal-sqrfree}. The primary objective of this section is to establish the following non-vanishing result.

\begin{proposition}\label{jacobi-theta-main}
     Let $n, k, N$ be positive integers, and $N$ be odd, and $L$ denote the square-free part of $N$. Let $T\in\Lambda_g^+$ be such that 
     \begin{enumerate}
     \item $T$ is diagonal modulo $L$, 
     \item  $d(T)$ is odd, square-free, and 
     \item $(d(T), L)=1$. 
     \end{enumerate}
     Let $0 \neq \phi\in J_{k, T}(\Gamma_\pm^{1,g, J}(N, L), \mc L_L)$ be such that the Fourier coefficients $c_{\phi}(n, \mf r)$ are supported on $\mf r\equiv (0,0,...,r)\bmod L$. Then there exists an $\mf r_0$ of the form $\mf r_0\equiv (0,0,...0)\bmod L$ such that $c_{\phi}(n_0, \mf r_0)\neq 0$.
\end{proposition}

\subsection{The Coset Representatives.} We require a suitable set of coset representatives for $\mc N_{2LT}^{g,1}$ when $(L, d(T))=1$. In this case, $\mc N_{2LT}^{g,1}\cong \mc N_{L}^{g,1}\times \mc N_{2T}^{g,1}$. Thus the set
\begin{align}\label{coset-reps}
  \{(2T)\bm{\mu_1}+L\bm{\mu_2}: \bm{\mu_1} \in  \mc N_{L}^{g,1}, \bm{\mu_2} \in \mc N_{2T}^{g,1} \}
\end{align}
can be taken as a set of coset representatives for $\mc N_{2LT}^{g,1}$ which we adhere to in the calculations given below.

Further, as in the previous section, let $\widetilde{T}=(2T)[U]$  and $\{\bm{\tilde\mu_2}\}=\{(0,0,..\mu_2)^t: \mu_2 \bmod \det(2T)\}$ denote the set of coset representatives for $\mbb Z^g/\widetilde{T}\mbb Z^g$. Then we choose $\{\bm \mu_2\}=\{(U^t)^{-1}\bm{\tilde\mu_2}\}$.

Write $\bm{\mu_1}=(\mu_{1,1}, \mu_{1,2},....,\mu_{1,g})^t$ with $\mu_{1,i}\in \mbb Z/L\mbb Z$. Our hypothesis in Proposition~\ref{jacobi-theta-main} implies that the theta components $h_{2T\bm{\mu_1}+L\bm{\mu_2}}$ of $\phi$ survive only for coset representatives of the form $ \displaystyle 2T\bm{\mu_1}+L\bm{\mu_2}\equiv (0,0,...,*)\bmod L$. But since $T$ is diagonal modulo $L$, the surviving coset representatives can be written as $2T\bm{\mu_1}+L\bm{\mu_2}$, where $\bm{\mu_1}=(0,0,...,\mu_1)^t$ and $\mu_1$ varies $\mod L$. Thus, the theta decomposition of $\phi$ can be written as 
     \begin{align}\label{phi-dec-4.1}
         \phi(\tau, z)=\sum_{\mu_1\bmod L}\sum_{\bm{\mu_2}\in\mc N_{2T}^{g,1}}h_{2T\bm{\mu_1}+L\bm{\mu_2}} \, \Theta_{LT, 2T\bm{\mu_1}+L\bm{\mu_2}}(L\tau, z).
     \end{align}  

\subsection{The related Gauss sums}
For $m, \ell\in \mbb Z$,  consider $M'=\smat{m}{-1}{1-\ell m}{\ell}\in \Gamma_1(N, L)$. That is, $\ell\equiv m\equiv 1\bmod L$ and $\ell m\equiv 1\bmod N$. Then we have that
\begin{align}
 \smat{L}{0}{0}{1}M'=\smat{m}{-L}{(1-\ell m)/L}{\ell} \smat{L}{0}{0}{1}.   
\end{align}
Let us write 
\begin{align}
 M=\smat{m}{-L}{(1-\ell m)/L}{\ell}\q \text{ and } \q M'\{\tau\}:=(1-\ell m)\tau +\ell.
\end{align}
Then using that $\phi|M'=\phi$, we can write
\begin{align}
    \phi(\tau, z) &=M'\{\tau\}^{-k}e\left(-(1-\ell m)T[z](M'\{\tau\})^{-1}\right)\phi\left(M'\tau, z(M'\{\tau\})^{-1}\right)\\
     &= M'\{\tau\}^{-k}e\left(-(1-\ell m)T[z](M'\{\tau\})^{-1}\right) \sum_{\bm\mu\in\mc N_{2LT}^{g,1}}h_{\bm\mu}(M'\tau)\Theta_{LT, \bm\mu}\left(LM'\tau, z(M'\{\tau\})^{-1}\right)\\
     &= M'\{\tau\}^{-k}e\left(-(1-\ell m)T[z](M'\{\tau\})^{-1}\right) \sum_{\bm\mu\in\mc N_{2LT}^{g,1}}h_{\bm\mu}(M'\tau)\Theta_{LT, \bm\mu}\left(M(L\tau), z(M\{L\tau\})^{-1}\right).
\end{align}

Using the transformation properties of $\Theta_{LT, \bm\mu}\left(M(L\tau), z(M\{L\tau\})^{-1}\right)$ from \eqref{theta-trans}, we get
\begin{align} \label{theta-phi-4.1}
\sum_{\bm\nu}h_{\bm\nu}(\tau) \Theta_{LT,\bm\nu}(L\tau, z)=M'\{\tau\}^{-k+1/2}\sum_{\bm\mu, \bm\nu\in \mc N_{2LT}^{g,1}}\varepsilon_{LT}(\bm\mu, \bm\nu; M) h_{\bm\mu}(M'\tau)   \Theta_{LT, \bm\nu}(L\tau,z) .
\end{align}
Next, by using the linear independence of the theta series in \eqref{theta-phi-4.1} and a change of variable $\tau\mapsto M'^{-1}\tau$, we get,
\begin{align}\label{hmutrans-sec5}
  \{M'^{-1} \tau\}^{k-1/2}  h_{\bm\nu}(M'^{-1}\tau)=\sumn_{\bm\mu\in\mc N_{2LT}^{g,1} }\varepsilon_{LT}(\bm\mu, \bm\nu;M) h_{\bm\mu}(\tau).  
\end{align}
Now let $m= 1+La$, $\ell =1+Lb$. Then we have $a+b\equiv -Lab \bmod N/L$ and we get $M=\smat{1+L a}{-L}{-( a+ b)-L a b}{1+L b}$. Further, one can easily check that,
\begin{align}
   M&= \mc S\mc T^{ b}\mc S\mc T^{-L}\mc S\mc T^{ a}\mc S; \q\text{ where } \mc T=\smat{1}{1}{0}{1} \text{ and } \mc S=\smat{0}{-1}{1}{0}.
\end{align}
We use this decomposition of $M$ and \eqref{epsilonmunu} to calculate the $\varepsilon_{LT}(\bm\mu, \bm\nu; M)$ associated with $M$. Towards this, we first recall that 
\begin{align}
\varepsilon_{LT}(\bm\mu, \bm\nu; \mc T^a)=e\left(\frac{aT^{-1}[\bm\mu]}{4L}\right)   \text{ if } \bm\nu=\bm\mu \text{ and zero otherwise.} 
\end{align}
 To see this, note that
\begin{align}
 \Theta_{LT,\bm\mu}(L\tau+ a, z)&=\sumn_{R}e\left(LT[R+\frac{(2T)^{-1}\bm\mu}{L}](L\tau+a)+2z^tLT(R+\frac{(2T)^{-1}\bm\mu}{L})\right) \\
 &= e\left(\frac{aT^{-1}[\bm\mu]}{4L}\right)  \sumn_{R}e\left(LT[R+\frac{(2T)^{-1}\bm\mu}{L}]L\tau+2z^tLT(R+\frac{(2T)^{-1}\bm\mu}{L})\right).
\end{align}
For $\mc S$, we have from \eqref{S-theta-trans} that,
\begin{align}
\varepsilon_{LT}(\bm\mu, \bm\nu; \mc S)=   e\left(-\frac{\bm\nu^tT^{-1}\bm\mu}{2L}\right).
\end{align}
Thus, we get that
\begin{align}
\varepsilon_{LT}(\bm\mu, \bm\nu;M)= \omega(M) \sumn_{\bm{\alpha}, \bm{\beta}\in \mc N_{2LT}^{g,1}} &e\left(\frac{b T^{-1}[\bm{\alpha}]}{4L}-\frac{\bm{\alpha}^tT^{-1}\bm\mu}{2L}-\frac{\bm{\beta}^tT^{-1}\bm{\alpha}}{2L}-\frac{T^{-1}[\bm{\beta}]}{4}\right)\\
&\q\q\times\sumn_{\gamma\in \mc N_{2LT}^{g,1}}e\left(\frac{aT^{-1}[\bm\gamma]}{4L}-\frac{\bm\gamma^tT^{-1}(\bm\nu+\bm{\beta})}{2L}\right).   
\end{align}
\begin{lemma}\label{lem:varepssplit}
    Let $T\equiv \mrm{diag}(t_1, ....,t_g)\mod L$ and $\varepsilon_{LT}(\bm\mu, \bm\nu;M)$ be as above. Then
\begin{align}
    \varepsilon_{LT}(\bm\mu, \bm\nu;M)= \omega(M)\varepsilon_{L}(t_g\bm\mu_1, \bm\nu_1;M)\varepsilon_{T}(L\bm\mu_2, \bm\nu_2;M),
\end{align}
where
\begin{align}
\varepsilon_{L}(t_g\bm\mu_1, \bm\nu_1;M)=\sum_{\alpha_1,\beta_1,\gamma_1\bmod L}e_L(t_g(b\alpha_1^2-2\alpha_1\mu_1-2\beta_1\alpha_1-L\beta_1^2+a\gamma_1^2-2\gamma_1(\nu_1+\beta_1))    
\end{align}
and
\begin{align}
\varepsilon_{T}(L\bm\mu_2, \bm\nu_2;M)= \sum_{\bm{\alpha}, \bm{\beta}\in \mc N_{2T}^{g,1}} &e_4\left(L(b T^{-1}[\bm{\alpha}]-2\bm{\alpha}^tT^{-1}\bm\mu_2-2\bm{\beta}^tT^{-1}\bm{\alpha}-LT^{-1}[\bm{\beta}])\right)\\
&\q\q\times\sum_{\bm\gamma\in \mc N_{2T}^{g,1}}e_4\left(L(aT^{-1}[\bm\gamma]-2\bm\gamma^tT^{-1}(\bm\nu_2+\bm{\beta}_2))\right).   
\end{align}
\end{lemma}
\begin{proof}
Let $\bm\alpha=2T\bm{\alpha}_1+L\bm\alpha_2$, $\bm\beta=2T\bm{\beta}_1+L\bm\beta_2$, then we have that
\begin{align}
e\left(-\frac{\bm\beta^tT^{-1}\bm\alpha}{2L}\right)=e_L(-\bm\beta_1^t(2T)\bm\alpha_1)e(-L\bm\beta_2^t(2T)^{-1}\bm\alpha_2).   
\end{align}
Further, noting that $\bm\alpha_1=(0,0,..\alpha_1)$, $\bm\beta_1=(0,0,..\beta_1)$ and $T\equiv \mrm{diag}(t_1, ....,t_g)\mod L$, we get that
\begin{align}
e\left(-\frac{\bm\beta^tT^{-1}\bm\alpha}{2L}\right)=e_L(-2t_g\beta_1\alpha_1)e(-L\bm\beta_2^t(2T)^{-1}\bm\alpha_2).   
\end{align}
Similarly,
\begin{align}
e\left(\frac{T^{-1}[\bm\alpha]}{4L}\right)=e_L(t_g\alpha_1^2)e_4(LT^{-1}[\bm\alpha_2]).   
\end{align}
We get the required decomposition for $\varepsilon_{LT}(\bm\mu, \bm\nu;M)$ by writing  $\bm\mu=2T\bm\mu_1+L\bm\mu_2$ and $\bm\nu=2T\bm\nu_1+L\bm\nu_2$. This completes the proof of the lemma.
\end{proof}
\begin{remark}
 This nested summation is not really a product of Gauss sums, whence we can't directly invoke known evaluations to immediately understand the coefficients $\varepsilon_{LT}(\bm\mu, \bm\nu;M)$. The quantities $\varepsilon_{LT}(...)$ will ultimately lead us to a matrix whose rank will play an important role in proving Theorem~\ref{diagonal-sqrfree}. The next lemma is thus devoted to understanding these coefficients.
\end{remark}
\begin{lemma}\label{lem:epsL}
    Let $\varepsilon_{L}(t_g\bm\mu_1, \bm\nu_1;M)$ be as above. For any $(a,L)=1$, there exists a  non-zero constant $A_M:=A_M(\bm\nu, L)$ such that
    \begin{align}
     \varepsilon_{L}(t_g\bm\mu_1, \bm\nu_1;M)=A_M  \cdot G_1(t_g(b-\overline{\ell}), -2t_g(\mu_1+\overline{\ell a}\nu_1), L), 
    \end{align}
where $\ell= -(L+\overline{a})$, $a\overline{a}\equiv 1\bmod L$, $\ell\overline{\ell}\equiv 1\bmod L$ and $a+b\equiv -Lab\bmod N/L$.
\end{lemma}
\begin{proof}
    Recall 
    \begin{align}
\varepsilon_{L}(t_g\bm\mu_1, \bm\nu_1;M)=\sum_{\alpha_1,\beta_1,\gamma_1\bmod L}e_L(t_g(b\alpha_1^2-2\alpha_1\mu_1-2\beta_1\alpha_1-L\beta_1^2+a\gamma_1^2-2\gamma_1(\nu_1+\beta_1)).    
\end{align}
Thus sum over $\gamma_1$ is nothing but the one-dimensional Gauss sum $G_1(t_ga, -2t_g(\nu_1+\beta_1), L)$. By choosing $(a, L)=1$, $G_1(t_ga, -2t_g(\nu_1+\beta_1), L)\neq 0$ and the one-dimensional Gauss sum simplifies to
\begin{align}
    G_1(t_ga, -2t_g(\nu_1+\beta_1), L)= e_L(-t_g\overline{a}(\nu_1+\beta_1)^2)G_1(t_ga, 0, L).
\end{align}
Thus the sum over $\beta_1$ reduces to
\begin{align}
    &\sum_{\beta_1\bmod L}e_L(t_g(-L\beta_1^2-2\beta_1\alpha_1)e_L(-t_g\overline{a}(\nu_1+\beta_1)^2)\\
    &= e_L(-t_g\overline{a}\nu_1^2)\sum_{\beta_1\bmod L}e_L(t_g(-(L+\overline{a})\beta_1^2-2\beta_1(\alpha_1+4\overline{a}\nu_1))\\
    &=e_L(-t_g\overline{a}\nu_1^2) G_1(-t_g(L+4\overline{a})), -2t_g(\alpha_1+\overline{a}\nu_1), L).
\end{align}
Let us put $\ell=-(L+\overline{a})$. Since $(a, L)=1$, we see that $(\ell, L)=1$. Thus $G_1(t_g\ell, -2t_g(\alpha_1+\overline{a}\nu_1), L)\neq 0$ and we have
\begin{align}
G_1(t_g\ell, -2t_g(\alpha_1+\overline{a}\nu_1), L)= e_L(-t_g\overline{\ell} (\alpha_1+\overline{a}\nu_1)^2) G_1(t_g\ell, 0, L).   
\end{align}
Put $A_M'(\bm \nu, L):=e_L(-t_g\overline{a}\nu_1^2)G_1(t_ga, 0, L) G_1(t_g\ell, 0, L)$. Then we can write
\begin{align}
 \varepsilon_{L}(t_g\bm\mu_1, \bm\nu_1;M)&=  A_M'(\bm \nu, L) \sum_{\alpha_1\bmod L}  e_L(t_g(b\alpha_1^2-2\alpha_1\mu_1)e_L(-t_g\overline{\ell} (\alpha_1+\overline{a}\nu_1)^2)\\
&=A_M'(\bm \nu, L)e_L(-t_g\overline{\ell}\overline{a}^2\nu_1^2) G_1(t_g(b-\overline{\ell}), -2t_g(\mu_1+\overline{\ell a}\nu_1), L).
\end{align}
We get the required expression by defining $A_M(\bm\nu, L):=A_M'(\bm \nu, L)e_L(-t_g\overline{\ell}\overline{a}^2\nu_1^2) $.
\end{proof}
\begin{lemma}\label{lem:epsT}
    Let $\varepsilon_{T}(\bm\mu, \bm\nu;M)$ be as above and $a$ be any integer such that $(a, d(T))=1$. In addition, let $a\equiv 0\bmod 4$ if $g$ is odd. Then there exist non-zero constants $C_M:=C_M(\bm\nu, d(T))$ and $D_M:=D_M(\bm\nu, d(T))$ such that
    \begin{align}
     \varepsilon_{T}(L\bm\mu, \bm\nu;M)=\begin{cases}
       C_M \cdot G_1(Lt(b-\overline{\ell}),2Lt ((\overline{\ell}L+1)\nu-\mu), 4d(T)) &\text{ when } g \text{ is odd;}\\
        D_M \cdot G_1(Lt(b-\overline{\ell}),2Lt (-(\overline{a\ell}+1)\nu-\mu), d(T))&\text{ when } g \text{ is even,}
     \end{cases}   
    \end{align}
where $(t, \det(2T))=1$, $\ell=-( L+\overline{a})$, $(\ell, d(T))=1$, $a\overline{a}\equiv 1\bmod d(T)$, $\ell\overline{\ell}\equiv 1\bmod d(T)$  and $a+b\equiv -Lab \bmod N/L$.
\end{lemma}
\begin{proof}
Let $U\in\slnz$ and $\widetilde{T}=$ be as in Section \ref{sec:HDGauss}. That is, $\widetilde{T}=(2T)[U]$. We choose the coset representatives $\mc N_{2T}^{g,1}$ as in Section \ref{sec:HDGauss}. Further, let $\widetilde{T}^*=(t_{i,j}^*)$ denote the adjoint of $\widetilde{T}.$ For any $\tilde{\bm{\alpha}}=(0,..,\alpha)$ and $\tilde{\bm{\beta}}=(0,..,\beta)$, we have 
\begin{align}
    \bm\alpha (2T)^{-1}\bm\beta=\tilde{\bm{\alpha}}\widetilde{T}^{-1}\tilde{\bm{\beta}}= \frac{1}{\det(2T)}t_{g,g}^*\alpha\beta.
\end{align}
Thus we get $\varepsilon_{T}(\bm\mu, \bm\nu;M)= \varepsilon_{T}(\mu, \nu;M)$ and
\begin{align} \label{eLT-1}
\varepsilon_{T}(\mu, \nu;M)= \omega(M)\sum_{\alpha, \beta,\gamma}e_{2\det(2T)}\left(Lt_{g,g}^*(b\alpha^2-2\alpha\mu-2\beta\alpha-L\beta^2+a\gamma^2-2\gamma(\nu+\beta)\right),    
\end{align}
where $\alpha, \beta,\gamma$ vary $\mod \det(2T)$. We also recall that $(t_{g, g}^*, \det(2T))=1$.  Note that we have used the fact that $\det(2T)=2 d(T)$ if $g$ is odd and $\det(2T)= d(T)$ if $g$ is even, which allows us to write \eqref{eLT-1} uniformly in terms of $\det(2T)$. To simplify the notations, for the rest of the section, we use $d:=d(T)$ and $t:=t^*_{g,g}$.

Next, consider the sum over $\gamma$. It is nothing but the one-dimensional Gauss sum given by
\begin{align}
    \begin{cases}
        \frac{1}{2}G_1(Lta, -2Lt(\nu+\beta), 4d) &\text{ if } g \text{ is odd;}\\
        G_1(Lta, -2Lt(\nu+\beta), d) &\text{ if } g \text{ is even}.
    \end{cases}
\end{align}
\textbf{ $g$ odd:} In this case, we choose $a=4a_1$ with $(a_1, d)=1$.  Then noting that $G_1(a,b,c)=0$ if $ (a,c)\nmid b$, the Gauss sum $G(4Lta_1, -2Lt(\nu+\beta), 4d)$ survives only when $2|(\nu+\beta)$ and in that case we have  (see \cite{berndt1998gauss})
\begin{align}
G(4Lta_1, -2Lt(\nu+\beta), 4d)&=4G(Lta_1, -L(\nu+\beta)/2, d)=4e_{d}(-Lt\overline{16a_1}(\nu+\beta)^2)G_1(Lta_1, 0, d).
\end{align}
Here $16a_1(\overline{16a_1})\equiv 1\bmod d$. Thus we get
\begin{align}
\varepsilon_{T}(\mu, \nu;M)&=4G_1(Lta_1, 0, d)\sum_\alpha\underset{\beta\equiv -\nu\bmod 2}{\sum_{\beta\bmod 2d}}e_{4d}\left(Lt(b\alpha^2-2\alpha\mu-2\beta\alpha-L\beta^2)\right)e_{d}(-Lt\overline{16a_1}(\nu+\beta)^2).
\end{align}
Putting $\beta= -\nu+2\delta$, the sum over $\beta$ transforms into
\begin{align}
&e_{4d}(2Lt\alpha\nu)\sum_{\delta_g\bmod d}e_{4d}\left(Lt(-4\delta\alpha-L(-\nu+2\delta)^2)\right)e_{d}\left(-Lt\overline{4a_1}\delta^2\right) \\
&=e_{4d}(2Lt\alpha\nu)e_{4d}(-L^2t\nu^2)\sum_{\delta\bmod d}e_{d}\left(Lt(-(L+\overline{a})\delta^2+(L\nu-\alpha)\delta)\right)\\
&=e_{4d}(2Lt\alpha\nu)e_{4d}(-L^2t\nu^2) G_1(-Lt(L+\overline{a}), Lt(L\nu-\alpha), d). 
\end{align}
We choose $a$ such that $(L+\overline{a}, d)=1$ and thus $G_1(-Lt(L+\overline{a}), Lt(L\nu-\alpha), d)\neq 0$. Moreover, writing $\ell=-( L+\overline{a})$, we have
\begin{align}
G_1(Lt\ell, Lt(L\nu-\alpha), d)= e_{d}(-Lt\overline{4\ell}(L\nu-\alpha)^2)G_1(Lt\ell, 0, d). 
\end{align}
Put $C_M(\bm\nu, d):= 4G_1(Lta_1, 0, d)e_{4d}(-L^2t\nu^2)G_1(Lt\ell, 0, d)$. Thus we get    
\begin{align}
\varepsilon_{T}(\mu, \nu;M)&=C_M(\bm\nu, d)\sum_{\alpha\bmod 2d}e_{4d}\left(Lt(b\alpha^2-2\alpha(\mu-\nu)-\overline{\ell}(L\nu-\alpha)^2)\right)\\
&=C_M(\bm\nu, d)G_1(Lt(b-\overline{\ell}), 2Lt ((\overline{\ell}L+1)\nu-\mu), 4d).
\end{align}

\textbf{$g$ even:} In this case, $t=t_{g,g}^*$ is even (see \cite{boche-das}) and we choose $a$ such that $(a, d)=1$.  Then the Gauss sum $G_1(Lta, -2Lt(\nu+\beta), d) $ is given by  (see \cite{berndt1998gauss})
\begin{align}
G(Lta, -2Lt(\nu+\beta), d)= e_{d}(-Lt\overline{a}(\nu+\beta)^2)G_1(Lta,0, d).
\end{align}
Thus we get
\begin{align}
\varepsilon_{T}(\mu, \nu;M)&=G_1(Lta,0, d)\sum_\alpha\sum_{\beta\bmod d}e_{d}\left(Lt(b\alpha^2-2\alpha\mu-2\beta\alpha-L\beta^2-\overline{a}(\nu+\beta)^2)\right).
\end{align}
The sum over $\beta$ can be simplified as below.
\begin{align}
&=\sum_{\beta\bmod d}e_{d}\left(Lt(-2\beta\alpha-L\beta^2-\overline{a}(\nu+\beta)^2)\right)\\
&=e_{d}(-Lt\overline{a}\nu^2)\sum_{\beta\bmod d}e_{d}\left(Lt(-(L+\overline{a})\beta^2-2(\alpha+\overline{a}\nu)\beta)\right)\\
&=e_{d}(-Lt\overline{a}\nu^2) G_1(-tL(L+\overline{a}), -2Lt(\alpha+\overline{a}\nu), d).
\end{align}
We again choose $a$ such that $(L+\overline{a}, d)=1$ as before, and thus $ G_1(-Lt(L+\overline{a}), -2Lt(\alpha+\overline{a}\nu), d)\neq 0$. Moreover, writing $\ell=-( L+\overline{a})$, we have
\begin{align}
G_1(Lt\ell, -2Lt(\alpha+\overline{a}\nu), d)= e_{d}(-Lt\overline{\ell}(\alpha+\overline{a}\nu)^2)G_1(Lt\ell, 0, d). 
\end{align}
Put $D'_M(\bm\nu, d):=e_{d}(-Lt\overline{a}\nu^2)G_1(Lta,0, d)G_1(Lt\ell, 0, d)$ Thus $\varepsilon_{T}(\mu, \nu;M)$ can be written as
\begin{align}
&=D'_M(\bm\nu, d)\sum_{\alpha\bmod d}e_{d}\left(Lt(b\alpha^2-2\alpha\mu-\overline{\ell}(\alpha+\overline{a}\nu)^2)\right)\\
&=D'_M(\bm\nu, d)e_{d}(-Lt(\overline{\ell}\overline{a}^2\nu^2))G_1(Lt(b-\overline{\ell}),-2Lt ((\overline{a\ell}+1)\nu+\mu), d).\qedhere
\end{align}
We get the required expression by putting $D_M(\bm\nu, d):=D'_M(\bm\nu, d)e_{d}(-Lt(\overline{\ell}\overline{a}^2\nu^2))$.
\end{proof}
\subsection{Proof of the Proposition~\ref{jacobi-theta-main}}\label{sec:thetamain}
Recall from \eqref{phi-dec-4.1}, the theta decomposition of $\phi$.
\begin{align}
    \phi(\tau, z)=\sum_{\mu_1\bmod L}\sum_{\bm{\mu_2}\in\mc N_{2T}^{g,1}}h_{2T\bm{\mu_1}+L\bm{\mu_2}} \, \Theta_{LT, 2T\bm{\mu_1}+L\bm{\mu_2}}(L\tau, z).
\end{align}
Assume that the theta components $h_{2T\bm{\mu_1}+L\bm{\mu_2}}$ are zero when $\bm{\mu_1}\equiv (0,0,...0)\bmod L$. That is, $\mu_1=0$. In other words, $h_{L\bm{\mu_2}}=0$ for all $\bm{\mu_2}\in \mc N_{2T}^{g,1}$. That is, we get (the other terms being zero anyway)
\begin{align}\label{lin-eq-1}
\sum_{\gamma\bmod L}\sum_{\bm{\delta}\in\mc N_{2T}^{g,1}}\varepsilon_{LT}(2T{\bm{ \gamma}}+L\bm\delta, L\bm{\mu_2};  M') h_{2T{\bm{ \gamma}}+L\bm{\delta}}(\tau)=0,
\end{align} 
%where $\bm{ \tilde \gamma}=(0,\ldots, \gamma) \in \z^{1,g}$.

 We remind the reader that in \eqref{lin-eq-1}, $\bm \gamma=(0.0,\ldots, \gamma)$.
Our aim is to prove that $h_{2T{\bm{ \gamma}}+L\bm{\delta}}=0$ for all $\gamma\mod L$ and all $\bm{\delta}\in \mc N_{2T}^{g,1}$. If $\gamma\equiv 0\bmod L$, then $h_{2T{\bm{ \gamma}}+L\bm{\delta}}=0$ by assumption. Thus, it is enough to consider the case when $\gamma\not \equiv 0 \mod L$.

Using the transformation properties of $h_{2T\bm{\gamma}+L\bm\delta}$ under $\tau \mapsto \tau+t$, we have
\begin{align}
\sum_{\gamma\bmod L}\sum_{\bm{\delta}\in\mc N_{2T}^{g,1}}\varepsilon_{LT}(2T{\bm{ \gamma}}+L\bm\delta, L\bm{\mu_2};  M) e(T^{-1}[L\bm\delta/2]t)h_{2T{\bm{ \gamma}}+L\bm{\delta}}(\tau)=0.
\end{align} 
The linear independence of the additive character $t \mapsto e(v \cdot t)$ ($t \in \z$, $v \in \Q$)shows that we have an equation of the form
\begin{align}
   \sum_{\gamma\bmod L} \underset{T^{-1}[L\bm\delta] - T^{-1}[L\bm\delta_0] \in 4 \z }{\sum_{\bm\delta\in\mc N_{2T}^{g,1}}}\varepsilon_{LT}(2T{\bm{ \gamma}}+L\bm\delta, L\bm{\mu_2};  M)h_{2T{\bm{ \gamma}}+L\bm{\delta}}(\tau)=0.
\end{align}
Before proceeding, we make the following notational simplifications.
\begin{align}
d=d(T);\q d_T=\begin{cases}
    4d & \text{ if } g \text{ is odd;}\\
    d & \text{ if } g \text{ is even.}
\end{cases} \q d_T^*=\begin{cases}
    2d & \text{ if } g \text{ is odd;}\\
    d & \text{ if } g \text{ is even.}
\end{cases}
\end{align}
Next, we note that, since $(L, d)=1$,  
\begin{align}\label{squareclass-vector}
T^{-1}[L\bm\delta] - T^{-1}[L\bm\delta_0] \in 4 \z \iff T^{-1}[\bm\delta] - T^{-1}[\bm\delta_0] \in 4 \z .
\end{align}
Recall the choice of coset representatives for $\mc N_{2T}^{g,1}$: $\{\bm\delta\}=\{(0,0,...,\delta)^t: \delta\bmod d_T^*\}$. Thus, the square-class in \eqref{squareclass-vector} is equivalent to
\begin{align}
    {\delta}^2\equiv {\delta}_{0}^2 \bmod d_T,
\end{align}
Thus, we are left with
\begin{align}\label{lin-eq-sqrcls}
    \sum_{\gamma\bmod L} \underset{{\delta}^2\equiv {\delta}_{0}^2 \bmod d_T}{\sum_{\delta\bmod d_T^*}}\varepsilon_{LT}(2T{\bm{ \gamma}}+L\bm\delta, L\bm{\mu_2};  M)h_{2T{\bm{ \gamma}}+L\bm{\delta}}(\tau)=0.
\end{align}
We split the proof into two cases depending on the nature of the square classes above. In the first case, we deal with those square-classes where $\delta_0\neq 0 \bmod d$. In the second case, we deal with the rest of the square-classes (i.e., $\delta_0\equiv 0\bmod d$).
\subsubsection{Square classes that are non-zero modulo $d$:}\label{sec:non-zerosqrcls}
We first choose $a$ and $b$ such that the conditions in Lemma \ref{lem:epsL} and \ref{lem:epsT} are satisfied. That is,
\begin{align}
(a, Ld)=1,\, a\equiv 0\bmod 4 \,(\text{ if } g \text{ is odd}),\, a+b\equiv -Lab\bmod N/L,\, \text{ and } (\ell, Ld)=1.   
\end{align}
Then from  Lemma \ref{lem:epsL} and \ref{lem:epsT},  $\varepsilon_{LT}(2T{\bm{ \gamma}}+L\bm\delta, L\bm{\mu_2};  M)$ can be written as
\begin{align}
&= \omega(M)\varepsilon_L(2t_g\bm\gamma, 0,; M)\varepsilon_T(L\bm\delta, \bm{\mu_2}; M)\\
&=  E_{M}G_1(t_g(b-\overline{\ell}), -2t_g\gamma, L)G_1(Lt(b-\overline{\ell}),2Lt (\ell_1\mu_2-\delta), d_T),
\end{align}
where $E_{M}$ is non zero constant. We define
\begin{align}
    \ell_1:=\begin{cases}
        \overline{\ell}L+1& \text { if } g \text{ is odd};\\
        -(\overline{a\ell}+1) & \text { if } g \text{ is even}.
    \end{cases}
\end{align}
We also note here that the inverses, $a\overline{a}\equiv \ell\overline{\ell}\equiv 1\bmod Ld$ (see Lemma \ref{lem:epsL} and \ref{lem:epsT}). Recall that $(d_T, L)=(2t_g, L)=(t, d_T)=1$. Let $B_1, B_2$ be any two integers.  We choose $b$ such that 
\begin{align}
b-\overline{\ell}\equiv B_1\bmod L \q\text{and}\q
b-\overline{\ell}\equiv B_2\bmod d_T.
\end{align}
Then
\begin{align} \label{modL}
G_1(t_g(b-\overline{\ell}),-2t_g \gamma, L)&= G(t_gB_1, -2t_g\gamma, L)\\
G_1(Lt(b-\overline{\ell}),2L t(\ell_1\mu_2-\delta), d_T)&=G_1(LtB_2, 2Lt (\ell_1\mu_2-\delta), d_T).\label{mod2M}
\end{align}
By choosing $(B_1, L)=1$, the Gauss sum $G(t_gB_1, -2t_g\gamma, L)$ above is non-zero. Thus, in \eqref{lin-eq-sqrcls} we are left with
\begin{align}\label{gammaA_gamma}
    \sum_{\gamma\mod L}G(t_gB_1, -2t_g\gamma, L)A(\gamma; \mu_2, B_2, 2T)=0,
\end{align}
where
\begin{align}
  A(\gamma; \mu_2, B_2, 2T)=  \underset{{\delta}^2\equiv {\delta}_{0}^2 \bmod d_T}{\sum_{\delta\bmod d_T^*}} G_1(LtB_2, 2Lt (\ell_1\mu_2-\delta), d_T) h_{2T{\bm{ \gamma}}+L\bm{\delta}}(\tau).
\end{align}
\textbf{Claim 1:} $A(\gamma; \mu_2, B_2, 2T)+(-1)^kA(\gamma; -\mu_2, B_2, 2T)=0$ for all $\gamma\bmod L$. 

\textit{Proof of the claim:} We first note that 
\begin{align}
G(t_gB_1, -2t_g\gamma, L)&=\sum_{x\bmod L}e_L(t_gB_1x^2-2t_g\gamma)=e\left(-\frac{t_g\overline{B_1}\gamma^2}{L}\right)G(t_gB_1, 0, L).
\end{align}
We have that $G(t_gB_1, 0, L)\neq 0$. Let us put $r=-t_g\overline{B_1} $. As $B_1$ varies such that $(B_1, L)=1$, $r$ also varies over coprime residue classes $\mod L$.  Thus, we are left with
\begin{align}
\sum_{\gamma\mod L}e\left(\frac{r\gamma^2}{L} \right)A(\gamma; \mu_2, B_2, 2T)=0 .   
\end{align}
Now note that $e_{L}\left(r\gamma^2 \right)$ is invariant under $\gamma\mapsto -\gamma$.  Further, note that as $\gamma\mapsto -\gamma$, $2T{\bm{ \gamma}}+L\bm{\delta}\mapsto -2T\bm{\gamma}+L\bm\delta$. Thus we get
\begin{align}
A(-\gamma; \mu_2, B_2, 2T)=   \underset{{\delta}^2\equiv {\delta}_{0}^2 \bmod d_T}{\sum_{\delta\bmod d_T^*}}G_1(LtB_2, 2Lt (\ell_1\mu_2-\delta), d_T) h_{-2T\bm{\gamma}+L\bm\delta}(\tau).  
\end{align}
Using  $c_\phi(n,-r)=(-1)^kc_\phi(n,r)$, we have that $h_{-2T\bm{\gamma}+L\bm\delta}(\tau)=(-1)^k h_{2T\bm{\gamma}-L\bm\delta}(\tau)$. Thus, we get
\begin{align}
 A(-\gamma; \mu_2, B_2, 2T)&= (-1)^k  \underset{{\delta}^2\equiv {\delta}_{0}^2 \bmod d_T}{\sum_{\delta\bmod d_T^*}} G_1(LtB_2, 2Lt (\ell_1\mu_2-\delta), d_T)h_{2T\bm{\gamma}-L\bm\delta}(\tau).   
\end{align}
Next, making a change of variable $\delta\mapsto -\delta$ and from the fact that $G_1(a,-b,c)=G_1(a,b,c)$,  we get
\begin{align}
A(-\gamma; \mu_2, B_2, 2T)&= (-1)^k   \underset{{\delta}^2\equiv {\delta}_{0}^2 \bmod d_T}{\sum_{\delta\bmod d_T^*}} G_1(LtB_2, 2L t(\ell_1\mu_2+\delta), d_T)h_{2T\bm{\gamma}+L\bm\delta}(\tau)\\
&= (-1)^k A(\gamma; -\mu_2, B_2, 2T).\label{AgammaA-gamma}
\end{align}
Thus, the equation in \eqref{gammaA_gamma} now reduces to the following equation.
\begin{align}\label{gammaA-gamma2}
\sum_{\gamma\mod L}e\left(\frac{r\gamma^2}{L} \right)\left(A(\gamma; \mu_2, B_2, 2T)+(-1)^kA(\gamma; -\mu_2, B_2, 2T) \right)=0 .   
\end{align}
Making a change of variable $\gamma\mapsto -\gamma$, we see that
\begin{align}
\sum_{\gamma\mod L}e\left(\frac{r\gamma^2}{L} \right)\left(A(-\gamma; \mu_2, B_2, 2T)+(-1)^kA(-\gamma; -\mu_2, B_2, 2T) \right)&=0\\
\implies\sum_{\gamma\mod L}e\left(\frac{r\gamma^2}{L} \right)\left((-1)^kA(\gamma; -\mu_2, B_2, 2T)+A(\gamma; \mu_2, B_2, 2T) \right)=0.
\end{align}
Thus the equation \eqref{gammaA-gamma2} is invariant w.r.t $\gamma\mapsto -\gamma$. Consequently, it is enough to consider the sum over $1\le \gamma\le (L-1)/2$ (when $\gamma=0$, $h_{2T\gamma+L\delta}=0$ and thus $A(\gamma; \mu_2, B_2, 2T)=0$).

Recall that $L$ is odd and square-free. Let $t$ denote the number of prime divisors of $L$. To prove \textbf{Claim 1}, we argue by induction on $t$.

\textit{The base case:} In this case, $L=q$, a prime. For $1\le r\le (q-1)/2$ and $1\le \gamma\le (L-1)/2$, consider the $(q-1)/2\times (q-1)/2$ matrix
\begin{align}
    \mc M = \left(e_{q}\left(r\gamma^2 \right) \right).
\end{align}
Since the quantities  $e_{q}\left(s\gamma_1^2 \right)$ are pairwise distinct, this is simply a Vandermonde matrix and is thus invertible. This proves the base case of induction. 

\textit{Concluding the Induction:} Now assume the claim to be true for all square-free integers $L'$ with $<t$ prime factors and let $L=L'q$. Let $r=L'r_1+qr_2$ and $\gamma=L'\gamma'+q\gamma''$. Then we see that
\begin{align}\label{Ind-Agamma}
    \sum_{\gamma' \bmod q}e\left(\frac{r_1L'^2\gamma'^2}{q} \right)B(\gamma'; \mu_2, B_2, 2T)=0,
\end{align}
where
\begin{align}
    B(\gamma'; \mu_2, B_2, 2T):=\left( \sum_{\gamma'' \bmod L'}e\left(\frac{r_2q^2\gamma''^2}{L'} \right)\left(A(\gamma; \mu_2, B_2, 2T)+(-1)^kA(\gamma; -\mu_2, B_2, 2T) \right)\right)
\end{align}
Making the change of variables $\gamma''\mapsto -\gamma''$ and using \eqref{AgammaA-gamma}, we see that, $B(-\gamma'; \mu_2, B_2, 2T)$ is given by
\begin{align}
&\sum_{\gamma'' \bmod L'}e\left(\frac{r_2q^2\gamma''^2}{L'} \right)\left(A(-L'\gamma'-q\gamma''; \mu_2, B_2, 2T)+(-1)^kA(-L'\gamma'-q\gamma''; -\mu_2, B_2, 2T) \right)\\
&= \sum_{\gamma'' \bmod L'}e\left(\frac{r_2q^2\gamma''^2}{L'} \right)\left((-1)^kA(L'\gamma'+q\gamma''; -\mu_2, B_2, 2T)+A(L'\gamma'+q\gamma''; \mu_2, B_2, 2T) \right)\\
&=B(\gamma'; \mu_2, B_2, 2T).
\end{align}
Thus we can conclude that \eqref{Ind-Agamma} is invariant w.r.t $\gamma'\mapsto -\gamma'$. Consequently, the sum in \eqref{Ind-Agamma} can be assumed to be over $0\le \gamma'\le (q-1)/2$.

For $0\le \gamma'\le (q-1)/2$ and $1\le r_1\le (q+1)/2$ (recall that we can freely vary $r_1$ over co-prime residue classes modulo $q$), 
consider the following $(q+1)/2\times (q+1)/2$ Vandermonde matrix.
\begin{align}
    \mc M := \left(e_{q}\left(r_1L'^2\gamma'^2 \right) \right).
\end{align}
Since $\mc M$ is invertible, $B(\gamma'; \mu_2, B_2, 2T)=0$ for all $\gamma'\bmod q$. The assertion in \textbf{Claim 1} now follows by induction.

Continuing with the proof of the proposition, we are left with
\begin{align}\label{only2Msum}
\underset{{\delta}^2\equiv {\delta}_{0}^2 \bmod d_T}{\sum_{\delta\bmod d_T^*}}  \left(G_1(LtB_2, 2Lt (\ell_1\mu_2-\delta), d_T) +(-1)^kG_1(LtB_2, 2Lt (-\ell_1\mu_2-\delta), d_T) \right)h_{2T\bm{\gamma}+L\bm\delta}(\tau) =0.  
\end{align}
By choosing $B_2$ such that $(B_2, d_T)=1$, we see that
\begin{align}
G_1(LtB_2, 2Lt(\ell_1\mu_2\pm\delta), d_T) = e\left(-\frac{Lt\overline{B_2}(\ell_1\mu_2\pm\delta)^2}{d_T} \right)G_1(LtB_2, 0, d_T).    
\end{align}
Since $G_1(LtB_2, 0, d_T)\neq 0$, we are left with
\begin{align}
\underset{{\delta}^2\equiv {\delta}_{0}^2 \bmod d_T}{\sum_{\delta\bmod d_T^*}} \left(e\left(-\frac{Lt\overline{B_2}(\ell_1\mu_2-\delta)^2}{d_T}  \right)+(-1)^ke\left(-\frac{Lt\overline{B_2}(\ell_1\mu_2+\delta)^2}{d_T} \right) \right)   h_{2T\bm{\gamma}+L\bm\delta}(\tau)=0.
\end{align}
\textbf{Reduction to  $1$-variable:}
For any fixed $\delta_{0}\in \mbb Z$ with $\delta_0\not\equiv 0 \bmod d(T)$, let us define the set $S_{\delta_{0}}:=\{\delta \bmod d_T^*:{\delta}^2\equiv {\delta}_{0}^2 \bmod d_T\}$. 
As $\bm{\mu_2}$ varies over elements in $\mc N_{2T}^{g,1}$ and $\delta\in S_{\delta_0}$, let us define the matrix
\begin{align}
    \mc M_{\delta_{0}}:=\left(e\left(-\frac{Lt\overline{B_2}(\ell_1\mu_2-\delta)^2}{d_T}  \right)+(-1)^ke\left(-\frac{Lt\overline{B_2}(\ell_1\mu_2+\delta)^2}{d_T} \right)\right),
\end{align}
We see that $\mc M_{\delta_0}$ has the same rank as the matrix
\begin{align}
   \mc S_{\delta_{0}} :=\left( e\left(\frac{2\ell_1Lt\overline{B_2}\delta\mu_2}{d_T} \right)+(-1)^k \right).
\end{align}
\textbf{Claim 2:} $\text{rank}(\mc S_{\delta_{0}})=|S_{\delta_{0}}|$. That is, the rank of $\mc S_{\delta_{0}}$ is maximal. 

\textit{Proof of the claim: }We first recall the following conditions for choosing $a$, $b$ and the various parameters involved in the previous calculations.
\begin{align}
&(a, Ld)=1,\, a\equiv 0\bmod 4 \text{ if } g \text{ is odd},\,a+b\equiv -Lab\bmod N/L,\\
&\ell=-(L+\overline{a}),\, (\ell, d)=1, \, b-\overline{\ell}\equiv B_2\bmod d_T,\q \ell_1=\begin{cases}
        \overline{\ell}L+1& \text { if } g \text{ is odd};\\
        -(\overline{a\ell}+1) & \text { if } g \text{ is even}.
    \end{cases}
\end{align}

We would need to choose $a$ and $b$ such that $(\ell, d)=(\ell_1, d)=1$. Towards this, let $d=p_1p_2..p_t$ and let $L\equiv L_j\bmod p_j$. We choose $a$ such that $\overline{a}\equiv -L_j+1\bmod p_j$ if $p_j\nmid -L_j+1$. Otherwise, we choose $\overline{a}\equiv -L_j-1\bmod p_j$. Suppose $p_j\nmid -L_j+1$ for some $j$, then
\begin{align}
   \ell= -(L+\overline{a})\equiv -1\bmod p_j; \q \ell_1=\begin{cases}
        -L_j+1\bmod p_j& \text { if } g \text{ is odd};\\
        -(-\overline{a}+1)\bmod p_j & \text { if } g \text{ is even}.
    \end{cases}
\end{align}
Thus, $(\ell, p_j)=1$ and since $p_j\nmid -L_j+1$, $(\ell_1, p_j)=1$ when $g$ is odd. Further, $-\overline{a}+1\equiv L_j\bmod p_j $, thus $(\ell_1, p_j)=1$ in $g$ even case also.

Similarly, if $p_j|-L_j+1$ for some $j$, then we have that $p_j\nmid -L_j-1$. In this case, $\overline{a}\equiv -L_j-1\bmod p_j$. Thus
\begin{align}
   \ell= -(L+\overline{a})\equiv 1\bmod p_j;\q \ell_1=\begin{cases}
        -L_j+1\bmod p_j& \text { if } g \text{ is odd};\\
        -(\overline{a}+1)\equiv L_j\bmod p_j & \text { if } g \text{ is even}.
    \end{cases}
\end{align}
Thus we get that $(\ell, p_j)=(\ell_1, p_j)=1$ in this case too. 

 The choice of $b$ can now be made using CRT as below.
\begin{align}
    b&\equiv -a(\overline{aL+1})\bmod N/L\\
    b&\equiv B_2+\overline{\ell}\bmod d_T.
\end{align}

In conclusion, we have chosen $a$ and $b$ such that $(\ell, d)=(\ell_1, d)=1$. We also choose $(B_2, d)=1$. For simplicity, we write 
\begin{align}
d^*=d_T^*,\, B:=
\begin{cases}
\ell_1Lt\overline{B_2} & \text{ if } g \text{ is odd;}\\ 2\ell_1Lt\overline{B_2} & \text{ if } g \text{ is even}.
\end{cases}
\end{align}
Since $(B_2,d)=1$, we see that $(B, d)=1$. Let $\zeta_{d^*} = e^{2\pi i / d^*}$. The entries of $\mc S_{\delta_{0}}$ can be written as
\begin{align}
    s_{\delta, \mu_2}:=\zeta_{d^*}^{B\delta\mu_2}+(-1)^k.
\end{align}
Thus the rows of the matrix $\mc S_{\delta_{0}}$, denoted by $R_\delta$ for $\delta \in S_{\delta_{0}}$, are given by:
\[
R_\delta = \chi_{B\delta} + (-1)^k\chi_0,
\]
where $\chi_n(\mu) = e^{2\pi i n \mu / d^*}$ are the additive characters of $\mathbb{Z}/d^*\mathbb{Z}$.

We aim to show that the rows $R_\delta$ are linearly independent. Towards this, it suffices to show that the indices $n_\delta = B\delta \bmod{d^*}$ are distinct for distinct $\delta \in S_{\delta_{0}}$. That is, we must show that the map 
\[ \varphi: S_{\delta_{0}} \to \mathbb{Z}/d^*\mathbb{Z} \text{  \q defined by \q  } \varphi(\delta) = B\delta\] 
is injective. We split the proof of the claim into two cases, depending on the parity of $g$. 

\textbf{$g$ odd}:
Suppose $\varphi(\delta_1) = \varphi(\delta_2)$ for some $\delta_1, \delta_2 \in S_{\delta_{0}}$. That is,
$B \delta_1 \equiv B \delta_2 \bmod{2d}.$ Since $B$ is coprime to $d$, we get
\[
\delta_1 \equiv \delta_2 \bmod d.
\]
In $\mathbb{Z}/2d\mathbb{Z}$, this congruence implies that either $\delta_1 = \delta_2$ or $\delta_2 = \delta_1 + d$.

We claim that $\delta_1$ and $\delta_1 + d$ cannot be simultaneously in the square-class of $\delta_0$ $\mod 4d$.  Suppose $\delta_1 \in S_{\delta_0}$, then
\[
(\delta_1 + d)^2 = \delta_1^2 + 2d\delta_1 + d^2.
\]
If $\delta_1 + d\in S_{\delta_0}$, then $(\delta_1 + d)^2 \equiv \delta_0^2\equiv\delta_1^2 \pmod{4d}$. Thus
\begin{align}
2d\delta_1 + d^2 \equiv 0 \pmod{4d}\implies 2\delta_1 + d \equiv 0 \pmod 4,
\end{align}
which is not possible since $d$ is odd. Thus $\delta_2\not\in S_{\delta_0}$ and we get that $\varphi$ is injective.

Since the map $\delta \mapsto B\delta \pmod{2d}$ is injective on $S_{\delta_0}$,  the set of rows $\{R_\delta\}_{\delta \in S_{\delta_0}}$ corresponds to a set of vectors $\Phi:=\{ \chi_{n} +(-1)^k \chi_0 \}_{n \in \text{Im}(\varphi)}$. To see that these are linearly independent, let $\mc C= \{ B\delta \bmod{2d} \mid \delta \in S_{\delta_0} \}$. Since $\delta_0\neq 0\bmod 0$, we see that $0\not\in \mc C$. That is, none of the characters $\chi_{B\delta}$ correspond to the trivial character $\chi_0$.

A linear relation among  the members of $\Phi$ looks like: 
\begin{align}
\sumn_{n \in \mc C} c_n \chi_n + (-1)^k \big( \sumn_{n \in \mc C} c_n \big) \chi_0 &= 0.
\end{align}
Using the orthogonality of characters, we see that the term $\chi_0$ is linearly independent from all $\chi_n$ ($n \neq 0$). Thus $c_n = 0$ for all $n$. Thus, the rows are linearly independent. Thus, the matrix $\mc S_{\delta_0}$ has maximal rank.

\textbf{$g$ even:} In this case, injectivity of $\varphi$ is straightforward since $(B,d)=1$. Thus the set of rows $\{R_\delta\}_{\delta \in S}$ corresponds to a set of vectors $\{ \chi_{n} +(-1)^k \chi_0 \}_{n \in \text{Im}(\varphi)}$. The linear independence now follows as in the previous case. 

Thus, we have shown that $\mc S_{\delta_0}$ has maximal rank. This completes the proof of $\textbf{Claim 2}$. As a consequence, we have proved that the theta components $h_{2T\gamma+L\bm\delta}$ with $\delta$ in a non-zero square class are zero.

\subsubsection{Square classes that are zero modulo $d$:}\label{sec:zerosqrcls} In this section, we prove that the remaining theta components $h_{2T\bm\gamma+L\bm\delta}=0$ when $\delta\equiv 0\bmod d$. Towards this, we use the vanishing of the theta components we obtained in the previous section. That is, from the previous section, we get that 
\begin{align}
    h_{2T\bm\mu_1+L\bm\mu_2}=0 \text{ for all } \mu_1\bmod L \text{ whenever } \mu_2^2\not\equiv 0\bmod d.
\end{align}
Fixing one such $\bm \mu_2$, from \eqref{hmutrans-sec5}, we get the following system of equations.
\begin{align}
    \sum_{\gamma\bmod L} \underset{{\delta}^2\equiv {\delta}_{0}^2 \bmod d_T}{\sum_{\delta\bmod d_T^*}}\varepsilon_{LT}(2T{\bm{ \gamma}}+L\bm\delta, 2T\mu_1+L\bm{\mu_2};  M)h_{2T{\bm{ \gamma}}+L\bm{\delta}}(\tau)=0 \text{ for all } \mu_1\bmod L,
\end{align}
where $\delta_0= 0$ or $\delta_0=d$. Note that in this case, the congruence $\delta^2\equiv {\delta}_{0}^2 \bmod d_T$ has only one solution modulo $d_T^*$ given by $\delta=\delta_0$.

Then from  Lemma \ref{lem:epsL} and \ref{lem:epsT},  $\varepsilon_{LT}(2T{\bm{ \gamma}}+L\bm\delta, 2t\mu_1+L\bm{\mu_2};  M)$ can be written as
\begin{align}
&= \omega(M)\varepsilon_L(2t_g\bm\gamma, \mu_1,; M)\varepsilon_T(L\bm\delta, \bm{\mu_2}; M)\\
&=  E_{M}G_1(t_g(b-\overline{\ell}), -2t_g(\gamma+\overline{\ell a}\mu_1), L)G_1(Lt(b-\overline{\ell}),2Lt (\ell_1\mu_2-\delta), d_T),
\end{align}
where $E_{M}$ is non zero constant, and we define
\begin{align}
    \ell_1:=\begin{cases}
        \overline{\ell}L+1& \text { if } g \text{ is odd};\\
        -(\overline{a\ell}+1) & \text { if } g \text{ is even}.
    \end{cases}
\end{align}
By choosing $b$ such that $(b-\overline{\ell}, d_T)=1$, $G_1(Lt(b-\overline{\ell}),2L t(\ell_1\mu_2-\delta), d_T)\neq 0$. Thus we are left with
\begin{align}
\sum_{\gamma\bmod L}G_1(t_gB_1, -2t_g(\gamma+\overline{\ell a}\mu_1), L) h_{2T\bm\gamma+L\bm{\delta_0}}(\tau)=0 \text{ for all } \mu_1\bmod L.    
\end{align}
The one dimensional Gauss sum $G(t_gB_1, -2t_g(\gamma+\overline{\ell a}\mu_1), L)$ can be written as below. 
\begin{align}
\sum_{x\bmod L}e_L(t_gB_1x^2-2t_gx(\gamma+\overline{\ell a}\mu_1))=e\left(-\frac{t_g\overline{B_1}(\gamma+\overline{\ell a}\mu_1)^2}{L}\right)G(t_gB_1, 0, L).
\end{align}
Since $G(t_gB_1, 0, L)\neq 0$, it is enough to consider the following system of equations.
\begin{align}\label{zerosqrclass-LEq}
\sum_{\gamma\bmod L}e\left(-\frac{t_g\overline{B_1}(\gamma+\overline{\ell a}\mu_1)^2}{L}\right) h_{2T{\bm{ \gamma}}+L\bm{\delta_0}}(\tau)=0 \text{ for all } \mu_1\bmod L.     
\end{align}
The vanishing of $h_{2T\bm\gamma+L\bm\delta_0}$ now follows from the following lemma.
\begin{lemma}
Let $L$ be an odd square-free integer, and let $B, C$ be integers coprime to $L$. Consider the $L \times L$ matrix $M$ defined by the entries:
\[
M_{r,s} = e\left( \frac{B(r+Cs)^2}{L} \right),
\]
where the indices $r, s$ range over $0, \dots, L-1$. Then the rank of $M$ is $L$.
\end{lemma}

\begin{proof}
Let $\zeta_L = e^{2\pi i / L}$. Expanding the quadratic term in the exponent, the $(r,s)$-th entry of the matrix $M$ can be factored as:
\begin{align*}
M_{r,s}= \zeta_L^{B(r+Cs)^2}= \zeta_L^{Br^2} \cdot \zeta_L^{2BCrs} \cdot \zeta_L^{BC^2s^2}.
\end{align*}
This factorization allows us to express $M$ as a product of three matrices:
\[
M = D_R\cdot V \cdot D_S,
\]
where $D_R = \operatorname{diag}(\zeta_L^{B \cdot 0^2}, \dots, \zeta_L^{B(L-1)^2})$, $D_S = \operatorname{diag}(\zeta_L^{BC^2 \cdot 0^2}, \dots, \zeta_L^{BC^2(L-1)^2})$ and $V$ is a matrix with entries $V_{r,s}=(\zeta_L^{2BC})^{rs}$.

Since $D_R$ and $D_S$ are diagonal matrices with roots of unity on the diagonal, they are non-singular (invertible). Thus, $\operatorname{rank}(M) = \operatorname{rank}(V)$.

The matrix $V$ is a Vandermonde matrix whose columns are the roots of unity. Its rank is $L$ if and only if $(2BC, L) = 1$. This is definitely the case as $L$ is odd and from \eqref{zerosqrclass-LEq}, $(t_g\overline{B_1\ell a}, L)=1$. Therefore, $V$ has full rank $L$, and implies that $\operatorname{rank}(M) = L$.
\end{proof}
\subsubsection{Completing the proof of the proposition.} Recall the we start with the hypothesis $h_{L\bm{\mu_2}}=0$ for all $\bm{\mu_2}\in \mc N_{2T}^{g,1}$. Our aim is to prove that $h_{2T{\bm{ \gamma}}+L\bm{\delta}}=0$ for all $\gamma\mod L$ and all $\bm{\delta}\in \mc N_{2T}^{g,1}$.
\begin{enumerate}
    \item When $\bm\delta$ is such that $\delta\not\equiv 0 \bmod d$, the vanishing of such theta components is covered in section \ref{sec:non-zerosqrcls}.
    \item When $\bm\delta$ is such that $\delta\equiv 0 \bmod d$, the vanishing of such theta components is covered in section \ref{sec:zerosqrcls}.
\end{enumerate}
Thus, we get that all the theta components of $\phi$ are zero and thus $\phi=0$, a clear contradiction. This completes the proof of the proposition. \qed
%===============================

%===============================
\section{Existence of non-zero primitive Theta components}

In this section, we prove a more general result that holds for the smaller group $\Gamma_1^{1,g,J}(N):= \Gamma_1^{(1)}(N)\ltimes \mc L_N \subset \Gamma_\pm^{1,g,J}(N, L)$. 

This proposition, given below, assures the existence of non-zero `primitive' theta components of a non-zero $\phi\in J_{k,T}(\Gamma_1^{1,g,J}(N), \mc L_N)$. The notion of `primitive' vectors for a matrix index $M \in \Lambda_g^+$ is recalled in Definition~\ref{prelim:prim}. This result is comparable to the results in \cite[Proposition 3.5]{boche-das} and \cite[Theorem 1]{anamby2024non}, although here we are dealing with Jacobi forms on a different subgroup and a smaller lattice. As a result, the calculations are more subtle.

This result forms the link between the discriminants of matrices indexing the Fourier coefficients in degrees $n-1$ and $n$, respectively, under the induction procedure. Here we remind the reader that the conditions $(N, d(T))=1$ and $d(T)$ odd and square-free are necessary to arrive at the forthcoming conclusions.

\begin{proposition}\label{prop:thetacom}
    Let $T\in \Lambda^+_g$ be such that  $(d(T), N)=1$ and $d(T)$ odd, square-free. Suppose $\phi\in J_{k, T}(\Gamma_1^{1,g, J}(N), \mc L_N)$ is non-zero, then there exists a $T$-primitive $\mu $ such that $h_{\mu}\neq 0$.
\end{proposition}

First, we show that Proposition \ref{prop:thetacom} implies the existence of non-zero primitive theta components for $\phi\in J_{k, T}(\Gamma_\pm^{1,g, J}(N, L), \mc L_L)$. This is not immediately a priori, as the theta components of $\phi$ can be different w.r.t. two groups $\Gamma_\pm^{(1)}(N, L)\ltimes \mc L_L$ and $\Gamma_1^{(1)}(N)\ltimes \mc L_N$. Towards this, for any divisor $M$ of $N$, we first choose the coset representatives $\mc N_{2MT}^{g,1}$ as in \eqref{coset-reps}. That is,
\begin{align}\label{coset-repsN}
  \mc N_{2MT}^{g,1}=\{(2T)\bm{\mu_1}+M\bm{\mu_2}: \bm{\mu_1} \in  \mc N_{M}^{g,1}, \bm{\mu_2} \in \mc N_{2T}^{g,1} \}.
\end{align}
Then we have the following lemma about the $T$-primitive representatives.
\begin{lemma}\label{lem:cosNNT}
    Let $\mu=(2T)\bm{\mu_1}+M\bm{\mu_2}\in \mc N_{2MT}^{g,1} $ with $ \bm{\mu_1} \in  \mc N_{M}^{g,1}, \bm{\mu_2} \in \mc N_{2T}^{g,1}$. Then $\mu$ is $T$-primitive iff  $\bm{\mu_2}$ is $T$-primitive.
\end{lemma}
\begin{proof}
    The proof is immediate by noting that,  $(2T)^{-1}[\mu]=(2T)[\bm{\mu_1}]+M^2(2T)^{-1}[\bm{\mu_2}]+ 2M \bm{\mu_1}^t\bm{\mu_2}$, and considering the denominator.
\end{proof}
Let $\phi\in J_{k, T}(\Gamma_\pm^{1,g, J}(N, L), \mc L_L) \subseteq J_{k, T}(\Gamma_1^{1,g,J}(N), \mc L_N)$. Then we can consider the theta expansion of $\phi$ with respect to both of these groups:
\begin{align}
    \phi(\tau, z)=\sum_{\bm\mu\in \mc N_{2LT}^{1,n}} g_{\bm\mu}(\tau)\cdot \Theta_{LT, \bm\mu}(L\tau,z), \q \phi(\tau, z)=\sum_{\bm\mu\in \mc N_{2NT}^{1,n}} h_{\bm\mu}(\tau)\cdot \Theta_{NT, \bm\mu}(N\tau,z).
\end{align}
We have the following lemma.
\begin{lemma}
    Let $\phi\in J_{k, T}(\Gamma_\pm^{1,g, J}(N, L), \mc L_L)$ and $g_{\bm\mu}$ and $h_{\bm\mu}$ be as above. If $h_{\bm\mu}\neq 0$ for some $T$-primitive $\bm\mu$, then $g_{\bm\mu'}$ is non-zero for some $T$-primitive $\bm\mu'$.
\end{lemma}
\begin{proof}
Next, we note that the coset representative $\mc N_{2NT}^{g,1}$ can also be taken as below.
\begin{align}
 \mc N_{2NT}^{g,1}=\{(2LT){\bm\mu_1}+\bm{\mu_2}: \bm{\mu_1} \in  \mc N_{N/L}^{g,1}, \bm{\mu_2} \in \mc N_{2LT}^{g,1} \}. 
\end{align}
Now consider 
\begin{align}
    \sum_{\bm\eta\bmod N/L}\Theta_{NT, \bm\mu+2LT\bm\eta}(N\tau, z)&=\sum_{\bm\eta\bmod N/L}\underset{R\equiv \bm\mu+2LT\bm\eta\bmod 2NT}{\sum_{R\in M_{g,1}(\z)}}
e(\tr( (NT)^{-1}[R/2]N\tau+z^tR).
\end{align}
Making a change of variable $R'=R-2LT\bm\eta$, we see that
\begin{align}
 \sum_{\bm\eta\bmod N/L}\Theta_{NT, \bm\mu+2LT\bm\eta}(N\tau, z)=\underset{R'\equiv \bm\mu\bmod  2LT}{\sum_{R'\in M_{g,1}(\z)}}
e(\tr( T^{-1}[R'/2]\tau+z^tR')=\Theta_{LT, \bm\mu}(L\tau, z).
\end{align}
Thus we get
\begin{align}
\phi(\tau, z)&=\sumn_{\bm\mu\in \mc N_{2LT}^{1,n}} g_{\bm\mu}(\tau)\cdot \Theta_{LT, \bm\mu}(L\tau,z)\\
&=\sumn_{\bm\eta\bmod N/L}\sumn_{\bm\mu\in \mc N_{2LT}^{1,n}} g_{\bm\mu}(\tau)\cdot \Theta_{NT, \bm\mu+2LT\bm\eta}(N\tau,z)\\
&\sumn_{\bm\mu'\bmod 2NT}\sumn_{\bm\eta\bmod N/L}g_{\bm\mu'-2LT\bm\eta}(\tau)\Theta_{NT, \bm\mu'}(N\tau,z).
\end{align}
Thus, by the linear independence of theta functions, we get
\begin{align}
    h_{\bm\mu}(\tau)=\sumn_{\bm\eta\bmod N/L}g_{\bm\mu-2LT\bm\eta}(\tau).
\end{align}
Further, from Lemma \ref{lem:cosNNT}, note that if $\bm\mu$ is $T$-primitive, then $\bm\mu-2LT\bm\eta$ is $T$-primitive for all $\bm\eta$. Thus if $h_{\bm\mu}\neq 0$ for some $T$-primitive $\bm\mu$, then $g_{\bm\mu-2LT\eta}\neq 0$ for some $\bm\eta$ and this completes the proof of the lemma.
\end{proof}
\subsection{Proof of Proposition~\ref{prop:thetacom}.}
\begin{proof}
We start with $M=-\smat{0}{-1}{1}{0}\smat{1}{N}{0}{1}\smat{0}{-1}{1}{0}= \smat{1}{0}{N}{1}$ and $M'=\smat{1}{0}{N^2}{1}$. Then $M, M'\in \Gamma_1(N)$ and
\begin{align}\label{commuteN}
    \smat{N}{0}{0}{1}M'= M\smat{N}{0}{0}{1}.
\end{align} 
Now we can write
\begin{align}
    \phi(M'\tau, z(N^2\tau+1)^{-1}) &= \sumn_{\bm\mu} h_\mu(M'\tau)\Theta_{NT, \bm\mu}(NM'\tau, z(N^2\tau+1)^{-1}))\\
    &= \sumn_{\bm\mu} h_{\bm\mu}(M'\tau)\Theta_{NT, \bm\mu}(M(N\tau), z(N(N\tau)+1)^{-1})).
\end{align}
Using the transformation properties of $\Theta_{ N T, \bm\mu}(M(N\tau), z(N(N\tau)+1)^{-1}))$ from \eqref{theta-trans} and that $\phi|M'=\phi$, we get
\begin{align} \label{theta-phi}
\sum_{\bm\mu}h_{\bm\mu}(\tau) \Theta_{NT,\bm\mu}(N\tau, z)= (N^2\tau+1)^{-k+1/2}\sum_{\bm\mu, \bm\nu\in \mc N_{2NT}^{g,1}}\varepsilon_{NT}(\bm\mu, \bm\nu; M) h_{\bm\nu}(M'\tau)   \Theta_{NT, \bm\mu}(N\tau,z) .
\end{align}
Next, by using the linear independence of the theta series in \eqref{theta-phi} and a change of variable $\tau\mapsto M'^{-1}\tau$, we get,
\begin{align}\label{hmutrans}
    h_{\bm\mu}(M'^{-1}\tau)= (N^2\tau+1)^{-k+1/2}\sumn_{\bm\nu\in\mc N_{2NT}^{g,1} }\varepsilon_{NT}(\bm\mu, \bm\nu;M) h_{\bm\nu}(\tau).
\end{align}
As in the previous section, using \eqref{epsilonmunu} we get the following expression for $\varepsilon_{NT}(\mu, \nu; M)$.
\begin{align}
    \varepsilon_{NT}(\mu, \nu;M)&= \omega(M)\sumn_{\bm{\eta}\in \mc N_{2NT}^{g,1}} e_N\left(\tr(\frac{1}{4} NT^{-1}[\bm\eta]+ \frac{1}{2}\bm\eta^t T^{-1}\mu- \frac{1}{2}\nu^t T^{-1}\bm\eta)\right)\\
    &= \omega(M)\sumn_{\bm\eta\in \mc N_{2NT}^{g,1}} e_4\left(\tr(N (NT)^{-1}[\bm\eta]+2\bm\eta^t(NT)^{-1}(\mu-\nu) )\right).
\end{align}
Here we use \eqref{epsilonmunu} and $\omega(M)$ is the complex number of absolute value $1$ as in \eqref{epsilonmunu}.

Since $N, M$ are fixed throughout this proposition, we just write $\varepsilon_T(\mu, \nu)$. If $N\nmid(\mu-\nu)$, then considering $i$ such that $N\nmid (\mu_i-\nu_i)$, we see  (by considering the vector $e_i=(0,.,1,.,0)$ with $1$ at $i$-th place and making a change of variable $\eta\mapsto \eta+2Te_i^t$) that,
\begin{align}
 \varepsilon_{NT}(\bm\mu, \bm\nu;M)= e((\mu_i-\nu_i)/N)  \varepsilon_{NT}(\bm\mu, \bm\nu;M).   
\end{align}
Thus $\varepsilon_{NT}(\bm\mu, \bm\nu)=0$ if $N\nmid(\mu-\nu)$.  When $N|(\bm\mu-\bm\nu)$, by writing $\bm\eta\in \mc N_{2NT}^{g,1}$ as $\eta=(2T) \bm{\eta_1}+N\bm{\eta_2}$ with $\bm{\eta_1}\in \mc N_{N}^{g,1}$ and $\bm{\eta_2}\in \mc N_{2T}^{g,1}$ as given in \eqref{coset-repsN}, we can compute:
\begin{align}
\tr\left(T^{-1}[\bm\eta]+\bm\eta^tT^{-1}\frac{(\bm\mu-\bm\nu)}{N} \right)= \bm\eta_1^t\bm\eta_1+N^2 T^{-1}[\bm\eta_2/2]+N\bm\eta_1^t\bm\eta_2+ \frac{\bm\eta_1^t(\bm\mu-\bm\nu)}{N}+N\bm\eta_2^tT^{-1}\frac{(\bm\mu-\bm\nu)}{N}. 
\end{align}
Since $(N, d(T))=1$, $N$ is invertible $\mod 2T$,  we can make a change of variable $\bm\eta_2= N^{-1}\bm\eta_3 \bmod 2T$, and since the $\bm\eta_1$ sum contributes $N^g$, we finally get
\begin{align}\label{epsmunufin}
    \varepsilon_T(\bm\mu, \bm\nu)= \omega(M)\sumn_{\bm\eta_3\in \mc N_{2T}^{g,1}} e_4\left(\tr(T^{-1}[\bm\eta_3]+\bm\eta_3^tT^{-1}(\bm\mu-\bm\nu)/N )\right).
\end{align}
Now suppose that $h_{\bm\mu}=0$ for all $T$-primitive $\bm\mu\in \mc N_{2NT}^{g,1} $ . Thus for any such $\bm\mu$, from \eqref{hmutrans}, we get
\begin{align}
\sumn_{\bm\nu\in \mc N_{2NT}^{g,1}}\varepsilon_T(\bm\mu, \bm\nu) h_{\bm\nu}(\tau)=0.
\end{align}
Substituting for $\varepsilon_T(\bm\mu, \bm\nu)$ from \eqref{epsmunufin}, and using the change of variable $\bm\eta\mapsto \bm\eta-(\bm\mu-\bm\nu)/N$, we get
\begin{align}
\underset{N|(\bm\mu-\bm\nu)}{\sum_{\bm\nu\in \mc N_{2NT}^{g,1}}}e\left(\tr(T^{-1}[(\bm\mu-\bm\nu)/2N])\right)\sum_{\bm\eta\in \mc N_{2T}^{g,1}} e_4\left(\tr(T^{-1}[\bm\eta])\right) h_{\bm\nu}(\tau)=0.
\end{align}
The sum over $\bm\eta$ is nothing but the quadratic Gauss sum $G_g(1,0, 2T)$ as in Section \ref{sec:HDGauss}
and it is non-zero by Lemma \ref{gauss-nonvan}. As a result, for each $T$-primitive $\bm\mu$, we get a linear equation
\begin{align}\label{hnueq}
    \underset{N|(\bm\mu-\bm\nu)}{\sum_{\bm\nu\in \mc N_{2NT}^{g,1}}}e\left(\tr(T^{-1}[(\bm\mu-\bm\nu)/2N])\right) h_{\bm\nu}(\tau)=0.
\end{align}
Now we use the coset representatives of $\mc N_{2NT}^{g,1}$ as in\eqref{coset-repsN}. Fixing a $\bm\mu= (2T)\bm{\alpha}+N\bm{\beta}$, let $\bm\nu= (2T)\bm\gamma+N\bm\delta$. Then $N|(\bm\mu-\bm\nu)$ iff $N|(\bm{\alpha}-\bm\gamma)$. Thus, we can write
\begin{align}
    \underset{N|(\bm{\alpha}-\bm\gamma)}{\sum_{\bm\gamma\in \mc N_{N}^{g,1}}}\sum_{\bm\delta\in \mc N_{2T}^{g,1}}e\left(\tr(T^{-1}[(\bm{\beta}-\bm\delta)/2])\right) h_{(2T)\bm\gamma+N\bm\delta}(\tau)&=0.
\end{align}
By our choice of coset representatives, we have that $\bm{\alpha}=\bm\gamma$. Thus, the $\bm\gamma$ sum collapses and we get
\begin{align}
    \sumn_{\bm\delta\in \mc N_{2T}^{g,1}}e\left(\tr(T^{-1}[(\bm{\beta}-\bm\delta)/2])\right) h_{(2T)\bm{\alpha}+N\bm\delta}(\tau)&=0.
\end{align}
Using the transformation properties of $h_{\nu}$ under $\tau \mapsto \tau+t$, we have
\begin{align}\label{hnu-char}
   \sumn_{\delta\in \mc N_{2T}^{g,1}}e\left(\tr(T^{-1}[(\bm{\beta}-\delta)/2])\right) e(T^{-1}[N\delta/2]t)h_{(2T)\bm{\alpha}+N\delta}(\tau)&=0.
\end{align}
The linear independence of the additive character $t \mapsto e(h \cdot t)$ shows that we have an equation of the form
\begin{align}\label{sqrclass}
   \underset{T^{-1}[N\bm\delta] - T^{-1}[N\bm\delta_0] \in 4 \z }{\sum_{\bm\delta\in \mc N_{2T}^{g,1}} } e\left(\tr(T^{-1}[(\bm\beta-\bm\delta)/2])\right) h_{(2T)\bm{\alpha}+N\bm\delta}(\tau)&=0.
\end{align}
Next, we note that, since $(N, d(T))=1$,  
\begin{align}
T^{-1}[N\bm\delta] - T^{-1}[N\bm\delta_0] \in 4 \z \iff T^{-1}[\bm\delta] - T^{-1}[\bm\delta_0] \in 4 \z .
\end{align}
To see this, we multiply both sides by $\det(2T)$ to write an equivalent statement: $ \det(2T) (2T)^{-1}[N\bm\delta] \equiv \det(2T) (2T)^{-1}[N\bm\delta_0] \bmod{\det(2T)}$. The claim now follows immediately using standard arguments.

As a consequence, we see that the equation in \eqref{sqrclass} is the same as the linear equation
\begin{align}
   \underset{T^{-1}[\bm\delta] - T^{-1}[\bm\delta_0] \in 4 \z }{\sum_{\bm\delta\in \mc N_{2T}^{g,1}} } e\left(\tr(T^{-1}[(\bm\beta-\bm\delta)/2])\right) h_{(2T)\bm{\alpha}+N\bm\delta}(\tau)&=0.
\end{align}
For any fixed $\bm\delta_0\in \mbb Z^g$, let us define the set $S_{\delta_0}:=\{\bm\delta \in \mc N_{2T}^{g,1}:T^{-1}[\bm\delta] - T^{-1} [\bm\delta_0] \in 4\mbb Z\}$. As $\bm{\beta}$ varies over $T$-primitive elements in $\mc N_{2T}^{g,1}$ and $\bm\delta\in S_{\delta_0}$, let us define the matrix
\begin{align}
    \mc M_{\bm\delta_0}:=\left(e\left(T^{-1}[(\bm\beta-\bm\delta)/2]\right) \right).
\end{align}
Let $a_{\bm\beta}= e(T^{-1}[\bm\beta/2])$ and $b_{\bm\delta}= e(T^{-1}[\bm\delta/2])$, then $\mc M_{\delta_0}$, can be written as 
\begin{align}
    \mc M_{\bm\delta_0}= \left(a_{\bm\beta} b_{\bm\delta} e\left(\bm{\beta}^{t}T^{-1}\bm\delta/2\right) \right).
\end{align}
Since $a_{\bm{\beta} } b_{\bm\delta}\neq 0$, by row and column reductions (see also section \ref{sec:zerosqrcls}), we see that $\mc M_{\bm\delta_0}$ has the same rank as the matrix
\begin{align}
   \mc S_{\bm\delta_0} :=\left( e\left(\bm{\beta}^{t}T^{-1}\bm\delta/2\right) \right).
\end{align}
Now, since $d(T)$ is odd and square-free, from \cite[Proposition~3.5]{boche-das}, the matrix $\mc S_{\bm\delta_0}$ has maximal rank and more importantly, the rank does not depend on $\bm{\alpha}\in \mc N_{N}^{g,1}$. Thus $h_{(2T)\bm{\alpha}+N\bm\delta}(\tau)=0$ whenever $\bm\delta$ is not $T$-primitive. That is, we get (using Lemma \ref{lem:cosNNT})
\begin{align}
h_{(2T)\bm\gamma+N\bm\delta}=\begin{cases}
 0 & \text{ if }  (N, (2T)[\bm\gamma])=1, \bm\delta \text{ is not } T\text{-primitive};\\
0 & \text{ if }  (N, (2T)[\bm\gamma])>1, \bm\delta \text{ is not } T\text{-primitive}.
\end{cases}   
\end{align}
By our initial assumption that $h_{\bm\mu}=0$ for all $T$-primitive $\bm\mu$, using Lemma \ref{lem:cosNNT}, we get that 
\begin{align}
h_{(2T)\bm\gamma+N\bm\delta}=    \begin{cases}
0 & \text{ if }  (N, (2T)[\bm\gamma])=1, \bm\delta \text{ is }  T\text{-primitive};\\
0 & \text{ if }  (N, (2T)[\bm\gamma])>1, \bm\delta \text{ is }  T\text{-primitive}.
\end{cases}
\end{align}
Since all such $(2T)\bm\gamma+N\bm\delta$ exhaust the set $\mc N_N^{g,1}\times \mc N_{2T}^{g,1}$, we get that $h_{\bm\nu}=0$ for all $\bm\nu\in \mc N_{2NT}^{g,1}$. This implies that $\phi=0$, a clear contradiction, and this completes the proof. 
\end{proof}
As a final result of this section, for any $T\in \Lambda_n$, we have the following non-vanishing criteria for the Fourier coefficients of $H_{\bm\mu}(\tau):=h_{\bm\mu}(d_T\tau)$, where $d_T$ is as in \eqref{def:d_T}.
\begin{lemma}\label{lem:primthetaFC}
    Let $\phi \in J_{k,T}(\Gamma_1^{1, n, J}(N), \mc L_N)$ and $h_{\bm\mu}$ be a $T$-primitive theta component. Then $H_{\bm\mu}$ has all non-zero Fourier coefficients away from $d_T$.
\end{lemma}
\begin{proof}
  As $\bm\mu$ is $T$-primitive, we have the property that $\displaystyle (2T)^{-1}[\bm\mu]= \frac\alpha{d_T}$ with $(\alpha,d_T)=1$. Now from the Fourier expansion of $h_{\bm\mu}$ in \eqref{hmu-fe}, we see that if  $\bm\mu$ is $T$-primitive, then its Fourier expansion can be written as
\begin{align} \label{prim-hmu}
    H_{\bm\mu}(\tau)=\sumn_{n' \geq 1,}c_\phi(\frac{n'+\alpha}{d_T}, \bm\mu) e(n' \tau).
\end{align}
Since $c_{\phi}(n/d_T, \bm\mu)=0$ when $d_T\nmid n$ and $(\alpha, d_T)=1$, we see that the Fourier expansion of $H_{\bm\mu}$ is supported away from $d_T$. Note that we have arrived at \eqref{prim-hmu} from \eqref{hmu-fe} by rewriting  $n- (2T)^{-1}[\bm\mu]= \frac{n'}{d_T}$. This completes the proof.
\end{proof}

%==============================

%==============================
\section{\texorpdfstring{Getting all the diagonals of $T$ away from $L$}{tnnimplytii}} \label{sec:tnn-diagonal}
%Proof of Proposition \ref{tnn-diagonal}.

In this section, we prove that the condition $(T(n,n), L)=1$ implies the condition $(T(i,i), L)=1$ for all $1\le i\le n$.  This serves as a linkage between Theorem \ref{mainthm} and Theorem \ref{diagonal-sqrfree} (first step in Figure \ref{fig:flow}). In fact, we prove a more general version for the smaller group $\Gamma_1^{(n)}(N)$. More precisely, we show the following.
%\tnndiagonal*
\begin{proposition}
\label{tnn-diagonal}
Let $F\in  M_{\rho}(\Gamma_1^{(n)}(N))$ be such that $A(F, T_0)\neq 0$ for some $T_0\in\Lambda_n$ with $(T_0(n,n), N)=1$. Then there exists a $T_1\in \Lambda_n$ with $(T_1(j,j), N)=1$ for all $1\le j\le n$ such that $A(F, T_1)\neq 0$.
\end{proposition}

\begin{proof}
We use induction on $n$. The proposition is obviously true when $n=1$, see Lemma~\ref{ind-basecase}. Assume the statement for all degrees $<  n$  and all levels. One also needs to check the conditions on the weights, but this is exactly the same as presented in \cite[Section~4.3]{boche-das}, and will not be repeated here.

Let $T_0\in \Lambda_n$ be such that $(T_0(n,n), N)=1$ and $A(F, T_0)\neq 0$. We consider the non-zero modular form $ F_{T_0}(Z)=\sum_{T\equiv T_0\bmod N}A(F, T)e(\tr(TZ)) \in M_{\rho}(\Gamma_1^{(n)}(N^2))$ and the corresponding degree $n-1$ form $0\neq (F_{T_0})_{(n-1)}\in M_{\rho_*}(\Gamma_\pm^{(n-1)}(N^2))$ (see Section \ref{sec:Fcirc}). By definition, $(F_{T_0})_{(n-1)}$ satisfies the induction hypothesis. Thus there exist $S_1 \in \Lambda_{n-1}$ with $(S_1(i,i), N)=1$ for $1\le i\le n-1$ and $A((F_{T_0})_{(n-1)}, S_1)\neq 0$. As a consequence,  from Proposition \ref{sca-j}, we get a non-zero, scalar-valued Jacobi form $\phi_{S_1}$ of index $S_1$ from the Fourier-Jacobi expansion of $F_{T_0}$.

Now from any non-zero Fourier coefficient of $\phi_{S_1}$ we get a non-zero Fourier coefficient $A(F_{T_0}, T_1)$ of $F_{T_0}$, where $T_1=\smat{*}{*}{*}{S_1}$. Then we look  at
\begin{align}
    F_{T_0,T_1}(Z)= \sum_{T\equiv T_1 \bmod N}A(F_{T_0}, T)e(\tr(TZ))\in M_{\rho}(\Gamma_1^{(n)}(N^3)).
\end{align}
From the congruence condition in the definition of $F_{T_0,T_1}$, it is clear that for all matrices in $\{ T: A(F_{T_0, T_1}, T)\neq 0\}$, the diagonals $(T(j,j), N)=1$ for all $2\le j\le n$.

Now we continue with the proof by considering $0\neq (F_{T_0,T_1})^{(n-1)}\in M_{\rho^*}(\Gamma_\pm^{(n-1)}(N^3))$. For any $S\in \Lambda_{n-1}$ with $A((F_{T_0,T_1})^{(n-1)},S)\neq 0$, by construction, $(S(n-1,n-1),N)=1$. Thus, $(F_{T_0,T_1})^{(n-1)}$ satisfies the induction hypothesis; as a result, we get that $(F_{T_0,T_1})^{(n-1)}$ has good diagonal property. That is, we get an $S_3\in \Lambda_{n-1}$ with $(S(i,i), N)=1$ for all $1\le i\le n-1$ such that $A((F_{T_0,T_1})^{(n-1)}, S_3)\neq 0$. Thus we get a $T_2=\smat{S_3}{*}{*}{*}$ such that $A(F_{T_0,T_1},T_2)\neq 0$.

Now, by definition of $F_{T_0,T_1}$, we have that $T_2\equiv T_1\bmod N$. Thus $(T_2(n,n), N)=1$. Since $T_2(i,i)=S_3(i,i)$ for $1\le i\le n-1$, we get $(T_2(i,i), N)=1$ for all $1\le i\le n$. This completes the proof of the proposition.
\end{proof}
%==================================

%==================================
\section{Proof of Theorem \ref{diagonal-sqrfree}}\label{sec:thmdiag}

\begin{lemma}\label{lem:psd-pd}
    Let $F \in \Mrho{n}{N}$ and $T_0\in \Lambda_n$ be as in the hypothesis of Theorem \ref{mainthm}. Then $F_{T_0}\neq 0$ and there exists a $T_1\in \Lambda_n^+$ with $(T_1(n,n),N)=1$ such that $A(F_{T_0}, T_1)\neq 0$.
\end{lemma}
\begin{proof}
The condition on weight,  $k(\rho) - \frac{n-1}{2} \ge 2 $, excludes the existence of `singular' forms (see \cite{weissauer1983vektorwertige}). Since the weight does not decrease upon the construction of $F_{T_0}$ (see Lemma \ref{lem:FTp}), $F_{T_0}$ is non-singular and hence there exists a $T_1\in \Lambda^+_n$ such that $A(F_{T_0}, T_1)\neq 0$. Further, since $T_1\equiv T_0 \bmod N$, we get that $(T_1(n,n), N)=(T_0(n,n),N)=1$ from Lemma \ref{lem:detcong}.
\end{proof}

\begin{remark}\label{rmk:psd-pd}
Lemma \ref{lem:psd-pd} ensures that the modular forms defined using the congruence condition w.r.t $T_0$ (see Lemma \ref{lem:FTN}) are not supported on $\Lambda_n \setminus
\Lambda^+_n$. This is crucial since $T$ in the conclusion of Theorem \ref{mainthm} are in $\Lambda^+_n$. Of course, there is no issue with cusp forms.
\end{remark}
We first establish the following non-vanishing proposition, which provides the crucial machinery required for the proof of Theorem \ref{diagonal-sqrfree}.
\begin{proposition}\label{alldiag-fund-prop}
Let $F\in M_\rho(\Gamma_1^{(n)}(N))$ be such that the Fourier expansion $F$ is supported on $T\equiv \mrm{diag}(t_1, t_2,...t_n)\bmod L$ with $(t_i, L)=1$ for all $1\le i\le n$. Suppose that $F$ has a non-zero Fourier-Jacobi coefficient $\varphi_S$  with $d(S)$ odd and square-free with $(d(S), L)=1$. Then there exist infinitely many $T\in \Lambda_n$ with $d(T)$ odd and square-free, $(d(T), L)=1$ and $A(F, T)\neq 0$.
\end{proposition}

\begin{proof}
Consider the non-zero Fourier-Jacobi coefficient $\varphi_S$ as in the hypothesis. Then from Proposition \ref{sca-j}, we get a non-zero component $\phi_S$ of $\varphi_S$, which is a scalar-valued Jacobi form of weight $k\ge k(\rho)$ for the group $\Gamma_1^{1, n-1, J}(N)$.

In this case, $d_S$ from \eqref{def:d_T} can be rewritten as below.
\begin{align} \label{d'-def}
    d_S = \begin{cases}
        d(S) & \text{ if   }  k - \frac{n-1}{2} \in \z ;\\
        4 d(S) & \text{ if   }  k - \frac{n-1}{2} \in \frac{1}{2}\z \setminus \z.
    \end{cases}
\end{align}
For $\mu \in \mc N_{2NS}^{1,n-1}$, let $h_\mu$ be a theta component of $\phi_S$. 
We first note that the level of $h_\mu$ divides $N d_S^2$, in fact it is in $M_{k - \frac{n-1}{2} }(\Gamma(N d_S^2))$. To see this, note that from \cite[\S 4, Part~B, p. 2022]{boche-das} $h_\mu$ has level $N d_S$ -- the $N$ contribution comes from the level $N$ of $\phi_S$, and $d_S$ from the Weil representation attached to $S$ -- which factors through $\Gamma(d_S)$. 

At this point, we need the following properties of the theta components $h_\mu$.

{\bf(P1)} 
There exists an index $\mu_0$ which is $S$-primitive such that $h_{\mu_0} \neq 0$. This is guaranteed by Proposition \ref{prop:thetacom}. Further, $\mu\equiv 0\bmod L$, since the Fourier coefficients of $F$ are supported on $T\equiv\mrm{diag}(t_1, t_2,...t_n)\bmod L$. We will work with this $\mu_0$.

{\bf(P2)} 
The Fourier expansion of $h_{\mu_0}$ is of the form $\sum_n a_n q^{n/ d_S}$, where $d_S$ is as defined in \eqref{d'-def}. In fact,
\begin{align}\label{hmu-FT}
        h_{\mu}(\tau) =\sumn_{n \ge (2S)^{-1}[\mu_0]} c_{\phi_S}(n ,\mu_0) e( (n- (2S)^{-1}[\mu_0]) \cdot \tau).
\end{align}

{\bf(P3)} 
We define $H_{\mu_0}(\tau):= h_\mu(d_S \tau)$. The Fourier expansion of $H_{\mu_0}$ is given by
\begin{align}
H_{\mu_0}(\tau) &=  \sumn_{n \ge (2S)^{-1}[\mu_0]} c_{\phi_S}(n ,\mu_0) e( d_S(n- (2S)^{-1}[\mu_0]) \cdot \tau) \\
    &=\sumn_{n \ge 0, \, n \equiv \alpha \bmod{d_S} } c_{\phi_S}(\frac{n}{d_S}+ (2S)^{-1}[\mu_0],\mu_0) e(n \cdot \tau) . \label{Hmu-fe}
\end{align}
Note that $c_{\phi_S}(n,\mu_0)=0$ if $n$ is not an integer.

{\bf(P4)} 
If we write the Fourier expansion of $H_{\mu_0}=\sum_{ n \equiv \alpha \bmod{d_S}} a(H_{\mu_0}, n)q^n$. 
For all $n \ge 1$ in the Fourier expansion of $H_{\mu_0}$, if $ a(H_{\mu_0}, n)\neq 0$, then $(n, d_S)=1$ (see Lemma \ref{lem:primthetaFC}). Here we have put, as before, $(2S)^{-1}[\mu_0] = \frac\alpha{d_S}$. Note that $(\alpha, d_S)=1$ by $S$-primitivity of $\mu_0$. 

{\bf(P5)} 
For each $n \geq 1$, one has $ c_{\phi_S}(\frac{n}{d_S}+S^{-1}[\mu_0/2], \mu) = A(G, \mc T )$, where $\mc T=\smat{\frac{n}{d_S}+S^{-1}[\mu_0/2]}{\mu_0/2}{\mu^t_0/2}{S}$ and has the following properties. $ d(\mc T) = n$. This is precisely the calculation in \cite[Section~3.4.1, 3.4.2]{boche-das} which follows without any changes. We see that $d(\mc T)=d_S(\frac{n}{d_S} +S^{-1}[\mu_0/2] - S^{-1}[\mu_0/2])=d_S\cdot \frac{n}{d_S}=n$.
  
{\bf(P6)} 
We claim that $H_{\mu_0} \in M_{k - \frac{n-1}{2} }(\Gamma_1(N' d_S^2))$. To see this, note that from the theta decomposition of $\phi_S$ as in \eqref{theta-decomp}, we see that $H_\mu$ is invariant under $\smat{1}{1}{0}{1}$. Thus so is $H_{\mu_0}$. 
Therefore it is invariant under the subgroup generated by $\smat{1}{1}{0}{1}, \Gamma(N'), \Gamma_1(d_S^2)$, which contains $\Gamma_1(N'd_S^2)$.
%One could also have seen this by observing that since $\phi_S$ and also the theta functions $\Theta_{NS, \mu}(N d' \tau,z)$ are all invariant under $\Gamma_1(N')$.
Finally, note that the change of variable $\tau \mapsto d_S \tau$ squared the $d_S$.

We can now put together all the ingredients assembled above to obtain the proof of the proposition. We will proceed in two stages.

\textbf{Stage-I: Away from $L$:}
From the hypothesis, the Fourier expansion of $F$ is supported on the congruence class $T\bmod L$, from Lemma \ref{lem:detcong}, we have that $n=d(\mc T)\equiv d(T)\bmod L$ and thus the Fourier expansion of $H_{\mu}$ (cf. \eqref{Hmu-fe}) is supported on $(n, L)=1$. That $n=d(\mc T)$ follows from {\bf (P5)}.

\textbf{Stage-II: Away from $d_S$:}
From {\bf(P6)} we also know that Fourier expansion of $H_{\mu}$ (cf. \eqref{Hmu-fe}) is supported on $(n, d_S)=1$. Since $(N,d_S)=1$, we therefore get that the Fourier expansion of $H_{\mu_0}$ is supported on $n$ with $(n,Nd_S)=1$, as desired. 

In addition, if $F$ is a cuspform, then so is $H_{\mu_0}$. Thus, from \cite[Theorem 4.5 and Theorem 4.6]{boche-das}, we have
\begin{align}\label{quant-cusp}
   \#\{m\le X: m \text{ odd, square-free}, (m,N)=1, a(H_{\mu_0}, m)\neq 0\} \gg\begin{cases} X(\log X)^{-1/2} &\text{ if } n \text{ is odd };  \\ 
X^{5/8 - \epsilon} &\text{ if } n \text{ is even}.
\end{cases}
\end{align}
When $F$ is non-cuspidal, we cannot guarantee at present that $\varphi_S$ is non-cuspidal. As a consequence, $H_{\mu_0}$ from \eqref{hmu-FT} can either be cuspidal or non-cuspidal (see \cite{boche-das} for a discussion on this). However, in both cases, as can be seen from \cite[Theorem 4.6]{boche-das}, the bound in \eqref{quant-cusp} is still valid. This gives the quantitative result in Theorem \ref{diagonal-sqrfree}.
\end{proof}
\subsection{Proof of the Theorem:} We now proceed to prove Theorem \ref{diagonal-sqrfree}. For the sake of convenience, we recall the statement of the theorem here. \diagonalsqrfree*
\begin{proof}
     We again use induction on $n$. For $n=1$, the result is trivially true. Now, assume that the result holds for degree $n-1$. We once again remind the reader that one needs to check the conditions on the weights; however, as mentioned in the previous section, these are identical to those presented in \cite[Section~4.3]{boche-das} and will not be repeated here.

\textbf{Step 1:}    Let $T_0$ be such that $(T_0(i,i), N)=1$ for all $1\le i\le n$ and $A(F, T_0)\neq 0$. As in the previous section, we consider the non-zero form $ F_{T_0}(Z) \in \Mrho{n}{N^2}$. Further, we have
\begin{align}
    A(F_{T_0}, T)=\begin{cases}
        A(F, T)& \text{ if } T\equiv T_0 \bmod N;\\
        0 & \text{ otherwise.}
    \end{cases}
\end{align}
In particular, $A(F_{T_0}, T)=0$ whenever $(T(i,i), N)\neq 1$ and since $0\neq A(F, T_0)=A(F_{T_0}, T_0)$, $F_{T_0}$ also satisfies the hypothesis of the Theorem.

Now consider $(F_{T_0})_{(n-1)}$ (see section \ref{sec:Fcirc}). Since $(F_{T_0})_{(n-1)}\neq 0$, by definition, $(F_{T_0})_{(n-1)}$ satisfies the induction hypothesis. Thus there exist $S_1 \in \Lambda_{n-1}$ with $S_1\equiv \mrm{diag}(s_1, s_2,..,s_{n-1})\bmod L$, $(d(S_1),L)=1$ such that $A((F_{T_0})_{(n-1)}, S_1)\neq 0$.  As a consequence,  from Proposition \ref{sca-j}, we get a non-zero, scalar-valued Jacobi form $\phi_{S_1}$ of index $S_1$ from the Fourier-Jacobi expansion of $F_{T_0}$. 

Now from any non-zero Fourier coefficient of $\phi_{S_1}$ we get a non-zero Fourier coefficient $A(F_{T_0}, T_1)$ of $F_{T_0}$, where $T_1=\smat{*}{*}{*}{S_1}$. Then we look  at
\begin{align}
    F_{T_0,T_1}(Z)= \sum_{T\equiv T_1 \bmod N}A(F_{T_0}, T)e(\tr(TZ))\in \Mrho{n}{N^3}.
\end{align}

From the congruence condition in the definition of $F_{T_0,T_1}$, it is clear that for all matrices in $\{ T: A(F_{T_0, T_1}, T)\neq 0\}$, the lower $n-1$ block is diagonal $\mod L$.

\textbf{Step 2:} Now we continue with the proof by considering $0\neq (F_{T_0,T_1})^{(n-1)}\in M_{\rho^*}(\Gamma_\pm^{(n-1)}(N^3, L))$.  For any $S\in \Lambda_{n-1}$ with $A((F_{T_0,T_1})^{(n-1)},S)\neq 0$, by construction, $(S(i,i),N)=1$ for all $1\le i\le n-1$. Thus, $(F_{T_0,T_1})^{(n-1)}$ satisfies the induction hypothesis; as a result, we get that $(F_{T_0,T_1})^{(n-1)}$ has the required properties. That is, we get an $S_2\in \Lambda_{n-1}$ with $S_2\equiv \mrm{diag}(s_1', s_2',..,s_{n-1}')\bmod L$, $d(S_2)$ odd, square-free and $(d(S_2),L)=1$ such that $A((F_{T_0,T_1})^{(n-1)}, S_2)\neq 0$. Thus we get a $T_2=\smat{S_2}{*}{*}{*}$ such that $A(F_{T_0,T_1},T_2)\neq 0$.

\textbf{Step 3:}\label{Step3:Theorem1.3} Now consider $0\neq F_{T_0,T_1,T_2}\in \Mrho{n}{N^4}$ and the non-zero degree $n-1$ form $(F_{T_0, T_1, T_2})_{(n-1)}$, which satisfies the induction hypothesis. Thus we get an $S_3 \in \Lambda_{n-1}$ with $S_3\equiv S_1\equiv \mrm{diag}(s_1,...,s_{n-1})\bmod L$, $d(S_3)$ odd, square-free  and $(d(S_1), L)=1$ such that $A((F_{T_0,T_1,T_2})_{(n-1)}, S_3)\neq 0$.  As a consequence,  from Proposition \ref{sca-j}, we get a non-zero, scalar-valued Jacobi form $\phi_{S_3}$ of index $S_3$ from the Fourier-Jacobi expansion of $F_{T_0,T_1,T_2}$.

We make the following observation with regard to $\phi_{S_3}$.
\begin{enumerate}
    \item $\phi_{S_3}\in J_{k, S_3}(\Gamma_\pm^{1,g, J}(N, L), \mc L_L)$ with $S_3$ diagonal $\bmod L$ and $(d(S_3), L)=1$;
    \item For any $r\in \mbb Z^{1,n- 1}/\mbb Z^{1,n- 1}$ we have that
\begin{align}
    c_{\phi_{S_3}}(n, r)= A\left(F_{T_0,T_1,T_2}, \smat{n}{r/2}{r^t/2}{S_3} \right).
\end{align}
Further, note that $\smat{n}{r/2}{r^t/2}{S_3}\equiv T_2\bmod L$. Since $T_2$ is diagonal $\mod L$ except for the entries $T_2(1,n)$ and $T_2(n,1)$, we get that $\phi_{S_3}$ is supported on $r$ of the form 
\begin{align}
r\equiv (0,0,...,r_{n-1})\bmod L.  
\end{align}
\end{enumerate}
In other words, $\phi_{S_3}$ satisfies the hypothesis of the Proposition \ref{jacobi-theta-main}. 
\begin{figure}[!htbp]
\begin{align}
% --- LEFT GRID ---
\begin{array}{@{}c@{}}
    \begin{array}{|c|c|c|c|c|c|}
    % BRACE ROW
    \multicolumn{1}{c}{} & \multicolumn{5}{c}{\overbrace{\hspace{6em}}^{ {r\equiv(0,0,\ldots,r_{n-1})  }}} \\ 
    \cline{1-6}
    * & 0 & 0 & 0 & 0 & \cellcolor{black!60} \\ \hline
    0 & * & 0 & 0 & 0 & 0 \\ \hline
    0 & 0 & * & 0 & 0 & 0 \\ \hline
    0 & 0 & 0 & * & 0 & 0 \\ \hline
    0 & 0 & 0 & 0 & * & 0 \\ \hline
    \cellcolor{black!60} & 0 & 0 & 0 & 0 & * \\ \hline
    \end{array}
\end{array}
% --- ARROW ---
\quad \longrightarrow \quad
% --- RIGHT GRID ---
\begin{array}{@{}c@{}}
    \begin{array}{|c|c|c|c|c|c|}
    % DUMMY ROW (For vertical alignment)
    \multicolumn{1}{c}{} & \multicolumn{5}{c}{\overbrace{\hspace{6em}}^{r_0\equiv(0,0,\ldots,0)  }} \\ 
    \cline{1-6}
    * & 0 & 0 & 0 & 0 & 0 \\ \hline
    0 & * & 0 & 0 & 0 & 0 \\ \hline
    0 & 0 & * & 0 & 0 & 0 \\ \hline
    0 & 0 & 0 & * & 0 & 0 \\ \hline
    0 & 0 & 0 & 0 & * & 0 \\ \hline
    0 & 0 & 0 & 0 & 0 & * \\ \hline
    \end{array}
\end{array}
\end{align}
\caption{The effect of Proposition~\ref{jacobi-theta-main}. The entries are to be read $\mod L$, $*$ denotes an unit $\mod L$.} \label{2cell-fig}
\end{figure}

Consequently, we get an $r_0$ of the form $r_0\equiv (0,0,...0)\bmod L$ such that $c_{\phi}(n_0, r_0)\neq 0$ (see Figure \ref{2cell-fig}). Let $T_3=\smat{n_0}{r_0/2}{r_0^t/2}{S_3}$. Then $A\left(F_{T_0,T_1,T_2}, T_3\right)\neq 0$ and $T_3\equiv \mrm{diag}(n_0, s_1, s_2,...s_{n-1})\bmod L $. Further, since $T\equiv T_0\bmod N$, we get that $(n_0, L)=1$. Thus $d(T_3) \equiv  n_0d(S_3)\bmod L$, which gives us that $(d(T_3), L)=1$. 

\textbf{Step 4:} In this step, we again refine the set of Fourier coefficients to obtain the odd and square-free Fourier coefficients to arrive at the conclusion of the Theorem. In this direction, we consider $0\neq F_{T_0,T_1,T_2,T_3}\in \Mrho{n}{N^5}$.  For ease of notation, let us put
\begin{align}
 G:= F_{T_0,T_1, T_2, T_3}.   
\end{align} 
Then clearly $G_{(n-1)}\in M_{\rho_*}(\Gamma_\pm^{(n-1)}(N^5, L))$ is non-zero and satisfies the induction hypothesis. Thus $G$ has the following properties.
\begin{enumerate}
    \item  The Fourier expansion $G$ is supported on $T\equiv \mrm{diag}(t_1, t_2,...t_n)\bmod L$ with $(t_i, L)=1$ for all $1\le i\le n$.
    \item From $G_{(n-1)}$, we get an $S\in \Lambda_{n-1}^+$ with $(d(S), L)=1$ and $d(S)$ odd and square-free such that $A(G_{(n-1)}, S)\neq 0$. In other words, the Fourier-Jacobi coefficient $\varphi_S$ of $G$ is non-zero.
\end{enumerate}
Thus $G$ satisfies the hypothesis of Proposition \ref{alldiag-fund-prop} and using it, we complete the proof of the theorem.
\end{proof}

\begin{remark} \label{deg-freedom}
The conditions on $T$ ($T \equiv \dia(t_1,\ldots, t_n) \mod L$) and $r$ ($r \equiv (0,0,\ldots, r_{n-1}) \mod L$) in Theorem~\ref{diagonal-sqrfree} and Proposition~\ref{jacobi-theta-main}, respectively, might seem a bit unusual, but they are very beneficial for us. The diagonal condition in Theorem~\ref{diagonal-sqrfree} is to ensure ultimately that $(d(T),L)=1$. This diagonal condition behaves well with induction, but along the way, it leads to a roadblock as shown in Figure~2. Proposition~\ref{jacobi-theta-main} is precisely to remove this roadblock.

Next, as briefly explained in the introduction, we need to descend to the scalar-valued Jacobi forms via induction. This is facilitated by Proposition \ref{sca-j}, which currently is valid only for the decomposition of type $(1, n-1)$. Thus, the resulting scalar-valued Jacobi forms satisfy transformation properties with respect to the {\bf degree $1$ group} $\Gamma_1^{(1)}(N)$. As a consequence, when proving the non-vanishing properties of the theta components $h_{\bm\mu}$, with $\bm\mu=(\mu_1, \mu_2,...,\mu_{n-1})\in \mbb Z^{n-1}$, only the last component $\mu_{n-1}$ can be dealt with. More precisely, using the matrices in $\Gamma_1^{(1)}(N,L)$ to arrive at linear equations, we have essentially one degree of freedom to work with. These matrices can be parametrized by integers $a,b \mod L$, but there is a relation among them. However, to handle all the $(n-1)$ entries, one would need $(n-1)L$ degrees of freedom, which is not possible in the degree one situation. It would have been possible if we could have worked with the $(n-1,1)$ type Jacobi forms (we would get $n-1$ tuples ${\bm a,\bm b}$ from $\Sp{n-1}{\z}$ matrices), but in this case, we do not know how to extract a scalar-valued Jacobi form out of it. This approach will be further dealt with in a forthcoming paper. There, we do not require the congruence conditions like in Proposition~\ref{jacobi-theta-main}.

In summary, we can say:

(i) The condition  $r \equiv (0,0,\ldots, r_{n-1}) \mod L$ is necessary, as it is what we are led to consider via induction.

(ii) The condition  $r \equiv (0,0,\ldots, r_{n-1}) \mod L$ is also sufficient, as it gives one variable, viz. $r_0$ to deal with, which matches the one degree of freedom that we get for the degree one groups $\Gamma_1^{(1)}(N,L)$. Let us add that even in the one-variable case, we need $L$ to be square-free to handle it.
\end{remark}
%\begin{remark}
%    We like to remark about the seemingly unusual conditions on $T$ and $r$ in Theorem \ref{diagonal-sqrfree} and Proposition \ref{jacobi-theta-main}, respectively. As briefly explained in the introduction, this is forced upon us by the need to descend to the scalar-valued Jacobi forms via induction. This is facilitated by Proposition \ref{sca-j}, which currently is valid only for the decomposition of type $(1, n-1)$. Thus, the resulting scalar-valued Jacobi forms satisfy transformation properties with respect to the degree $1$ group $\Gamma_1^{(1)}(N)$. As a consequence, when proving the non-vanishing properties of the theta components $h_{\bm\mu}$, with $\bm\mu=(\mu_1, \mu_2,...,\mu_{n-1})\in \mbb Z^{n-1}$, only the last component $\mu_{n-1}$ can be dealt with. More precisely, using the matrices in $\Gamma_1^{(1)}(N)$, we can only get linear equations as in \eqref{lin-eq-1}. Combining this with the necessity of passing through the induction hypothesis,  we are forced to prove the non-vanishing of theta components of the form $h_{0,0,...,\mu_{n-1}}$.
    
 %   this assures a solution to our linear algebra problem...(how)
    
%    This forces the assumption on $r$ in the hypothesis of Proposition \ref{jacobi-theta-main}. This, in turn, forces the condition $T=\mrm{diag}(t_1, t_2,...,t_n)\bmod L$  in Theorem \ref{diagonal-sqrfree}. 
%\end{remark}}

%==================================

%==================================
\section{Applications} \label{applns}
Our results have two primary applications. First, Theorem~\ref{eq-cond} provides a simple criterion to determine if a Siegel modular form possesses a non-zero ``fundamental'' Fourier coefficient at higher levels. Second, Corollary~\ref{spl-val} establishes the unconditional nature of the special values of spinor $L$-functions, removing the conditions previously required in \cite{eischen2024algebraicity}. We also discuss an application to the theta series.

We begin by establishing a preliminary result regarding the non-vanishing of certain Fourier coefficients, analogous to \cite{yamana2009determination}.
\begin{lemma}\label{lem:yamana}
   Let $\chi$ be a Dirichlet character $\mod N$ of conductor $m_\chi$ and $F \in M_\rho(\Gamma_0^{(n)}(N), \chi)$. If $A(F, T)=0$ for all $T\in\Lambda_n$ with $\mf c(T)|(N/m_\chi)$, then $F=0$.
\end{lemma}
\begin{proof}
Let $F(\mc Z)= \sum_{S\in \Lambda_{n-1}}\varphi_T(\tau, z) e(\tr(SZ))$. Suppose $F$ is non-zero, let $S\in \Lambda_{n-1}$ be such that the Fourier-Jacobi coefficient $\varphi_S \neq 0$. Then using similar arguments as in \cite{boche-das}, we get a non-zero component $\phi_S$ of $\varphi_S$, a scalar-valued Jacobi form of weight $k\ge k(\rho)$ for the group $\Gamma_0(N)$ with character $\chi$. We write the Taylor expansion for $\phi_S$ around $z=0$ as below
\begin{align}
    \phi_S(\tau, z)= \sumn_{\lambda\in \mbb N^{n}}\frac{1}{\lambda!}\left.\frac{\partial^\lambda}{\partial z^\lambda}\phi_S(\tau, z)\right |_{z=0} z^\lambda
\end{align} 
and let
$\nu_0=\min \{\nu(\lambda): \left.\frac{\partial^\lambda}{\partial z^\lambda}\phi_S(\tau, z)\right |_{z=0}\neq 0\}$. For some $\lambda$ with $\nu(\lambda)=\nu_0$, let $f(\tau):= \left.\frac{\partial^\lambda}{\partial z^\lambda}\phi_S(\tau, z)\right |_{z=0}$. Then $f\neq 0$ and since $\phi_S\in J_{k,S}(\Gamma_0(N), \chi)$, we see that $f\in M_{k+\nu_0}(\Gamma_0(N), \chi)$. The Fourier coefficient of $f$ is given by \begin{align}
    a_f(n)= \frac{1}{\lambda!}\sumn_{r}A\left(F, \smat{n}{r/2}{r/2}{S}\right) (2\pi i r)^\lambda.
\end{align}

Suppose $A(F, T)=0$ for all $T$ with $\mf c(T)|(N/m_\chi)$. Let $l=(r, \mf c(S))$, then $a_f(n)=0$ for $(n, l)|(N/m_\chi)$. The lemma now follows by \cite[Lemma 3]{yamana2009determination}.
\end{proof}

\subsection{Application to the spaces \texorpdfstring{$M_{\rho}(\Gamma_{0}^{(n)}(N), \chi)$}{MrhoCHi}}
Let $N$ be odd, and  $\chi \bmod N$ be an even Dirichlet character ($\chi(-1)=1$). Let $M_\rho(\Gamma_0^{(n)}(N), \chi)$ denote the space of Siegel modular forms of degree $n$, weight $k$, level $N$ and character $\chi$. That is, the space of holomorphic functions $F:\mbb H_n\longrightarrow \mbb C$ that satisfy $F|_\rho M =\chi(\det(D))F$ for $M=\smat{A}{B}{C}{D}\in \Gamma_0^{(n)}(N)$.

Let $M=\smat{A}{B}{C}{D}\in \Gamma_{\pm}^{(n)}(N, L)$. Then clearly $\det(D)\equiv \pm 1\bmod N$. Thus for any $F\in M_\rho(\Gamma_0^{(n)}(N), \chi)$, we have
\begin{align}\label{chi-inv}
    (F|_\rho M)(Z)&= \chi(\det(D))  F(Z)= F (Z)\text{ if } M\in \Gamma_1^{(n)}(N,L)\\
    (F|_\rho M )(Z)&= \chi(\det(D)) F(Z)= F(Z) \text{ if } M\in \mc E.
\end{align}
Thus for $N$ odd and even characters $\chi\bmod N$, we get that $M_\rho(\Gamma_0^{(n)}(N), \chi)\subset M_\rho(\Gamma_\pm ^{(n)}(N,L)$.

%Note that if $\rho |_{\mc E}$ is irreducible, the condition above is simply $\rho(\pm 1)= \chi(\pm 1)$.

\begin{theorem} \label{eq-cond}
   Let $N$ be odd, and  $\chi \bmod N$ be an even Dirichlet character. The following are equivalent conditions on an $F \in M_{\rho}(\Gamma_{0}^{(n)}(N), \chi)$.
    \begin{enumerate}
        \item 
        $F$ has a non-zero primitive Fourier coefficient.
        \item 
        $F$ has  infinitely many non-zero Fourier coefficients whose discriminants are pair-wise distinct, odd, square-free, and away from $N$ and the quantitative estimates as in Theorem~\ref{mainthm} hold.
        \item 
        $F$ has  infinitely many non-zero Fourier coefficients whose discriminants are pair-wise distinct, odd, square-free.
    \end{enumerate}
      In (2) and (3) above, we can replace `odd, square-free' with `fundamental'.
\end{theorem}
\begin{proof}
    $(1) \implies (2)$: If $F$ has a non-zero primitive Fourier coefficient $A(F,T)$, then we can use the $\glnz$ equivariance to get $A(F,T[U])=\chi(\det(U))\rho(U)A(F,T)$ for all $U \in \glnz$ (see \ref{glnz-eq}). We also know that such a $T'=T[U]=\smat{*}{*}{*}{p}$ for infinitely many primes $p$, and hence for a $p$ with $(p,N)=1$ (see e.g., \cite{pleasants1966representation} or more recently \cite{martin2025prime}). If we note that $M_\rho(\Gamma_0^{(n)}(N), \chi)\subset M_\rho(\Gamma_\pm ^{(n)}(N,L)$, we can then invoke Theorem~\ref{mainthm} with $T'(n,n)$ satisfying its hypothesis, to get $(2)$.

All of the implications $(2) \implies (3)$, $(3) \implies (1)$ are trivial.
\end{proof}

%We have the following corollaries from Theorem \ref{eq-cond}.
\begin{corollary}\label{cor:Nchi}
    Let $N, 
    \chi$ be as above. The equivalent conditions stated in Theorem~\ref{eq-cond} hold true in the following cases.
    \begin{enumerate}
        \item 
        $F \in M_{\rho}(\Gamma_{0}^{(n)}(N), \chi)$ with $\chi \mod N$ primitive.

        \item 
        $F \in M_{\rho}(\Gamma_{0}^{(n)}(N), \chi)$ which is an eigenfunction of  the $U(p)$ operators, for all $p|N$.

        \item A non-zero $F\in S_k^0(\Gamma_0^{(n)}(N), \chi)$, the Atkin-Lehner new-space defined as in \cite{ibukiyama2012atkin}.

         %\item $F \in M_{\rho}(\Gamma_{1, \det}^{(n)}(N))$ and $A(F,T_0) \neq 0$ for some $T_0 \in \Lambda_n$ such that $(T_0(i,j),N)=1$ for some $1\le i, j\le n$.
    \end{enumerate}
\end{corollary}

\begin{proof}
\begin{enumerate}
    \item The proof of $(1)$ follows from Lemma \ref{lem:yamana}. 

    \item  Note that since $F \neq 0$, by  Lemma~\ref{lem:yamana}, we get a $T$ with $\mf c(T)=1$ such that $A(F, M\cdot T) \neq 0$ for some $M \mid N/m_\chi$. But we know that $F|U(m)= \lambda_m F$ for all $m|N^\infty$, so that $A(F, M\cdot T) = \lambda_M A(F,T)$, which shows that $A(F,T) \neq 0$ where $T$ is primitive.
    
    \item Any such $F$ has a non-zero primitive Fourier coefficient from \cite[Theorem 2]{ibukiyama2012atkin} and thus the proof follows. \qedhere
\end{enumerate}
\end{proof}
%%%
\begin{remark}
 When the character $\chi\bmod N$ is odd, the spaces  $M_\rho(\Gamma_0^{(n)}(N), \chi)$ and $M_\rho(\Gamma_\pm ^{(n)}(N,L)$  have trivial intersection, as can be seen from \eqref{chi-inv}. Thus, we cannot use Theorem \ref{mainthm} to conclude anything about the forms in  $M_\rho(\Gamma_0^{(n)}(N), \chi)$.
\end{remark}

\begin{remark} \label{rmk:sydney}
    We note here that in (1) of Corollary~\ref{cor:Nchi}, we do not require $N$ to be square-free. Thus it improves upon the main result in \cite{washburn2025certain}, where other technical conditions on $\chi$ were also present.
\end{remark}

\begin{remark} \label{der}
    We mention here that the primitive property of the Fourier coefficients may not be preserved when we pass from $F$ to $F_{(n-1)}$ etc. This is because the Fourier coefficients of $F_{(n-1)}$ are the derivatives of the Fourier-Jacobi coefficients $\phi_S(\tau,z)$ evaluated at $z=0$; and even if $T= \smat{*}{*}{*}{S}$ was primitive, all the $S$ in the support of the Fourier expansion of $F_{(n-1)}$ may be imprimitive. At least there seems to be no immediate way to rule this out.
\end{remark}

\subsection{Prime discriminants and application to special values of \texorpdfstring{$L$}{l}-functions}

Let $\Gamma^{(n),0}(N)$ be the group defined by
\begin{align}
    \Gamma^{(n),0}(N):= \{\smat{A}{B}{C}{D}\in\mrm{Sp}_n(\mbb Z): B\equiv 0\bmod N \}.
\end{align}
For any integer $k$, the space of modular forms of weight $k$ for the group $\Gamma^{(n),0}(N)$ is defined analogously as in Section \ref{prelim-smf} and we denote it by $M_k(\Gamma^{(n),0}(N))$. For $p|N$, let $U^0(p)$ be the operator on $M_k(\Gamma^{(n),0}(N))$ given by
\begin{align}
    \sumn_{T} A(F,T) e(\frac{\tr(T Z)}{N})\mapsto \sumn_{T} A(F,pT) e(\frac{\tr(T Z)}{N}).
\end{align}
Then we have that
\begin{align}\label{U^0pdef}
    F|U^0(p)= p^{-n(n+1)/2}\sum_{S=S^t, S \bmod p} F|\smat{1_n}{NS}{0_n}{p\cdot 1_n}.
\end{align}

The above definition comes from the Hecke algebra. To see this for the group in question, note that $\Gamma^{(n),0}(N)\cap \smat{1_n}{0_n}{0_n}{p_n}^{-1}\Gamma^{(n),0}(N)\smat{1_n}{0_n}{0_n}{p_n}=\Gamma^{(n),0}(Np)$ and $\Gamma^{(n),0}(N)=\cup_{S=S^t\bmod p} \Gamma^{(n),0}(Np)\smat{1_n}{NS}{0_n}{1_n}$. Thus $\Gamma^{(n),0}(N)\smat{1_n}{0_n}{0_n}{p\cdot 1_n}\Gamma^{(n),0}(N)=\cup_{S=S^t\bmod p}\Gamma^{(n),0}(Np)\smat{1_n}{NS}{0_n}{p\cdot 1_n}$. 
As a consequence, we get 
\begin{align}
p^{-n(n+1)/2}F|\Gamma^{(n),0}(N)\smat{1_n}{0_n}{0_n}{p\cdot 1_n}\Gamma^{(n),0}(N)= p^{-n(n+1)/2}\sum_{S=S^t, S \bmod p} F|\smat{1_n}{NS}{0_n}{p\cdot 1_n}=F|U^0(p). 
\end{align}

Let $\mc M= \smat{N\cdot 1_n}{0_n}{0_n}{1_n}$ and consider $G(Z)=F(N Z)=F|\mc M$. Since $\smat{N\cdot 1_n}{0_n}{0_n}{1_n}^{-1}\Gamma^{(n),0}(N) \smat{N\cdot 1_n}{0_n}{0_n}{1_n}= \Gamma_0^{(n)}(N)$, it is easy to see that $G \in M_k(\Gamma_0^{(n)}(N))$. 

\begin{lemma}\label{U^0-U_0}
    Let $p|N$ and $F\in M_k(\Gamma^{(n),0}(N))$ be an eigenfunction of $U^0(p)$. Then $G=F|\mc M$ is an eigenfunction of $U(p)$.
\end{lemma}
\begin{proof}
   The proof is by computation. Suppose the eigenvalue of $F$ under $U^0(p)$ is $\lambda_p$. We have
    \begin{align}
       G|U(p)&= p^{-n(n+1)/2}\sumn_{S} F|\smat{N\cdot 1_n}{0_n}{0_n}{1_n}|\smat{1_n}{S}{0_n}{p\cdot 1_n}\\
       &=  p^{-n(n+1)/2}\sumn_{S} F|\smat{1_n}{NS}{0_n}{p\cdot 1_n}|\smat{N\cdot 1_n}{0_n}{0_n}{1_n}.
    \end{align}
Thus, we get $G|U(p)=F|U^0(p)|\mc M=\lambda(p) G$, and this completes the proof.  
\end{proof}

\begin{theorem} \label{oddthm}
Let $n$ be odd. Let $F \in S_\rho(\Gamma_{\pm}^{(n)}(N, L))$ be non-zero and $k(\rho) - \frac{n-1}{2} \ge 2 $. Suppose that there exists a $T_0\in \Lambda_n$ with $(T_0(n,n), N)=1$ such that $A(F,T_0)\neq 0$.  Then there exist $T \in \Lambda_n^+$ with $d(T)$ assuming infinitely many odd prime values, such that $A(F, T)\neq 0$. Moreover, the following stronger quantitative result holds: for any given $\epsilon>0$,
\begin{align*}
\# \left( {\mk S}_F (X) \cap \mf P \right) \gg X/\log X,
\end{align*}
where the implied constant depends only on $F$ and $\epsilon$.
\end{theorem}
\begin{proof}
For $F$ as in the statement, from Theorem \ref{mainthm}, we get a $T_1$ with $d(T_1)$ odd, square-free and $(d(T_1), N)=1$ such that $A(F, T_1)\neq 0$. Now consider $F_{T_1}$ from Lemma \ref{lem:FTN}. Again using Theorem \ref{mainthm},  we get a non-zero scalar-valued Jacobi form $\varphi_{S_0}$ of weight $k'\ge k(\rho)$ from the Fourier-Jacobi expansion of $F_{T_1}$ with $d(S_0)$ odd, square-free and $(d(S_0), N)=1$. Now, from Proposition \ref{prop:thetacom}, we get a non-zero theta component $h_\mu$ of $\varphi_{S_0}$, where $\mu$ is either $N^2S_0$ or $S_0$-primitive. As in the proof of Proposition \ref{alldiag-fund-prop}, define $H_\mu(\tau)= h_\mu((d_{S_0})\tau)$. Then $H_\mu \in S_{k'-\frac{n-1}{2}}(\Gamma_1^{(1)}(d_{S_0}^2 N^2))$. Further,  note that $H_\mu$ has all non-zero Fourier coefficients away from $d_{S_0}N$ (see Section \ref{sec:thmdiag} for details).

The rest of the proof now follows as in \cite[Theorem~5.1]{boche-das} and will not be repeated here.
\end{proof}
\begin{proposition}\label{prop:FU^0-prime}
    Let $n$ be odd and  $F\in M_k(\Gamma^{(n),0}(N)) $ be an eigenfunction of $U^0(p)$ for all $p|N$. Then there exist infinitely many $T$ with odd and prime discriminant such that $A(F, T)\neq 0$.
\end{proposition}
\begin{proof}
Let $G=F|\mc M$ be as above. Then, the Fourier expansion of $G$ is given by
    \begin{align}
        G(Z)= \sumn_{T\in \Lambda_n^+}A(F,T)e(\tr(TZ)).
    \end{align}

Next, since $F$ is an eigenfunction $U^0(p)$, from Lemma \ref{U^0-U_0} we get that $G$ is also an eigenfunction of $U(p)$ for all $p|N$. Now Corollary \ref{cor:Nchi} shows that $G$ has a non-zero primitive Fourier coefficient, and thus using the $\glnz$-equivalence of Fourier coefficients of $G$, there exists a $T_0$ with $(T_0(n,n), N)=1$ such that $A(G, T_0)\neq 0$.

Now invoking Theorem $\ref{oddthm}$, we see that there exist $T \in \Lambda_n^+$ with $d(T)$ assuming infinitely many odd prime values, such that $A(G, T)\neq 0$. Since $A(G, T)=A(F,T)$, we see that the same conclusion holds for $F$. This completes the proof of the proposition.
\end{proof}

\begin{corollary} \label{spl-val}
    The statement of the special-value result stated in \cite[Theorem~1.1.2]{eischen2024algebraicity} is unconditional, i.e., \cite[Condition~1.1.1]{eischen2024algebraicity} is true.
\end{corollary}

\begin{proof}
    In \cite{eischen2024algebraicity}, the authors consider a newform for the group $\Gamma^{(3),0}(N)$. Such an $F$ has a Fourier expansion of the form
    \begin{align} \label{f-fe}
        \sumn_{T \in \Lambda_3^+} A(F,T) e(\frac{T Z}{N}).
    \end{align}
It follows from the discussion in \cite[Section~5.1]{boche-das} and \cite[Condition~1.1.1]{eischen2024algebraicity} that the special-value result in \cite{eischen2024algebraicity} holds whenever $F$ has a non-zero Fourier coefficient $T$ with odd, prime discriminant. Indeed, in this case, via the Clifford correspondence (see e.g., \cite[Proposition~3.3]{pollack2017spin}), $T$ corresponds to a maximal order in a quaternion algebra over $\Q$.
Since such an $F$ is an eigenfunction of $U^0(p)$ for all $p|N$, the existence of such a $T$ follows from Proposition \ref{prop:FU^0-prime}, and this gives us the corollary.
\end{proof}

\begin{remark}
    It is quite plausible that $F$ and $G$ in the above proof have the same Satake parameters and hence the same spinor $L$-function. Then the algebraicity result should also hold verbatim for $G$.
\end{remark}

\subsection{Representation numbers, theta series}

One frequently used application of representation numbers to modular forms is to show the existence of an odd prime $p$ such that $A(F,T)$ has the property that $T(n,n)=p$.
For the group $\Gamma^{(n)}_0(N)$, this can be proved much more elementarily, see \cite{martin2025prime} or can be deduced at once from \cite{pleasants1966representation} or \cite{iwaniec1974primes}, as the entire unimodular group can be embedded in it.  In \cite{martin2025prime}, it was shown more generally that a primitive half-integral matrix $T \in \Lambda_n^+$ represents another primitive matrix $S$ of lower dimension $r$. We can show the following generalization, which may have some implication in the theory of quadratic forms. We assume that $r \ge 2$, as the case $r=1$ is covered in the references mentioned above: when $r=1$, the convention is that $T$ represents a scalar primitive matrix if it represents a prime.

Suppose $T_0 \in \Lambda_n^+$ and $2 \le r \le n$. Define the theta series $\theta_{T_0}$ associated to $T_0$ by
\begin{align}
    \theta_{T_0}(Z) = \sum_{X \in M_{n,r}(\z)} e(\tr T_0[X] \cdot Z) = \sum_{S \in \Lambda_r} r(T_0;S) e(\tr S Z) \q (Z \in \h_r),
\end{align}
where $\displaystyle r(T_0;S)=\{ X \in M_{n,r}(\z) \mid T_0[X]=S\}$
denotes the number of integral representations of $S$ by $T_0$. Further put
\begin{align}
N&:=\mrm{level}(T_0) := \min \{M \ge 1 \mid MT_0^{-1} \in \Lambda_n^+ \} \q \text{ and }\\
{\mc S}_{T_0} (X) &:= \{ d \leq X, \,  d \, \mbox{\rm{odd, square-free}}, \, (d,N)=1 \, \mid  d(S)=d \, \mbox{\rm{ for some }} S\in \Lambda_r^+ \mbox{\rm{ and }}  r(T_0,S) \neq 0 \}.  
\end{align}

\begin{theorem}\label{theta}
    With the above notation, let  $T_0$ and $n$ be such that $\mrm{level}(T_0)$ is odd,    $n \ge r+3$, $r \ge 2$. If, $T_0$ is primitive and $4|n$, one has for any $\epsilon>0$,
    $\displaystyle {\mc S}_{T_0} (X) \gg_{T_0, \epsilon} \begin{cases} X^{1-\epsilon} &\text{ if } r \text{ is odd };  \\ 
X^{5/8 - \epsilon} &\text{ if } r \text{ is even }.
\end{cases}$
\end{theorem}

\begin{proof}
    The proof follows by combining Theorem \ref{mainthm} and Theorem ~\ref{eq-cond}. First, it is well-known that $\theta_{T_0} \in M_{n/2}(\Gamma_0^{(r)}(N), \chi_Q)$, where $Q$ is the quadratic form defined by $T_0$, and $\chi_Q$ being the quadratic character attached to $Q$. We know that $\chi_Q(\cdot) = \left( \frac{ \mrm{disc}(2T_0) }{\cdot} \right)= \left( \frac{ (-1)^{n/2} d(T_0) }{\cdot} \right)$, which shows that $\chi_Q$ is even.
    
    Second, from \cite{martin2025prime}, we get that $A(T_0;S_0) \neq 0$ for some primitive $S_0 \in \Lambda_n^+$. Then the theorem follows from part $(1)$ (and hence part $(2)$) of Theorem~\ref{eq-cond}. Since ${\mc S}_{T_0} (X) = {\mk S}_{\theta_{T_0}} (X)$, the theorem is proved.
\end{proof}
Note that for any $S$ represented by $T_0$ as in the Theorem, the vector $X$ representing it is necessarily `primitive' in the sense that it can be completed to an unimodular matrix. See also \cite{martin2005dimensions} on this, for example. We remark that Theorem~\ref{theta} may not be easy to obtain without the use of modular forms.

\printbibliography
\end{document}